\makeatletter
\newcommand{\dontusepackage}[2][]{%
  \@namedef{ver@#2.sty}{9999/12/31}%
  \@namedef{opt@#2.sty}{#1}}
\makeatother
\dontusepackage{subfigure}

\documentclass[10pt,onecolumn]{IEEEtran}
\usepackage{lmodern}
\usepackage{amssymb,amsmath}
\usepackage{ifxetex,ifluatex}
\usepackage[usenames,dvipsnames]{color}
\usepackage{fixltx2e} 
\ifnum 0\ifxetex 1\fi\ifluatex 1\fi=0 
  \usepackage[T1]{fontenc}
  \usepackage[utf8]{inputenc}
  \usepackage{eurosym}
\else 
  \ifxetex
    \usepackage{mathspec}
    \usepackage{xltxtra,xunicode}
  \else
    \usepackage{fontspec}
  \fi
  \defaultfontfeatures{Mapping=tex-text,Scale=MatchLowercase}
  
\fi
\IfFileExists{upquote.sty}{\usepackage{upquote}}{}
\IfFileExists{microtype.sty}{%
\usepackage{microtype}
\UseMicrotypeSet[protrusion]{basicmath} 
}{}
\usepackage[margin=1.0in,bottom=1.5in]{geometry}
\usepackage{listings}
\usepackage{graphicx}
\makeatletter
\def\maxwidth{\ifdim\Gin@nat@width>\linewidth\linewidth\else\Gin@nat@width\fi}
\def\maxheight{\ifdim\Gin@nat@height>\textheight\textheight\else\Gin@nat@height\fi}
\makeatother
\setkeys{Gin}{width=\maxwidth,height=\maxheight,keepaspectratio}
\usepackage{caption}
\usepackage{float}
\usepackage{subfig}
\usepackage{algorithm} 
\let\scholmdAlgorithm\algorithm
\let\endscholmdAlgorithm\endalgorithm
\let\algorithm\relax \let\endalgorithm\relax

{
 \catcode`\*=11\relax
 \global\let\scholmdAlgorithm*\algorithm*
 \global\let\endscholmdAlgorithm*\endalgorithm*
 \global\let\algorithm*\relax 
 \global\let\endalgorithm*\relax
}
\ifxetex
  \usepackage[setpagesize=false, 
              unicode=false, 
              xetex]{hyperref}
\else
  \usepackage[unicode=true]{hyperref}
\fi
\hypersetup{breaklinks=true,
            bookmarks=true,
            pdfauthor={},
            pdftitle={Total-variation regularization strategies in full-waveform inversion},
            colorlinks=true,
            citecolor=black,
            urlcolor=blue,
            linkcolor=black,
            pdfborder={0 0 0}}
\usepackage[normalem]{ulem}

\DeclareMathOperator*{\argmin}{arg\,min}

\providecommand{\I}{}
\renewcommand{\I}{\mathrm{I}}

\providecommand{\L}{}
\renewcommand{\L}{\mathcal{L}}

\providecommand{\bbm}{}
\renewcommand{\bbm}{\begin{bmatrix}}
\providecommand{\ebm}{}
\renewcommand{\ebm}{\end{bmatrix}}
\providecommand{\dm}{}
\renewcommand{\dm}{\Delta m}

\providecommand{\comments}{}
\renewcommand{\comments}[1]{}
\providecommand{\e}{}
\renewcommand{\e}[1]{\ensuremath{\times 10^{#1}}}

\title{Total-variation regularization strategies in full-waveform inversion}
\author{Ernie Esser$^{\dagger,1,2}$,
Lluis Guasch$^2$, 
Tristan van Leeuwen$^3$,
Aleksandr Y. Aravkin$^4$,
and Felix J. Herrmann$^1$\\ 
$^1$The University of British Columbia, Department of Earth, Ocean and Atmospheric Sciences;\\
$^2$Sub Salt Solutions Limited;\\
$^3$Mathematical Institute, Utrecht University;\\
$^4$Applied Mathematics, University of Washington;\\
}
\date{}

\begin{document} 
\maketitle 

\begin{abstract} 

We propose an extended full-waveform inversion formulation that includes general convex constraints on the model. Though the full problem is highly nonconvex, the overarching optimization scheme arrives at geologically plausible results by solving a sequence of relaxed and warm-started constrained convex subproblems. The combination of box, total-variation, and successively relaxed asymmetric total-variation constraints allows us to steer free from parasitic local minima while keeping the estimated physical parameters laterally continuous and in a physically realistic range. For accurate starting models, numerical experiments carried out on the challenging 2004 BP velocity benchmark demonstrate that bound and total-variation constraints improve the inversion result significantly by removing inversion artifacts, related to source encoding, and by clearly improved delineation of top, bottom, and flanks of a high-velocity high-contrast salt inclusion. The experiments also show that for poor starting models these two constraints by themselves are insufficient to detect the bottom of high-velocity inclusions such as salt. Inclusion of the one-sided asymmetric total-variation constraint overcomes this issue by discouraging velocity lows to buildup during the early stages of the inversion. To the author's knowledge the presented algorithm is the first to successfully remove the imprint of local minima caused by poor starting models and band-width limited finite aperture data.

%
%

\footnotetext[2]{John "Ernie" Esser passed away on March 8, 2015 while preparing this manuscript. The original is posted here: \\ {\scriptsize\href{https://www.slim.eos.ubc.ca/content/total-variation-regularization-strategies-full-waveform-inversion-improving-robustness-noise}{https://www.slim.eos.ubc.ca/content/total-variation-regularization-strategies-full-waveform-inversion-improving-robustness-noise}}.}

\end{abstract}
 
\section{Introduction}\label{introduction}

Full-waveform inversion reconstructs high-resolution gridded models of
subsurface medium parameters from seismic measurements by solving a
large-scale inverse problem \cite{tarantola82,Tarantola1984}. The classic approach uses the following non-linear least-squares formulation:
\begin{equation} \label{eq:NLLS} \min_{m}
\frac{1}{2}\sum_{j=1}^{N_s}\|F(m)q_j - d_j\|_2^2, 
\end{equation} 
where the vector $m \in \mathbb{R}^M$ contains a discrete set of parameters describing the medium (e.g., spatially varying soundspeed), the vector $q_j$ is the source term for the $j^{\mathrm{th}}$ experiment, $F(m)$ is the non-linear forward operator and $d_j \in \mathbb{C}^{N_r}$ are the corresponding measurements \cite{tarantola82}. We denote the total number of measurements by $N = N_r \times N_s$. In the large-scale setting, problem size can vary from $M \sim 10^6$, $N \sim 10^4$ up to $M \sim 10^9$, $N\sim 10^6$.

Application of the forward operator typically involves solving a partial differential equation (PDE) with coefficients $m$ and right-hand-sides $q_j,\, j=1\cdots N_s$. The operator $F(m)$ can be formally written $PA(m)^{-1}$, and in practice requires the solution of a PDE cast into a discretized linear system of equations $A(m)u = q$ (e.g. Helmholz), followed by applying the measurement operator $P$ to the solution $u\in \mathbb{R}^M$. Solving these PDEs is the main computational cost in evaluating the objective~\eqref{eq:NLLS}. Calculating the gradient of the objective can be done using the so-called \textit{adjoint-state} approach and requires additional PDE solves for each experiment.

A well-known problem with such approaches is that the objective may have \textit{parasitic} stationary points or local minima that are not informative about the true parameters. In seismic inversion in particular, the so-called loop- or cycle-skipping phenomenon is a common source of such local minima \cite{Symes2009}. In practice this means that the inversion is very sensitive to the initial value of $m$. To mitigate this problem, many alternative formulations to the standard least-squares problem~\eqref{eq:NLLS} have been proposed \cite{Symes2008,Biondi2013,warner2014adaptive}, extending the parameter space to avoid local minima. The general strategy 
is to design formulations that are useful in controlling the model space, moving away from these detrimental parasitic stationary points and toward models with realistic physical properties. The main challenge is to make such formulations {\it computationally efficient}, so that they can be used on large-scale problems.

\smallskip \noindent {\bf Contributions.} Our primary focus is to design robust formulations by controlling the model space. We propose a general framework for PDE-constrained optimization that allows multiple convex constraints to be applied. Though the overarching problem is non-convex, it is solved using a sequence of constrained convex subproblems. The generality of the scheme allows the simultaneous application of several constraints that turn out to be particularly well suited to seismic inversion: bound constraints on the model parameters (slowness squared), discontinuity preserving total-variation constraints, and a novel asymmetric variant developed specifically for seismic inverse problems. The best results obtained use all three constraints.

When given accurate starting models, the final scheme removes inversion artifacts and improves the delineation of high-velocity and high-contrast salt inclusions. More importantly, the asymmetric total-variation constraint prevents the buildup of detrimental artifacts related to parasitic local minima for poor starting models. To our knowledge, this scheme is a first instance of a {\it hands-free} inversion methodology that produces high-fidelity reproducible results where other wave-equation based inversions fail.

In the remainder of the introduction, we give a brief illustration of two key techniques that help control the model space: a penalized reformulation of~\eqref{eq:NLLS}, and the systematic incorporation of prior knowledge via multiple convex constraints.

\smallskip \noindent {\bf Avoiding local minima.} An important penalized reformulation of PDE-constrained optimization that helps avoid local minima was recently proposed by \cite{VanLeeuwen2013}:
\begin{equation}
\label{eq:ENLLS}
\min_{m,\Delta q} \frac{1}{2}\sum_{j=1}^{N_s}\|F(m)(q_j + \Delta q_j) - d_j\|_2^2 + \lambda^2 \|\Delta q_j\|_2^2,
\end{equation}
where $\Delta q = [\Delta q_1; \Delta q_2, \ldots, \Delta q_{N_s}]$ can be thought of as slack variables that allow some freedom in fitting the data even for a wrong set of parameters. As $\lambda \uparrow \infty$, $\|\Delta q\|_2\downarrow 0$, and~\eqref{eq:ENLLS} coincides with the original problem~\eqref{eq:NLLS} \cite{vanLeeuwen2013Penalty2}. The above problem can be solved by \textit{projecting} out the slack variable~\cite{Aravkin2012}, which can be done efficiently since $\Delta q$ has a closed form solution for every fixed $m$. The result is a modified objective depending on $m$ alone, and well-behaved as $\lambda$ increases~\cite{aravkin2016qp}. To illustrate the benefit of this extension, consider the following example.

\begin{figure}
\centering
\subfloat[\label{fig:example1a}]{\includegraphics[width=0.330\hsize]{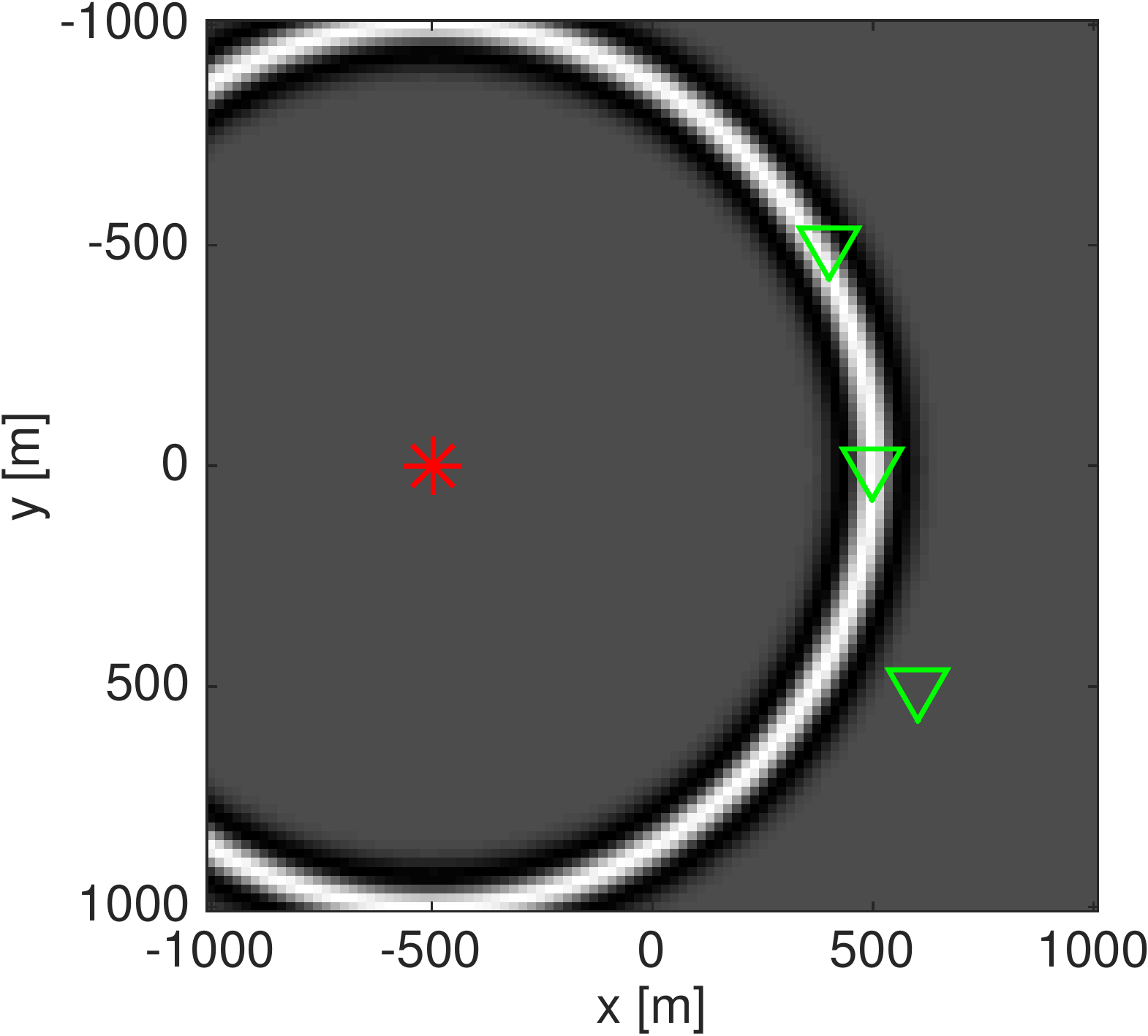}}
\subfloat[\label{fig:example1b}]{\includegraphics[width=0.330\hsize]{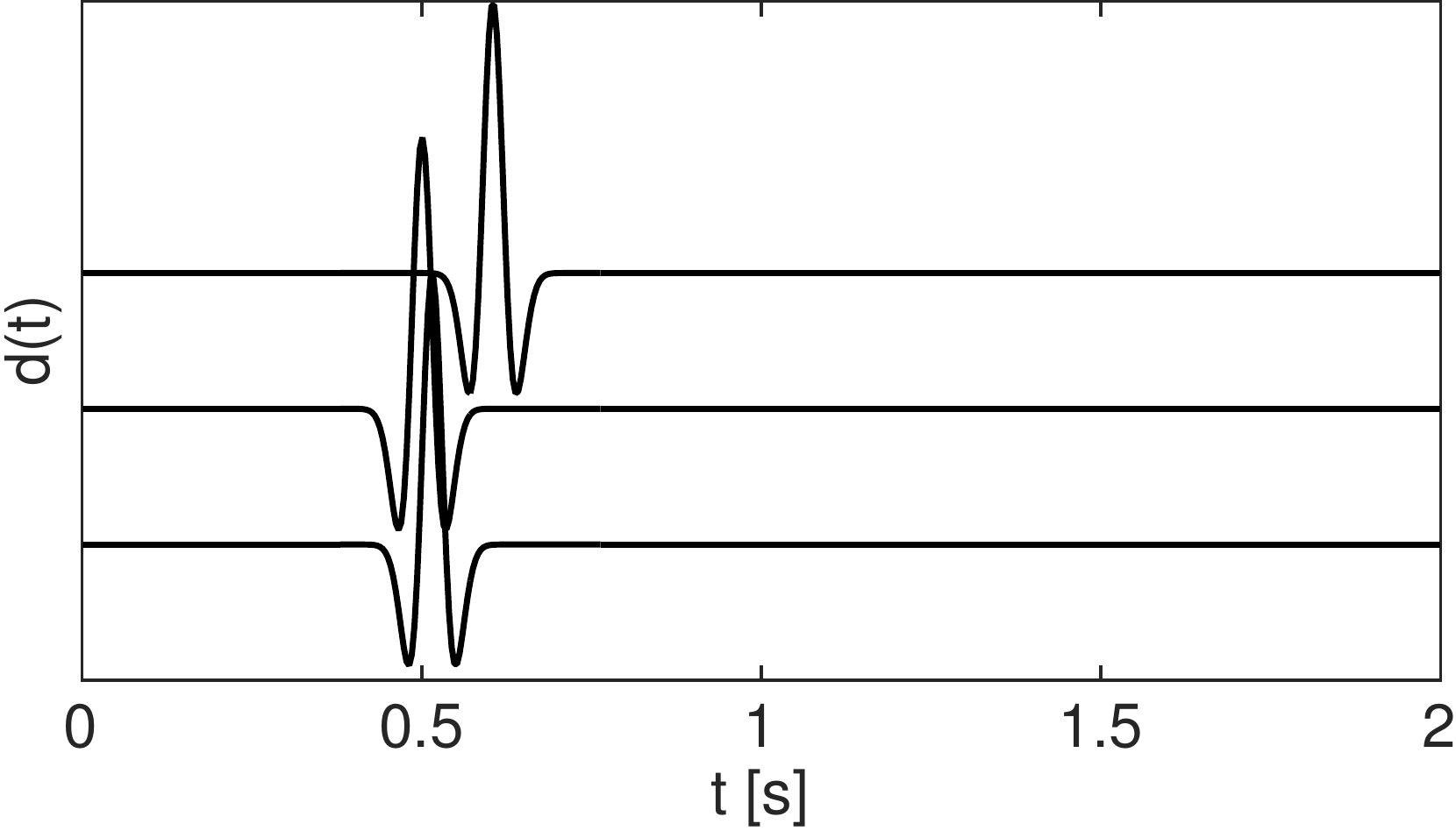}}
\subfloat[\label{fig:example1c}]{\includegraphics[width=0.330\hsize]{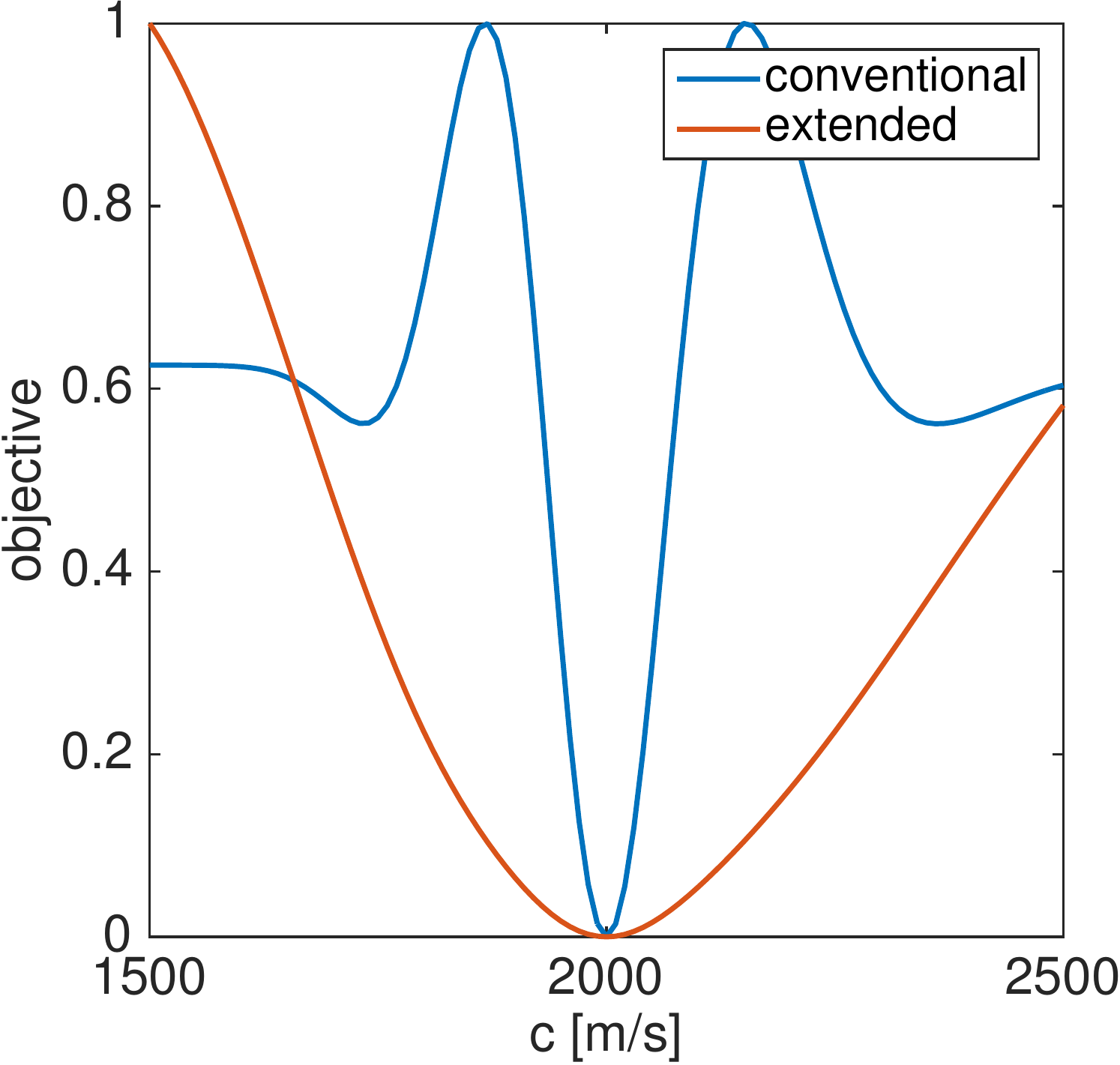}}
\caption{Example 1. (a) Snapshot of time-domain wavefield generated by a single
source ($*$); (b) recorded time-domain data at three receiver locations
($\nabla$); (c) objective function for conventional
($\lambda\rightarrow\infty$) and extended ($\lambda=0.1$)
formulations.}\label{fig:example1} \end{figure}

We aim to retrieve the soundspeed of a medium from measurements of the response
of a bandlimited source. The experimental setup is depicted in figure
\ref{fig:example1a}; a single source ($*$) and three receivers ($\nabla$) are
located as shown. The source emits a pulse that is subsequently recorded,
leading to three time series as depicted in figure \ref{fig:example1b}. The
forward operator, acting on the temporal source signature $q(t)$, is defined as
\[
F(c)q = 
\left(
\begin{matrix}
	\mathcal{F}^{-1}e^{\imath\omega\|x_1 - x_s\|_2/c}\mathcal{F}q\\
	\mathcal{F}^{-1}e^{\imath\omega\|x_2 - x_s\|_2/c}\mathcal{F}q\\
	\mathcal{F}^{-1}e^{\imath\omega\|x_3 - x_s\|_2/c}\mathcal{F}q\\
\end{matrix}
\right),
\]
where $\mathcal{F}$ denotes the temporal Fourier transform, $\omega$ the angular frequency, $\imath=\sqrt{-1}$, $x_i$ is the receiver location, $x_s$ is the source location, and $c$ is the soundspeed of the medium.

Figure \ref{fig:example1c} shows the conventional ($\lambda\rightarrow\infty$) and extended ($\lambda=0.1$) objectives as a function of the velocity $c$. The local minima are clearly visible in the conventional objective. By enlarging the search space, the extensions help to mitigate these local minima effectively.

In addition to exhibiting local minima, the above toy problem is severely ill-conditioned due to the band-limited nature of the data. Large- and small-scale variations of $c$ are difficult to retrieve because low and high frequencies are missing. 
Full-waveform inversion is also hampered by missing spatial frequencies at both ends of the spectrum due to physical 
constraints on seismic surveys, which include restriction of sources and receivers to the surface, bandwidth limitation of sources and receivers, and finite aperture (maximal distance between sources and receivers).  

Since our main objective is to retrieve the global velocity structure in geologic areas with high-velocity and high-contrast (sharp) inclusions, we need to recover both the low frequencies (responsible for the kinematics and therefore placement of velocity perturbations), and the high frequencies to ensure accurate delineations of the high-velocity inclusions. A complicating factor is that all scales in $c$ are intrinsically coupled through the wave equation--i.e., they appear as coefficients nonlinearly in the PDE, so failure to capture the large-scale variations in turn leads to a failure to retrieve useful information on the medium-scale features. As in many other inverse problems, invoking prior knowledge in the form of certain rudimentary constraints can have a significant impact on the quality of recovered models. 

\begin{figure}
\centering
\subfloat[\label{fig:example2a}]{\includegraphics[width=0.330\hsize]{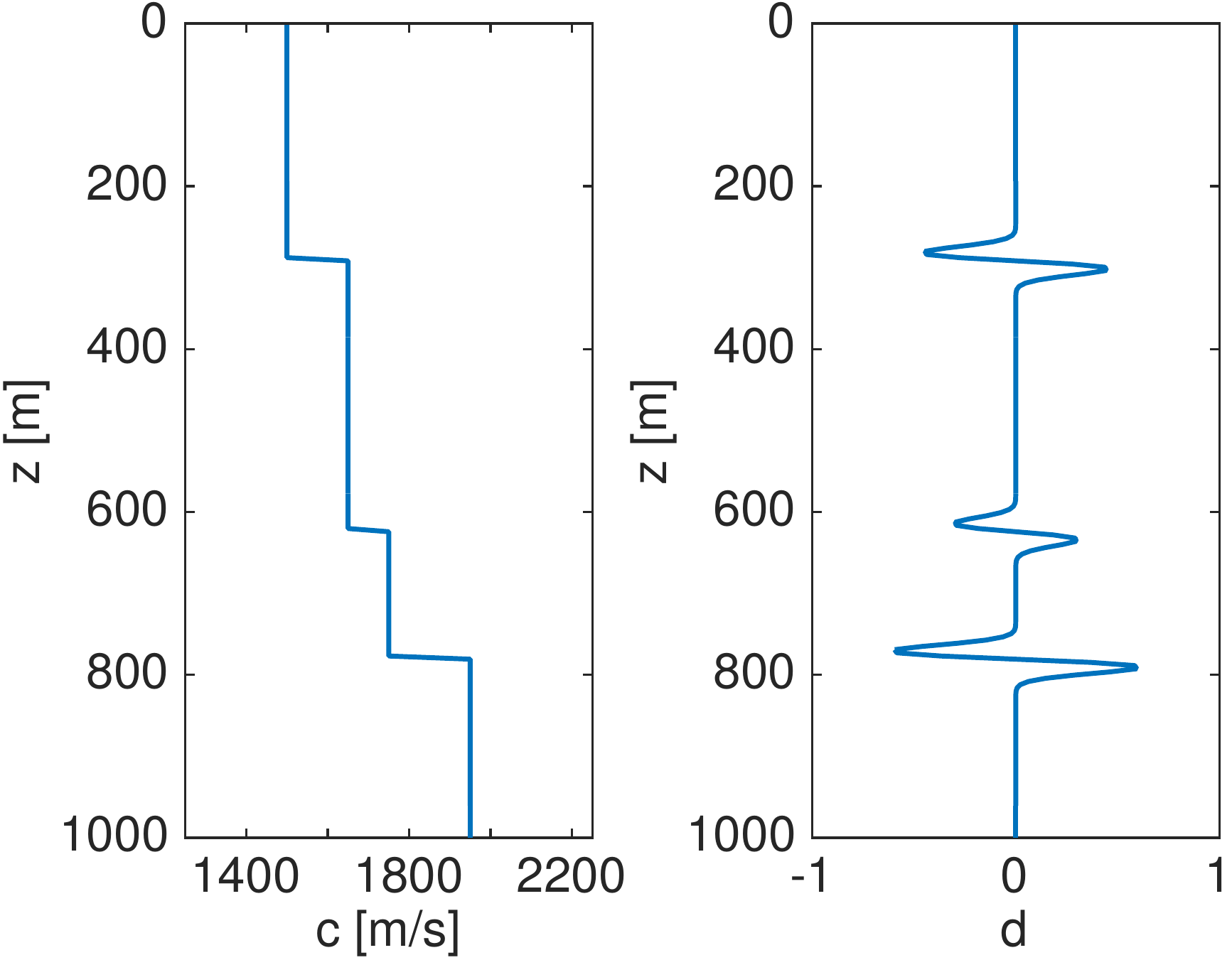}}
\subfloat[\label{fig:example2b}]{\includegraphics[width=0.330\hsize]{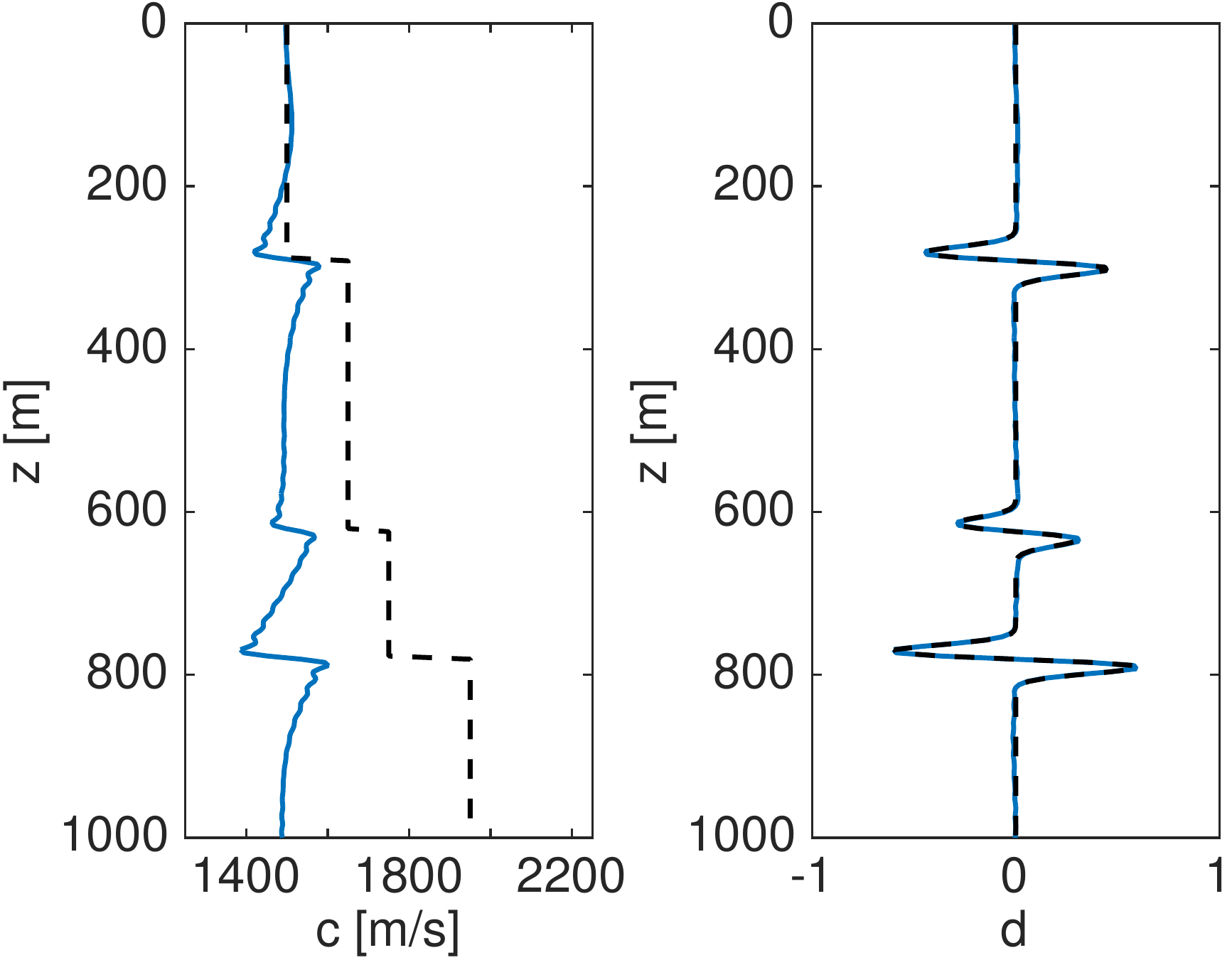}}
\subfloat[\label{fig:example2c}]{\includegraphics[width=0.330\hsize]{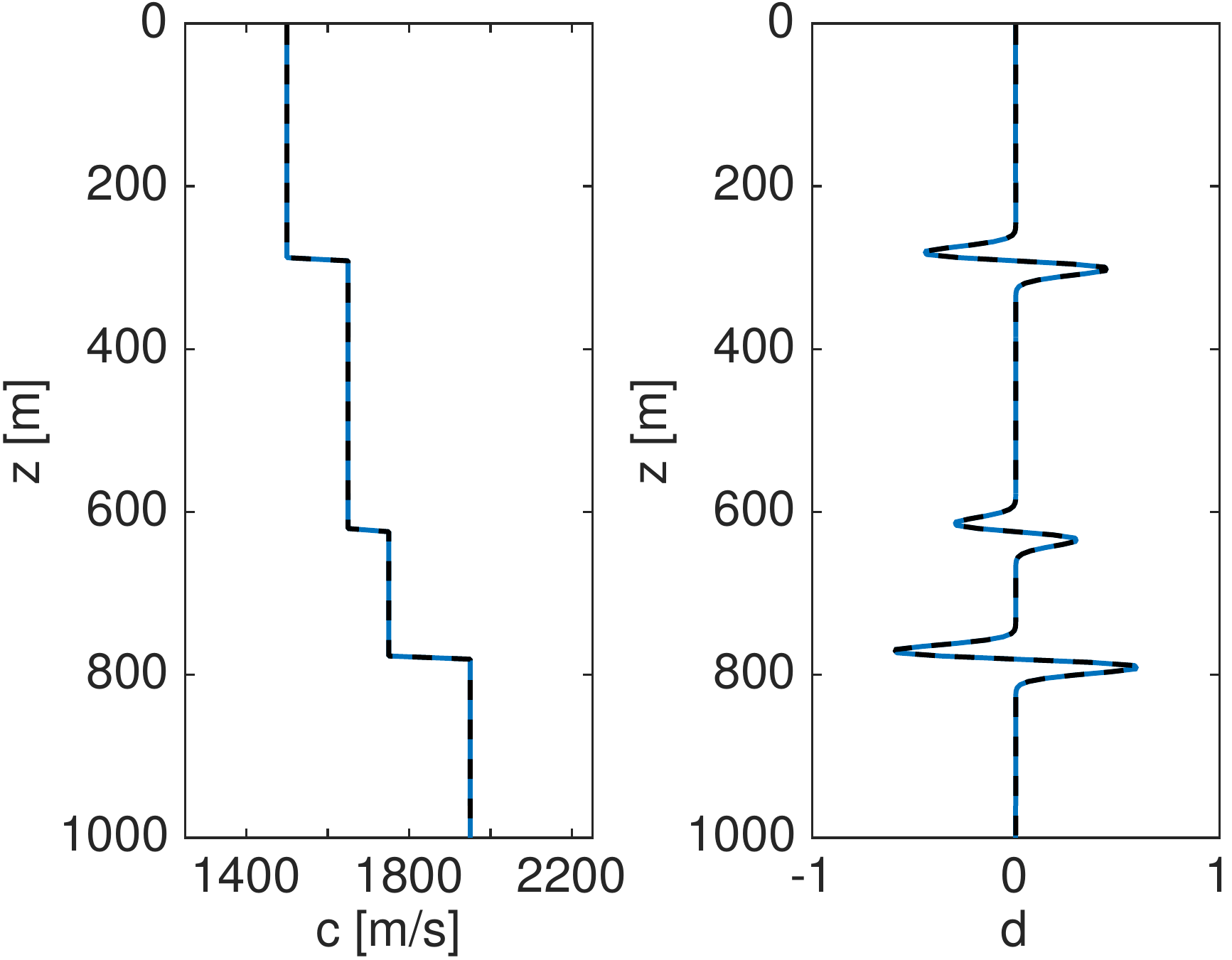}}
\caption{Example 2. (a) Monotonously increasing velocity profile and observed data. (b) Unconstrained inversion result fits the observations but misses the global trend. (c) Constrained inversion results fits the data and recovers the velocity model accurately.}\label{fig:example2}
\end{figure}

\smallskip\noindent{\bf Incorporating prior knowledge using constraints.} 
We aim to retrieve a one-dimensional function $c(z)$ from band-pass filtered measurements $d = F(c)q \equiv c \ast q$, where $\ast$ denotes convolution and $q$ is the band-pass filter in the time domain. Figure \ref{fig:example2a} shows the ground-truth discrete velocities and corresponding data. A least-squares reconstruction is shown in \ref{fig:example2b}. This result clearly illustrates the failure to retrieve the large-scale trends in $c$. Since $c$ represents a physical parameter, failure to reproduce these large-scale trends changes the physical interpretation (e.g. the two-way travel time) completely. We can incorporate 
prior information through constraints by imposing the lower bound $c_0 = 1500\, \mathrm{ms^{-1}}$, and requiring that $c$ increase monotonically. The resulting formulation is 
\begin{equation}
\label{eq:ConEx}
\min_{c\geq c_0} \|F(c)q - d\|_2^2 \quad \mbox{subject to}\quad Dc\geq 0,
\end{equation}
where $D$ is the finite difference matrix. The results from this constrained formulation are shown in figure \ref{fig:example2c}. The additional constraints allow us to retrieve the ground-truth exactly. 
This example is simplistic, but effectively illustrates the important role constraints can play in restricting 
the feasible model space, with potential to recover broadband velocity profiles from bandwidth limited data.

Motivated by these ideas, we develop an inversion framework 
able to avoid local minima and mitigate ill-posedness by incorporating 
three kinds of constraints: (1) bounds, (2) total variation (TV), and (3) a novel one-sided TV extensions based on the hinge loss. Both the classic formulation~\eqref{eq:NLLS} as well as the relaxed formulation~\eqref{eq:ENLLS} fit naturally into this framework. Formulation~\eqref{eq:ENLLS} has additional modeling advantages (see Fig.~\ref{fig:example1}) as well as computational advantages highlighted in Section~\ref{including-convex-constraints}.

The seismic setting comprises imaging sedimentary basins with the following features: 
(1) depth-increasing velocities, (2) high-velocity and high-contrast inclusions (i.e. salt bodies), and
(3) potential over-pressured reservoirs that yield velocity lows. We show that the velocity structure over this complicated geology can be recovered using the proposed framework. The 
one-sided TV constraint is especially important, since it penalizes the negative vertical derivative more heavily than the positive vertical derivative making it less likely for the velocity to decrease with depth. 
The final workflow gradually relaxes this asymmetric constraint, recovering velocity structures with complex high-velocity and high-contrast inclusions from poor starting models.

%

\smallskip
\noindent
{\bf Outline of the paper.} The paper is organized as follows. In Section~\ref{including-convex-constraints}, we specify a scaled gradient projection framework. In Section~\ref{total-variation-regularization}, we show how to incorporate bound constraints and the TV constraint, and how to solve the resulting convex subproblems. Numerical experiments for bound and TV constraints are presented in Section~\ref{numerical-experiments}. A new asymmetric TV variant, implementation details, and numerical results for joint convex constraints are presented in Section~\ref{one-sided-tv-constraint}. Sections~\ref{final-discussion} and \ref{conclusions-and-future-work} present discussion and conclusions. All necessary details needed to reproduce the numerical examples are given in the appendix.

\section{Scaled Gradient Projection for Large-Scale
Programs}\label{including-convex-constraints}

Considering formulations~\eqref{eq:NLLS} and~\eqref{eq:ENLLS} with constraints such as~\eqref{eq:ConEx}, 
our goal is to develop a framework for large-scale optimization where we distinguish between relatively simple convex constraints and more computationally expensive PDE constraints, discussed later. 

First consider problems of the following form:
\begin{equation}
\min_m f(m) \ \text{ subject to }\ m \in C \ ,
\label{Fmincon}
\end{equation}
where $f(m)$ is smooth, but may be nonconvex, and expensive to evaluate. The set $C$ represents the convex constraint, designed to mitigate the ill-posedness of the inverse problem. Box constraints are common to seismic inversion, with $C = \{m : m_i \in [b_i , B_i]\}$ for $i=1, \cdots, M$, where the $b_i$'s are the lower and the $B_i$ are the upper bounds for each entry in $m$. The set $C$ may also be constructed as the intersection of several simpler sets. In the following sections, 

The scaled projected gradient method for~\eqref{Fmincon}
is well-adapted for the PDE constrained setting, 
and can incorporate different constraint sets $C$ (e.g. intersection of box and TV constraints). 
We first consider scaled gradient descent, and then describe scaled projected gradient. 

\subsection*{Scaled gradient descent}

Following \cite{Bertsekas1999}, we start by considering the following unconstrained problem:
\[
\min_m f(m),
\]
and apply an iterative algorithm of the form 
\begin{equation}
\begin{aligned}
\dm & = \argmin_{\dm} \dm^\top \nabla f(m^n) + \frac{1}{2} \dm^\top H^n \dm \\
m^{n+1} & = m^n + \dm.
\end{aligned}
\label{SGD}
\end{equation}
In this quadratic minimization problem, $\nabla f(m^n)$ is the gradient of the objective with respect to the discretized model at the $n^{\text{th}}$ iteration and $H^n$ is a positive definite approximation to the Hessian of $f$. The Hessian approximation can range from a simple gradient descent with a step size of magnitude $\alpha$ if $H^n=\frac{1}{\alpha}I$, to more sophisticated forms such as $H^n = \nabla^2f(m^n)$, which corresponds to Newton's method using the full Hessian (when it is positive definite). Among positive definite Hessian approximations, the choice is typically shaped by a tradeoff between quality of approximation and efficiency of the scaled gradient iteration. Problem structure plays a key role here. For example, the Hessian for the conventional formulation~\eqref{eq:NLLS} is a dense matrix whose evaluation requires PDE solves, while the extended formulation~\eqref{eq:ENLLS} produces accurate sparse Hessian approximations.

\subsection*{Scaled gradient projections}
Adding constraints $m\in C$ to solve~\eqref{Fmincon} requires a simple modification to the model updates $\dm$: 
\begin{equation}
\begin{aligned}
\dm & = \argmin_{\dm} \dm^\top \nabla f(m^n) + \frac{1}{2} \dm^\top H^n \dm \\
& \text{ s.t. }\ m^n + \dm \in C. \\
m^{n+1} & = m^n + \dm.
\end{aligned}
\label{dmcon}
\end{equation}
This iteration is known as scaled gradient {\it projection} method \cite{Bertsekas1999, Bonettini2009}. 
It ensures model iterates $m^{n+1}$ are always in $C$, and admits provable convergence guarantees,
discussed below. Note that   
solving for the unconstrained minimizer $\dm$ in~\eqref{SGD} and then projecting onto $C$
\begin{equation}
\label{eq:projGrad}
m^{n+1} = \Pi_C(m^n - (H^n)^{-1}\nabla(f(m^n)),
\end{equation}
will not converge to a solution of~\eqref{Fmincon} for general $H^n$~\cite{Bertsekas1999}. 

We show that iterations~\eqref{dmcon} remain computationally tractable as
long as $C$ is easy to project onto, or is an intersection several simple simple convex constraints. 
In PDE constrained optimization, the 
main computational burden is obtaining the gradient and Hessian approximation, 
which remain fixed within~\eqref{dmcon}. 

The scaled gradient projection framework includes a variety of methods depending on the choice of $H^n$. For example, if $H^n=\frac{1}{\alpha}I$, we can solve \eqref{dmcon} with~\eqref{eq:projGrad} where the iterations correspond to projected gradients with stepsize of $\alpha$. When the $H^n$ are chosen to approximate the Hessian of $f$ at $m^n$, we arrive at projected Newton-like methods that incorporate second-order information. A good summary of some of the possible choices can be found in \cite{Schmidt2012}, which uses a projected quasi-Newton method proposed based on a limited-memory Broyden-Fletcher-Goldfarb-Shannon (l-BFGS) approximation of the Hessian, and solves the convex subproblems for each update with a spectral-projected gradient method. The extended formulation~\eqref{eq:ENLLS} fits naturally into this framework, 
and allows for accurate sparse Gauss-Newton Hessian approximations. We exploit this structure with the development of a dedicated solver.


\subsection*{Implicit trust-region approach}

As we mentioned earlier, objective function evaluations in PDE-constrained optimization can be expensive. Therefore, we use an implicit trust region method, which avoids expensive linesearches. For this purpose, we replace $H^n$ by $H^n + c_nI$ while adjusting the damping parameter $c_n$ at each iteration adaptively and rejecting iterations that do not lead to a sufficient decrease in the objective.

When $\nabla f$ is Lipschitz continuous, i.e. for some constant $K$ 
\begin{equation*}
\|\nabla f(x) - \nabla f(y)\| \leq K \|x - y\| \ \text{for all } x,y \in C ,
\end{equation*}
and $H$ is a symmetric matrix, then $c$ can be chosen large enough so that
\begin{equation}
f(m+\dm) - f(m) \leq \dm^\top \nabla f(m) + \frac{1}{2} \dm^\top (H + c I)\dm
\label{majmin}
\end{equation}
for any $m \in C$ and $\dm$ such that $m + \dm \in C$. When solving~\eqref{dmcon}, $\Delta m = 0$ is feasible for the right-hand side, so if~\eqref{majmin} is satisfied we have $f(m+\dm) - f(m) \leq 0$. If moreover the solution of~\eqref{dmcon} satisfies $\|\Delta m\| > 0$, we have
\[
f(m+\dm) - f(m) < 0,
 \]
so strict descent (or stationarity) is guaranteed for method~\eqref{dmcon} under condition~\eqref{majmin}.  

We can derive explicit conditions that ensure~\eqref{majmin} is satisfied. Note that since 
\begin{equation*}
f(m+\dm) - f(m) \leq \dm^\top \nabla f(m) + \frac{K}{2}\|\dm\|^2 \ ,
\end{equation*}
it follows that
\begin{equation*}
\begin{aligned}
f(m+\dm) - f(m) & \leq \frac{1}{2}(K - \lambda_H^{\text{min}} - c)\|\dm\|^2 +  \dm^\top \nabla f(m) + \frac{1}{2} \dm^\top (H + c I)\dm \ ,
\end{aligned}
\end{equation*}
where $\lambda_H^{\text{min}}$ denotes the smallest eigenvalue of $H$, so choosing $c > K - \lambda_H^{\text{min}}$ ensures~\eqref{majmin}. However, this extremely conservative choice of the damping parameter $c$ may lead to a slow rate of convergence, so instead one can choose $c$ to be as small as possible while still decreasing the objective by a sufficient amount:
\begin{equation}
f(m+\dm) - f(m) \leq \sigma \left(\dm^\top \nabla f(m) + \frac{1}{2} \dm^\top (H + c I)\dm\right) \ ,
\label{suffdec}
\end{equation}
for some $\sigma \in (0,1]$. Using the same framework as in \cite{Esser2013}, the resulting method is summarized in Algorithm~\ref{SGPalg}.

\begin{scholmdAlgorithm}
~~~~~~$n=0$;~$m^0 \in C$;~$\rho > 0$;~$\epsilon > 0$;~$\sigma \in (0,1]$;\\\hspace*{0.333em}\hspace*{0.333em}\hspace*{0.333em}\hspace*{0.333em}\hspace*{0.333em}\hspace*{0.333em}$H$~symmetric~with~eigenvalues~between~$\lambda_H^{\text{min}}$~and~$\lambda_H^{\text{max}}$;

~~~~~~$\xi_1 > 1$;~$\xi_2 > 1$;~$c_0 > \max(0,\rho - \lambda_H^{\text{min}})$;\\\hspace*{0.333em}\hspace*{0.333em}\hspace*{0.333em}\hspace*{0.333em}\hspace*{0.333em}\hspace*{0.333em}while~$n=0$~or~$\frac{\|m^n - m^{n-1}\|}{\|m^n\|} > \epsilon$

~~~~~~~~~~$\displaystyle\dm = \argmin_{\dm + m^n \in C} \dm^\top \nabla f(m^n) + \frac{1}{2} \dm^\top(H^n + c_n I)\dm$

~~~~~~~~~~if~$f(m^n + \dm) - f(m^n) > \sigma (\dm^\top \nabla f(m^n) + \frac{1}{2} \dm^\top(H^n + c_n I)\dm)$\\\hspace*{0.333em}\hspace*{0.333em}\hspace*{0.333em}\hspace*{0.333em}\hspace*{0.333em}\hspace*{0.333em}\hspace*{0.333em}\hspace*{0.333em}\hspace*{0.333em}\hspace*{0.333em}\hspace*{0.333em}\hspace*{0.333em}\hspace*{0.333em}\hspace*{0.333em}$c_n = \xi_2 c_n$\\\hspace*{0.333em}\hspace*{0.333em}\hspace*{0.333em}\hspace*{0.333em}\hspace*{0.333em}\hspace*{0.333em}\hspace*{0.333em}\hspace*{0.333em}\hspace*{0.333em}\hspace*{0.333em}else~\\\hspace*{0.333em}\hspace*{0.333em}\hspace*{0.333em}\hspace*{0.333em}\hspace*{0.333em}\hspace*{0.333em}\hspace*{0.333em}\hspace*{0.333em}\hspace*{0.333em}\hspace*{0.333em}\hspace*{0.333em}\hspace*{0.333em}\hspace*{0.333em}\hspace*{0.333em}$m^{n+1} = m^n + \dm$\\\hspace*{0.333em}\hspace*{0.333em}\hspace*{0.333em}\hspace*{0.333em}\hspace*{0.333em}\hspace*{0.333em}\hspace*{0.333em}\hspace*{0.333em}\hspace*{0.333em}\hspace*{0.333em}\hspace*{0.333em}\hspace*{0.333em}\hspace*{0.333em}\hspace*{0.333em}$c_{n+1} = \begin{cases} \frac{c_n}{\xi_1} & \text{if } \frac{c_n}{\xi_1} > \max(0,\rho - \lambda_H^{\text{min}}) \\ c_n & \text{otherwise} \end{cases}$\\\hspace*{0.333em}\hspace*{0.333em}\hspace*{0.333em}\hspace*{0.333em}\hspace*{0.333em}\hspace*{0.333em}\hspace*{0.333em}\hspace*{0.333em}\hspace*{0.333em}\hspace*{0.333em}\hspace*{0.333em}\hspace*{0.333em}\hspace*{0.333em}\hspace*{0.333em}Define~$H^{n+1}$~to~be~symmetric~Hessian~approximation~\\\hspace*{0.333em}\hspace*{0.333em}\hspace*{0.333em}\hspace*{0.333em}\hspace*{0.333em}\hspace*{0.333em}\hspace*{0.333em}\hspace*{0.333em}\hspace*{0.333em}\hspace*{0.333em}\hspace*{0.333em}\hspace*{0.333em}\hspace*{0.333em}\hspace*{0.333em}\hspace*{0.333em}\hspace*{0.333em}\hspace*{0.333em}\hspace*{0.333em}with~eigenvalues~between~$\lambda_H^{\text{min}}$~and~$\lambda_H^{\text{max}}$\\\hspace*{0.333em}\hspace*{0.333em}\hspace*{0.333em}\hspace*{0.333em}\hspace*{0.333em}\hspace*{0.333em}\hspace*{0.333em}\hspace*{0.333em}\hspace*{0.333em}\hspace*{0.333em}\hspace*{0.333em}\hspace*{0.333em}\hspace*{0.333em}\hspace*{0.333em}$n = n + 1$\\\hspace*{0.333em}\hspace*{0.333em}\hspace*{0.333em}\hspace*{0.333em}\hspace*{0.333em}\hspace*{0.333em}\hspace*{0.333em}\hspace*{0.333em}\hspace*{0.333em}\hspace*{0.333em}end~if\\
\hspace*{0.333em}~~~~~~end~while
\caption{A Scaled Gradient Projection Algorithm for
(\ref{Fmincon})}\label{SGPalg}
\end{scholmdAlgorithm}


If $\nabla f$ is Lipschitz continuous and level sets of $f(m)$ intersected with $C$ are bounded, then any limit point $m^*$ of the sequence of iterates $\{m^n\}$ defined by Algorithm~\ref{SGPalg} is a stationary point of (\ref{Fmincon}), i.e. $(m - m^*)^\top \nabla f(m^*) \geq 0$ for all $m \in C$.

To implement Algorithm~\ref{SGPalg}, we need to specify both the Hessian approximation and a way to solve subproblem~\eqref{dmcon}. In the remainder of the paper, we show how to incorporate (1) box, (2) TV and (3) one-sided TV (hinge-loss) constraints, and explain how to solve~\eqref{dmcon} in each of the proposed formulations. We end this section with an important example where~\eqref{dmcon} is very simple to implement.

\subsection*{Example: diagonal Hessian approximation with box constraints}
When $H^n$ is diagonal and positive (e.g. approximation of the Gauss-Newton Hessian of~\eqref{eq:ENLLS}), Algorithm~\ref{SGPalg} simplifies. 
In this case, the subproblem 
\begin{equation}
\begin{aligned}
\dm & = \argmin_{\dm} \dm^\top \nabla f(m^n) + \frac{1}{2}\dm^\top (H^n + c_n I) \dm  \\
& \text{subject to} \ m^n_i + \dm_i \in [b_i, B_i]
\end{aligned}
\label{bounds}
\end{equation}
has the following closed form solution:
\begin{equation*}
\dm_i = \max\left(b_i - m^n_i,\min\left(B_i - m^n_i,-\widetilde{\dm}_i\right) \right),\quad i=1,\cdots, M, 
\end{equation*}
with $\widetilde{\dm}=(H^n + c_n I)^{-1}\nabla f(m^n)$. In this expression, $b_i$ and $B_i$ are the lower and upper bounds for $m_i$.

\section{Total-Variation constraints}\label{total-variation-regularization}
In FWI, even with box constraints, the recovered model can still contain artifacts and spurious oscillations. As we will demonstrate below, inaccuracies in $m$ can be reduced via convex constraints that bound the size of the total-variation norm ball to some positive value $\tau$.

TV-norm regularization, via penalties, constraints, or objectives is widely used in image processing to remove noise while preserving discontinuities \cite{Rudin1992}. TV-norm minimization also plays an important role in a variety of other inverse problems, especially when the unknown model parameters can be represented by piecewise constant or piecewise smooth functions. Examples include electrical impedance tomography \cite{Chung2005}, inverse wave propagation \cite{Akcelik2002}, and more recently in FWI promoting blockiness \cite{Guitton2012}, shape optimization \cite{GUO2012}, and time-lapse data \cite{Maharramov2014}. These approaches use a TV-penalized formulation, while we introduce the regularization as a constraint.
%
%
The constrained formulation keeps all iterates confined to a pre-defined convex set $C$, restricting 
the feasible model space. Penalized formulations do not offer this guarantee, which proves to be essential in FWI. 

If we represent $m$ as a $M_x\times M_y$ array, we can define the discrete TV-norm as 
\begin{equation}
\begin{aligned}
\|m\|_{TV} & = \frac{1}{h}\sum_{k,l}\sqrt{(m_{k+1,l}-m_{k,l})^2 + (m_{k,l+1}-m_{k,l})^2} \\
& = \sum_{k,l}\frac{1}{h}\left\| \bbm m_{k,l+1}-m_{k,l}\\ m_{k+1,l}-m_{k,l} \ebm \right\| \ ,
\end{aligned}
\label{TVim}
\end{equation}
which is the sum of the $\ell_2$ norms of the discrete gradient vectors at each point in the discretized model. We assume Neumann boundary conditions so that these differences are zero at the boundary. We arrive at a more compact expression for $\|m\|_{TV}$ if we define the finite difference operator $D$ such $(Dm)_i$ is the discrete gradient at location indexed by $i = 1,..., M$, where  $M=M_x\times M_y$. We now define 
\begin{equation}
\|m\|_{TV} = \|Dm\|_{1,2} := \sum_{i=1}^M \|(Dm)_i\| \ .
\label{TVvec}
\end{equation}
If we impose this TV-norm constraint in addition to the box constraints defined earlier, the set $C$ is an intersection
 \[
 C =  [b_i , B_i] \cap \{\|m\|_{TV} \leq \tau\}.
 \] 
The corresponding model updates~\eqref{dmcon} required by Algorithm~\ref{SGPalg} are given by 
\begin{equation}
\begin{aligned}
\dm & = \argmin_{\dm} \dm^\top  \nabla f(m^n) + \frac{1}{2}\dm^\top  (H^n + c_n I) \dm  \\
& \text{subject to } m^n_i + \dm_i \in [b_i , B_i] \text{ and } \|m^n + \dm\|_{TV} \leq \tau \\
m^{n+1} & = m^n + \dm \ .
\end{aligned}
\label{pTVb}
\end{equation}
\subsection*{Example: Projecting the Marmousi model on the intersection of box and TV-norm constraints}

\begin{figure}
\centering
\subfloat[\label{fig:marm_a}]{\includegraphics[width=0.330\hsize, natwidth = 650 ,natheight=642]{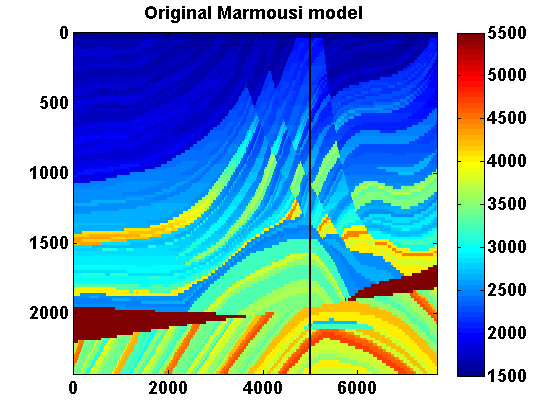}}
\subfloat[\label{fig:marm_b}]{\includegraphics[width=0.330\hsize, natwidth = 650 ,natheight=642]{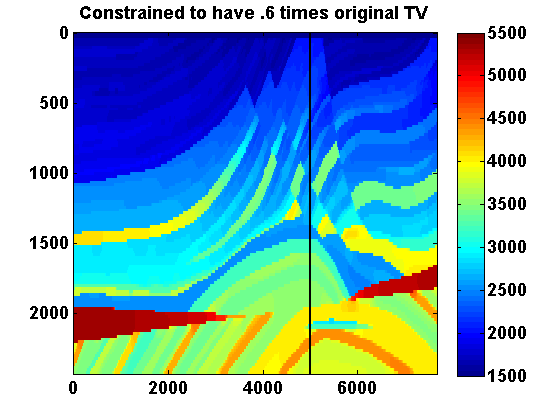}}
\subfloat[\label{fig:marm_c}]{\includegraphics[width=0.330\hsize, natwidth = 650 ,natheight=642]{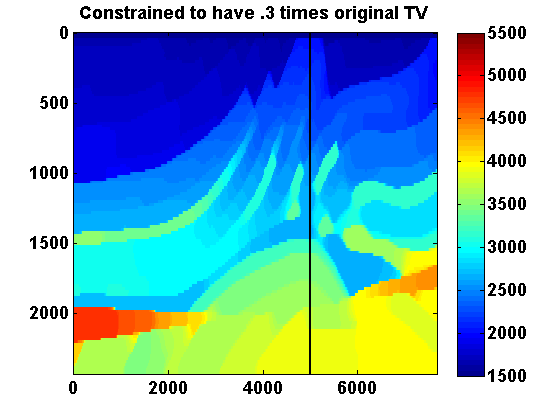}}
\caption{Marmousi model (a) and projected Marmousi model for
$\tau = .6\tau_0$ (b) and $\tau = .3\tau_0$ (c).}\label{fig:marm_proj}
\end{figure}

\begin{figure}
\centering
\subfloat[\label{fig:marm_slowness_slices}]{\includegraphics[width=0.490\hsize, natwidth = 650 ,natheight=642]{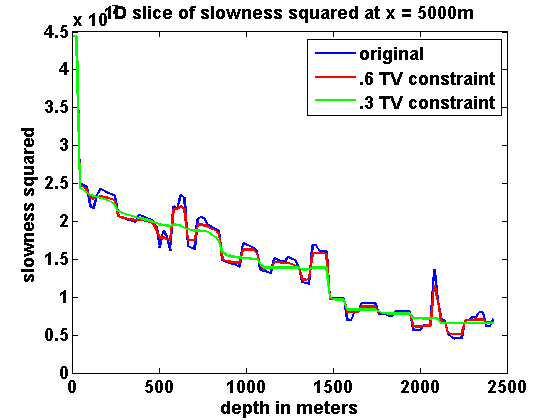}}
\subfloat[\label{fig:marm_velocity_slices}]{\includegraphics[width=0.490\hsize, natwidth = 650 ,natheight=642]{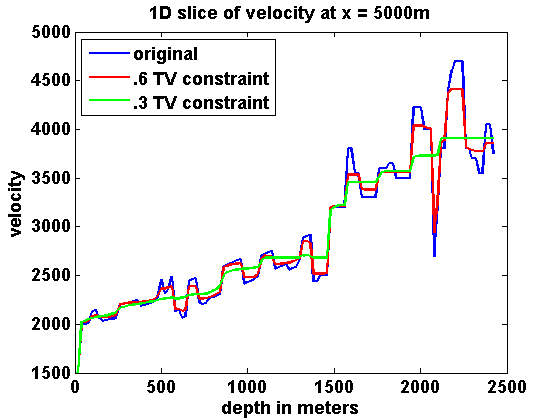}}
\caption{Comparison of slices from the Marmousi model and its
projections onto different TV constraints both in terms of slowness
squared (a) and velocity (b).}\label{fig:marm_slices}
\end{figure}

Before we discuss how to minimize \eqref{Fmincon} in the seismic setting, we motivate the use of TV-norm constraints in seismic inverse problems. Consider projecting the Marmousi model \cite{Bourgeois91}, shown in Figure~\ref{fig:marm_a}, onto two sets with TV-norm constraints using decreasing values of $\tau$. Let $m_0$ denote the original Marmousi model and let $\tau_0 = \|m_0\|_{TV}$. For the bound constraints on the slowness squared, set $B_i = 4.4444\e{-7}\,\mathrm{s^2m^{-2}}$ everywhere, which corresponds to a lower bound of $1500\,\mathrm{ms^{-1}}$ on the velocity. Taking advantage of the fact that these constraints can vary spatially, let $b_i = 4.4444\e{-7}\,\mathrm{s^2m^{-2}}$ in the water layer and $b_i = 3.3058\e{-8}\,\mathrm{s^2m^{-2}}$ everywhere else, which corresponds to an upper bound of $5500\,\mathrm{ms^{-1}}$ on the velocity. The orthogonal projection of $m_0$ onto the intersection of these box and TV-norm ball constraints is defined by
\begin{equation}
\begin{aligned}
\Pi_C(m_0) & = \argmin_m \frac{1}{2} \|m - m_0\|^2 \\
& \text{subject to } m_i \in [b_i , B_i] \text{ and } \|m\|_{TV} \leq \tau \ .
\end{aligned}
\label{proj_C}
\end{equation}
Results with $\tau = .6\tau_0$ and $\tau = .3\tau_0$ are shown in Figure~\ref{fig:marm_proj}. The vertical lines at $x = 5000\,\mathrm{m}$ indicate the location of the 1D vertical slices shown in Figure~\ref{fig:marm_slices} for both slowness squared (a) and the velocity (b). As we decrease the size of the TV-norm ball, the feasible models that are close to the original model ($m_0$) in the 2-norm become more and more "cartoon like" with fewer and fewer unconformities. The size of the TV-norm ball ($\tau$) controls the complexity of the model while still allowing for discontinuous unconformities. We exploit this property of the TV-norm extensively.

\subsection*{Solving the convex subproblems}\label{solving-the-convex-subproblems}

The proposed approach requires solving the quadratic approximation in \eqref{dmcon} 
with a given gradient and Hessian approximation. 
A computationally efficient approach for the convex subproblems of type~\eqref{pTVb}
is the primal-dual hybrid gradient (PDHG) method \cite{Zhu2008} studied in \cite{Esser2010,Chambolle2011,He2012,Zhang2010a}. To develop 
the method, we first write down the Lagrangian for a dualization of~\eqref{pTVb}: 
\begin{equation}
\begin{aligned}
\L(\dm,p) & = \dm^\top  \nabla f(m^n) + \frac{1}{2}\dm^\top  (H^n + c_n I) \dm + g_B(m^n + \dm)\\
& + p^\top  D(m^n + \dm) - \tau \|p\|_{\infty,2}. 
\end{aligned}
\label{Lagrangian}
\end{equation}
In~\eqref{Lagrangian}, $p$ is the vector of Lagrange multipliers, and $g_B$ is an indicator function for the bound constraints
\begin{equation*}
g_B(m) = \begin{cases} 0 & \quad\text{if}\quad m_i \in [b_i , B_i] \\ 
\infty & \quad\text{otherwise.} \end{cases}
\end{equation*}
Here, $\|\cdot\|_{\infty,2}$ uses mixed norm notation to denote the
dual norm of $\|\cdot\|_{1,2}$. It takes the $\max$ instead of the sum of
the $\ell_2$ norms so that $\|Dm\|_{\infty,2} = \max_i \|(Dm)_i\|$ in the
notation of~\eqref{TVvec}. This saddle point problem can be
derived from the convex subproblem in (\ref{pTVb}) using 
the conjugate representation of the TV constraint:
\begin{equation}
 \sup_p \left\{p^\top  D(m^n + \dm) - \tau \|p\|_{\infty,2}\right\} \ ,
\label{TVsp}
\end{equation}
which equals the indicator function
\begin{equation*}
\begin{cases}
0 & \quad \text{if}\quad \|D(m^n + \dm)\|_{1,2} \leq \tau \\
\infty & \quad \text{otherwise.}
\end{cases}
\end{equation*}
To find a saddle point of (\ref{Lagrangian}), the modified PDHG 
iterations are given by 
\begin{equation}
\begin{aligned}
p^{k+1} & = \argmin_p \tau \|p\|_{\infty,2} - p^\top  D(m^n + \dm^k) + \frac{1}{2 \delta}\|p - p^k\|^2 \\
\dm^{k+1} & = \argmin_{\dm} \dm^\top  \nabla f(m^n) + \frac{1}{2}\dm^\top  (H^n + c_nI) \dm \\
& + \dm^\top  D^\top (2p^{k+1}-p^k) + \frac{1}{2 \alpha}\|\dm - \dm^k\|^2 \\
& \text{subject to } m^n_i + \dm_i \in [b_i , B_i] \ .
\end{aligned}
\label{PDHGM}
\end{equation}
These iterations can be written more explicitly as

\begin{equation}
\begin{aligned}
p^{k+1} & = p^k + \delta D(m^n + \dm^k) - \Pi_{\|\cdot\|_{1,2} \leq \tau \delta}(p^k + \delta D(m^n + \dm^k)) \\
\dm^{k+1}_i & =  \max\left( b_i - m^n_i, \min\left(B_i - m^n_i, \widetilde{\dm}_i\right) \right) \ ,
\end{aligned}
\label{PDHGexplicit}
\end{equation}
where 
\[
\widetilde{\dm}=(H^n + (c_n+\frac{1}{\alpha})I)^{-1}(-\nabla f(m^n) + \frac{\dm^k}{\alpha} - D^\top (2p^{k+1}-p^k))
\] 
and $\Pi_{\|\cdot\|_{1,2} \leq \tau \delta}(z)$ denotes the orthogonal projection of $z$ onto the ball of radius $\tau \delta$ in the $\|\cdot\|_{1,2}$ norm. This requires projecting the vector of $\ell_2$ norms of the spatial gradient vectors onto a simplex.  A simple approach to project $z$ onto the unit simplex $\{x:\, x_i \geq 0,\, \sum_i x_i = 1\}$ is to use bisection to find the threshold $a$ such that $\sum_i \max(0,z_i - a) = 1$, in which case $\max(0,z_i - a)$ is the $i^\text{th}$ component of the projection. An efficient linear time $O(n)$ implementation based on this idea is developed and analyzed in~\cite{Brucker1984}.

The step size restriction required for convergence is $\alpha \delta \leq \frac{1}{\|D^\top D\|}$ (see \cite{Esser2010} for details). If $h$ is the mesh width, then the eigenvalues of $D^\top D$ are between $0$ and $\frac{8}{h^2}$ by the Gershgorin Circle Theorem, so it suffices to choose positive $\alpha$ and $\delta$ such that $\alpha \delta \leq \frac{h^2}{8}$.

The relative scaling of $\alpha$ and $\delta$ can have a large effect on the convergence rate of the method \cite{Esser2010}. A reasonable choice for the fixed step size parameters is
\begin{equation*}
\alpha = \frac{1}{\max(H^n + c_n I)} \text{ and } \delta = \frac{h^2 \max(H^n +c_n I)}{8} \leq \frac{\max(H^n + c_nI)}{\|D^\top\!D\|}. 
\end{equation*}
However, this choice may be too conservative. The convergence rate of the method can be improved by using iteration-dependent step sizes as proposed by \cite{Chambolle2011}. The adaptive backtracking strategy introduced by \cite{Goldstein2013} can also be a practical way of choosing step size parameters that improve the convergence rate.

\section{Numerical Experiments}\label{numerical-experiments}

We consider four 2D numerical experiments based on synthetic data. The PDE in
this case is a scalar Helmholtz equation
\[
\left(\omega^2 m + \nabla^2\right)u = q \quad \mbox{in the interior of the domain}
\]
with radiation boundary conditions
\[
\left(\frac{\partial }{\partial n} - \imath\omega\right)u = 0, 
\]
where $\frac{\partial }{\partial n}$ denotes the derivative in the direction normal to the boundary. A standard finite difference discretization of the Helmholtz operator leads to a sparse banded matrix $A(m)$. The data are collected by an array of receivers for several point-sources and a range of angular frequencies $\omega$. For details
on discretization see the appendix.

In the experiments, we primarily use the extended formulation~\eqref{eq:ENLLS} and minimize the
objective defined in \eqref{eq:WRI}. The corresponding gradient and Hessian
expressions are given in \eqref{eq:gWRI} and \eqref{eq:HWRI}. We
refer to this approach as \emph{Wavefield Reconstruction Inversion} (WRI). We
also show the benefits of adding constraints to the conventional approach
of minimizing the objective defined in \eqref{eq:FWI}, where we use the
gradient and pseudo-Hessian are defined in \eqref{eq:gFWI} and \eqref{eq:HFWI}. 
We refer to the latter approach as \emph{Full-Waveform Inversion}
(FWI). A standard approach in FWI/WRI is to start the inversion from low
frequencies and gradually move to higher frequencies. This continuation strategy
 helps to avoid cycle-skipping related local minima \cite{Bunks1995}.

%

To illustrate the performance of the proposed constrained formalism, we consider two examples derived from the 2004 BP velocity benchmark data set \cite{Billette2005}. This data set was designed to evaluate the capabilities of velocity-analysis techniques for complex geologies that contain sedimentary basins, with increasing velocities, interspersed with high-contrast and high-velocity salt bodies. 

The first experiment tries to recover the top left portion of this synthetic model from an accurate smooth starting model. Even in this situation, the resulting inverted model has noisy artifacts related to the applied source encoding, and the deeper part of the estimated model tends to have incorrect discontinuities. The TV constraint helps to remove some of these artifacts while still recovering the most significant discontinuities. This example also demonstrates that relaxing the TV constraint 
over multiple passes through the frequency batches leads to improvements. A stronger TV constraint 
gives an oversmoothed model estimate, but this estimate then serves as a good initial model for 
future passes as the TV constraint is relaxed. 

The second experiment tries to recover the top middle portion of the same BP velocity model from a poor starting model. This experiment illustrates that the combination of edge-preserving TV constraint and box-constraints are inadequate, and motivates an additional asymmetric TV-norm constraint, tailored to geologies with sedimentary basins with salt inclusions. Our approach here also uses a continuation strategy that at first strongly discourages downward jumps in the estimated velocity, and then gradually relaxing this constraint. 
This strategy helps prevent the method from getting stuck in bad local minima at early passes, eventually 
allowing downward jumps in velocity for a better data fit.

We first present details on frequency continuation and {\it simultaneous shot} strategies used in the experiments. 

\subsection*{Frequency Continuation}\label{frequency-continuation}

We work with small batches of frequency data at a time, moving from low to high frequencies in overlapping batches of two. This frequency continuation strategy does not guarantee that we solve the overall problem after a single pass from low to high frequencies, but is far more computationally tractable than minimizing over all frequencies simultaneously. Moreover, the continuation strategy of moving from low to high frequencies helps prevent the iterates from tending towards bad local minima \cite{Bunks1995}.

For example, if the data consists of frequencies starting at $3\,\mathrm{Hz}$ sampled at intervals at $1\,\mathrm{Hz}$, then we would start with the $3$ and $4\,\mathrm{Hz}$ data, use the computed $m$ as an initial guess for inverting the $4$ and $5\,\mathrm{Hz}$ data and so on. For each frequency batch, we will compute at most $25$ outer iterations, each time solving the convex subproblem to convergence, stopping when
\[
\max\left(\frac{\|p^{k+1}-p^k\|}{\|p^{k+1}\|},\frac{\|\dm^{k+1}-\dm^k\|}{\|\dm^{k+1}\|}\right) \leq 1\e{-4}.
\]
Since the magnitude of the data depends on the frequency, one can incorporate frequency dependent weights in the definition of the objective. However, as we only work with small frequency batches in practice, these weights do not have a significant effect.

In some of our experiments we perform multiple passes through the frequencies to further improve the results.

\subsection*{Simultaneous sources}
The data are typically collected for many sources, resulting in a redundant
data set. Since the number of PDE-solves required at each iteration is dictated
by the number of sources, it makes sense to try to compress the data volume and
thus reduce the computational cost. One way of doing this is by randomly
projecting the data, also known as "source encoding", resulting in a smaller number of terms
\[
\widetilde{d}_i = \sum_{j=1}^{N_s} w_{ij} d_j, \quad \mathrm{for} \quad i = 1, 2, \ldots, \widetilde{N}_s, \quad \widetilde{N}_s\ll N_s
\]
where the $w_{ij}$'s are drawn from a standard normal distribution \cite{Krebs2009}. The source terms are similarly reduced (using the same weights)
\[
\widetilde{q}_i = \sum_{j=1}^{N_s} w_{ij} q_j, \quad \mathrm{for} \quad i = 1, 2, \ldots, \widetilde{N}_s.
\]
Various alternative strategies have been proposed to choose the random weights and tailor the optimization to deal specifically with the resulting stochasticity in the problem \cite{Haber2010,Friedlander2012,VanLeeuwen2013a}.

\subsection*{Example -- accurate starting model}

\begin{figure}
\centering
\subfloat[\label{fig:BPTL_vtrue}]{\includegraphics[width=0.480\hsize, natwidth = 650 ,natheight=642]{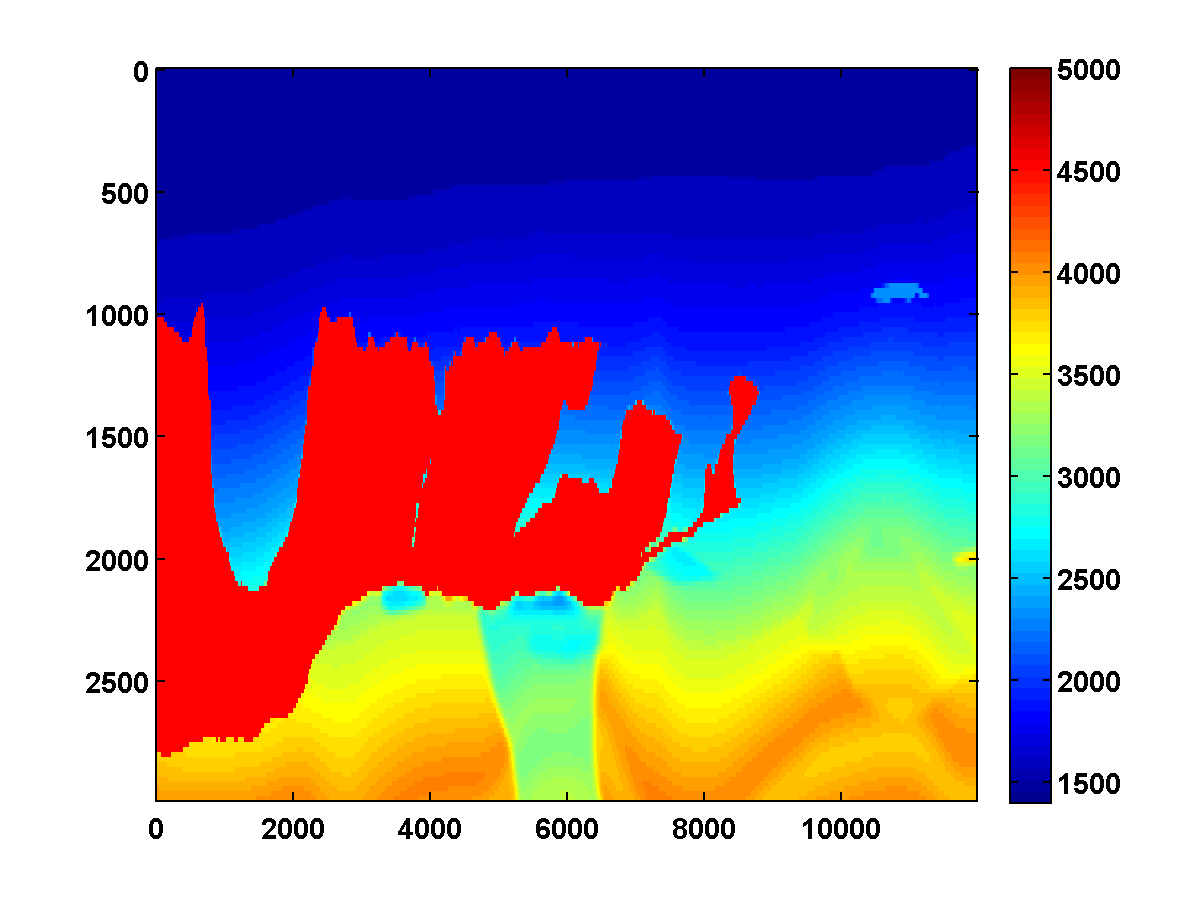}}
\subfloat[\label{fig:BPTL_vinit}]{\includegraphics[width=0.480\hsize, natwidth = 650 ,natheight=642]{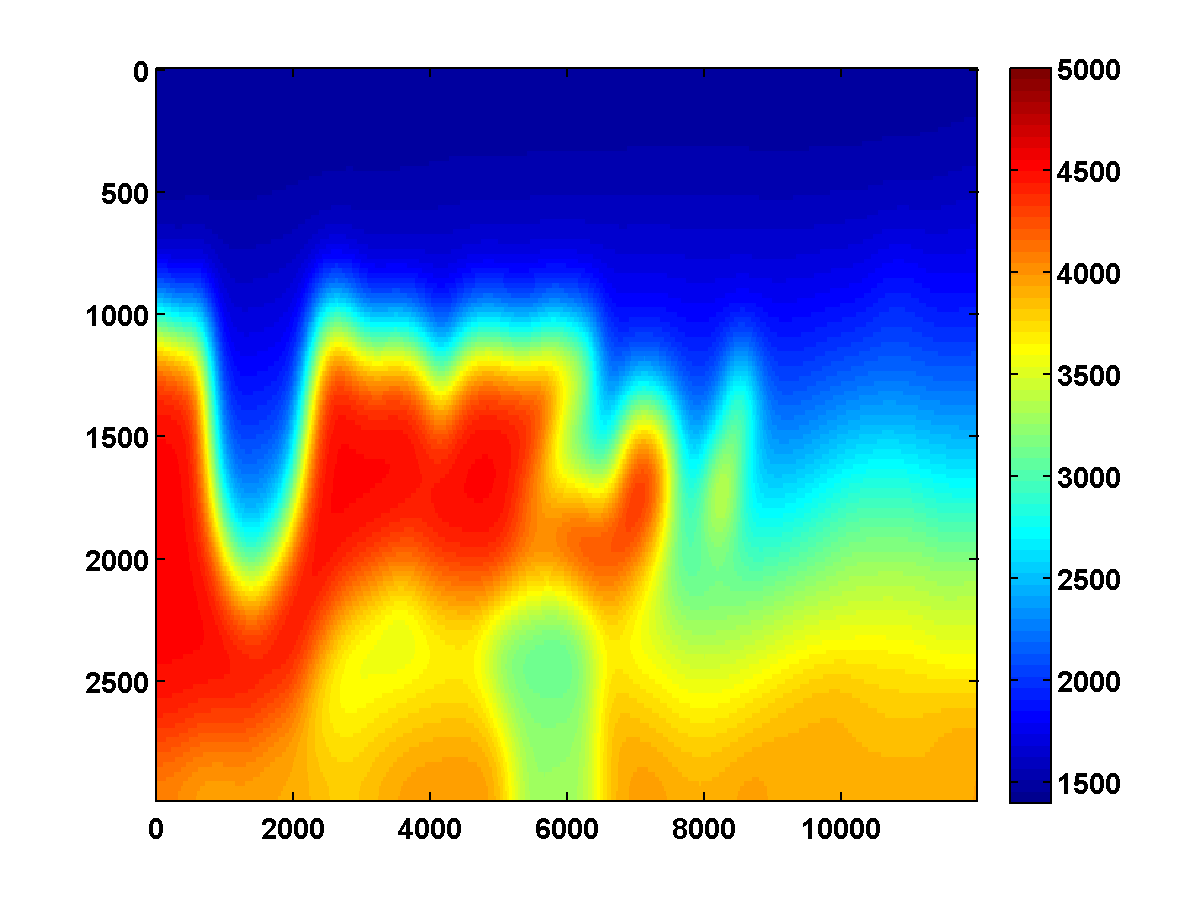}}
\caption{Top left portion of BP 2004 velocity model (a), and initial velocity (c).}\label{fig:BPTLsetup}
\end{figure}

\begin{figure}
\centering
\subfloat[\label{fig:BPTL_a}]{\includegraphics[width=0.480\hsize, natwidth = 650 ,natheight=642]{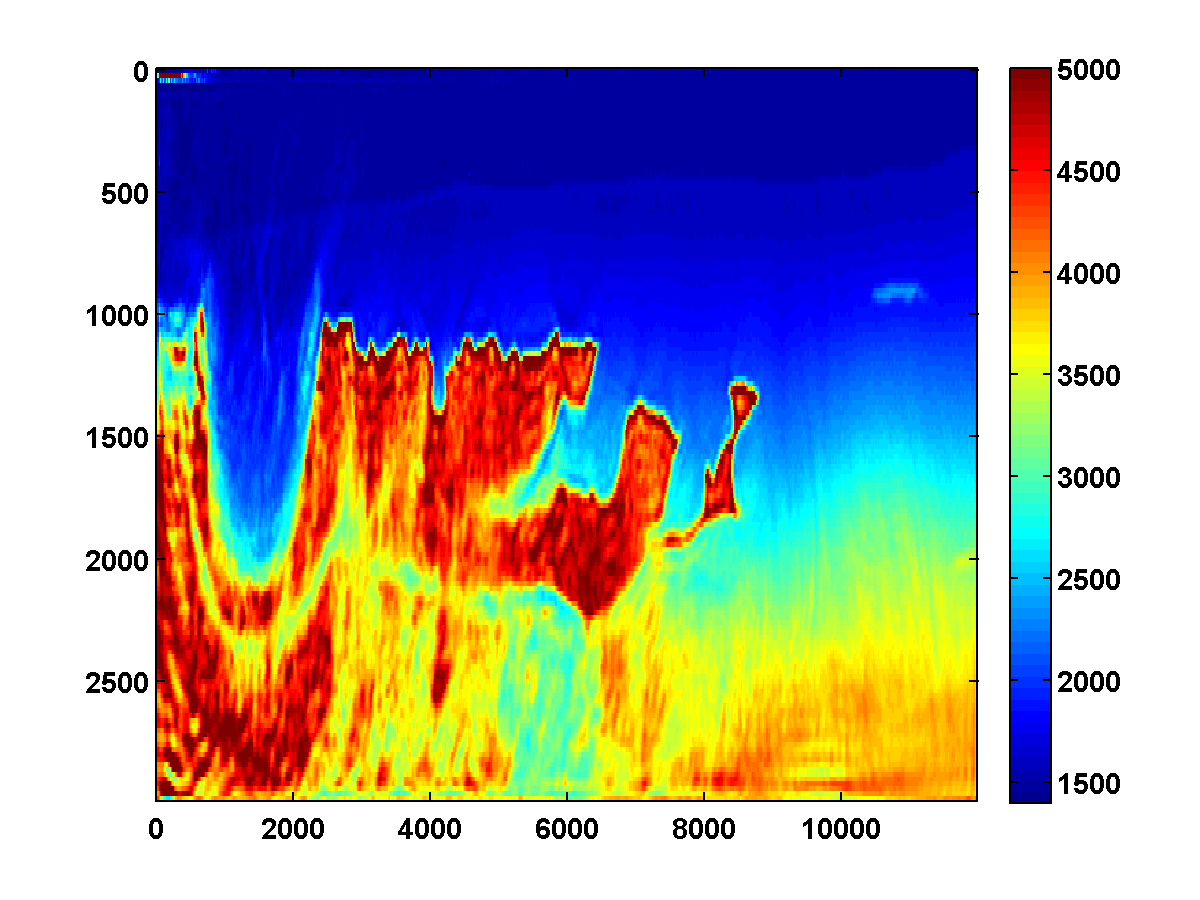}}
\subfloat[\label{fig:BPTL_b}]{\includegraphics[width=0.480\hsize, natwidth = 650 ,natheight=642]{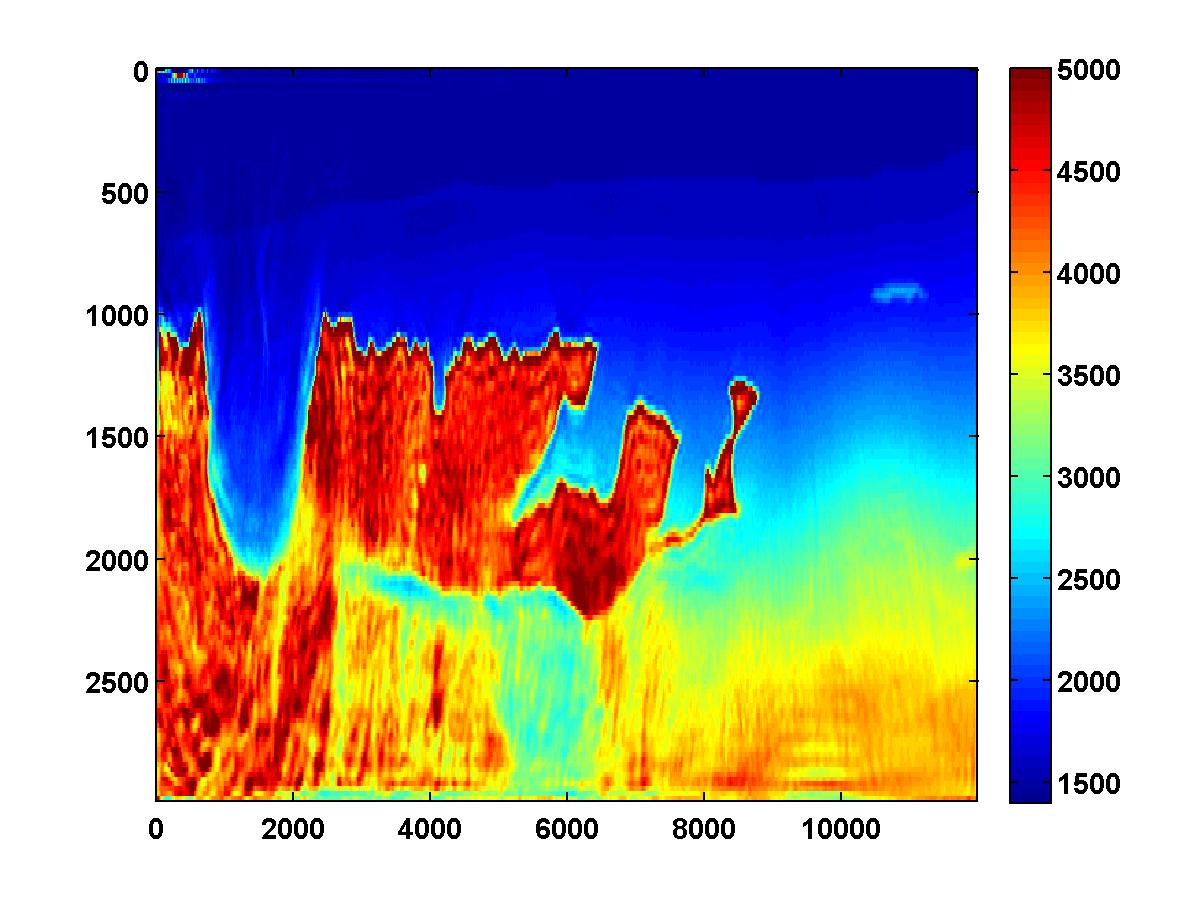}}
\caption{Recovered velocity without TV constraint from a good smooth initial model after one pass (a) and after two passes
(b) using small frequency batches from $3$ to $20$
Hz.}\label{fig:BPTL}
\end{figure}

To demonstrate the performance of the proposed formulation, we first consider a $3$ by $12\,\mathrm{km}$ portion of the the 2004 BP velocity benchmark model with a gridspacing of $20\,\mathrm{m}$, as shown in Figure~\ref{fig:BPTL_vtrue}, for a relative accurate smooth starting model shown in~\ref{fig:BPTL_vinit}. We use a Ricker wavelet with a peak frequency at $15\,\mathrm{Hz}$ as a for the source signature. There are 126 sources every $80\,\mathrm{m}$ between $1$ and $11\,\mathrm{km}$ at a depth of 40 m below the top of the model; 299 receivers are placed at $60\,\mathrm{m}$ below the surface and covering all the model in the horizontal dimension, in an equally spaced distribution (every $40\,\mathrm{m}$). We employ frequency continuation by looping through the frequencies ranging from $3$ to $20\,\mathrm{Hz}$ from low to high in overlapping batches of two frequencies each. We define bound constraints on the slowness squared that correspond to minimum and maximum velocities of $1400$ and $5000\,\mathrm{ms^{-1}}$, respectively. The inversion result for the velocity after one pass through the frequencies is included in Figure~\ref{fig:BPTL_a} and shows that the velocity model is reasonably well recovered (\textrm{cf.} Figures~\ref{fig:BPTL_vtrue} and \ref{fig:BPTL_a}). 
This result is obtained using WRI with box constraints, working with only two sources ($\widetilde{N}_s=2\ll N_s$). 
While encouraging, there are still visible noisy artifacts due to the simultaneous sources and from structural velocity lows in the salt. 
Without including additional constraints, the conventional adjoint-state approach to FWI is not able to produce a tangible result for this experiment.

\subsection*{Benefit of Multiple Passes}\label{benefit-of-multiple-passes}
Despite its evident shortcomings, after one pass through the frequencies (Figure~\ref{fig:BPTL_a}) 
the WRI result can serve as input for a second pass. Using multiple multiscale passes is used frequently in scientific computing (e.g.  V-cycles in multigrid), and \cite{peters2014EAGEweb} proposed this strategy for wave-equation based inversion. Comparing the inversion results after one and two passes shows significant improvements after the second pass (see Figure~\ref{fig:BPTL_b}) where the oscillation just below the top salt on the left side of the model mostly disappears after the second pass. 
The WRI approach with multiple passes allows useful results in situations where conventional FWI fails. 
However, the result remains noisy and suffers from artifacts near the boundary of the model due to poor illumination.

\subsection*{Including the TV constraint}

\begin{figure}
\centering
\subfloat[\label{fig:TVBPTL_a}]{\includegraphics[width=0.480\hsize, natwidth = 650 ,natheight=642]{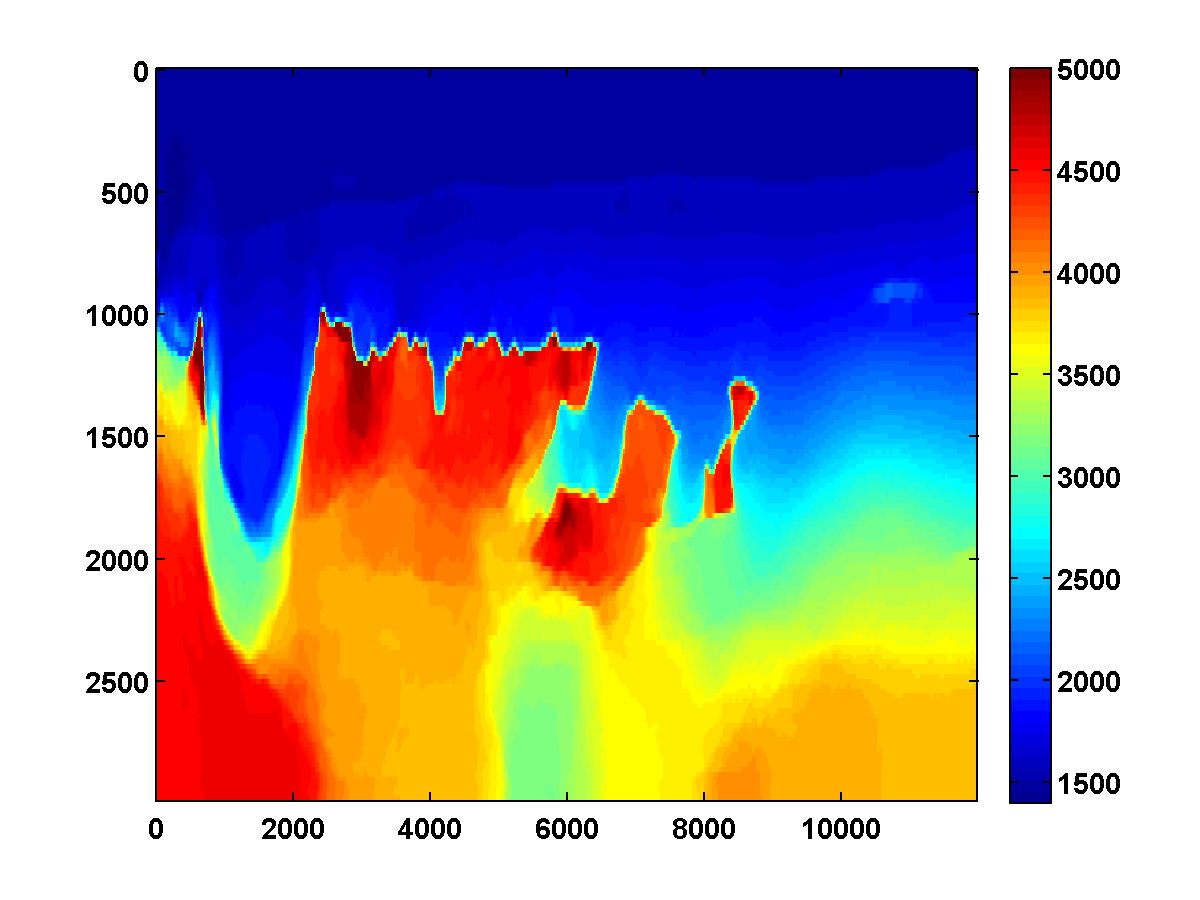}}
\subfloat[\label{fig:TVBPTL_b}]{\includegraphics[width=0.480\hsize, natwidth = 650 ,natheight=642]{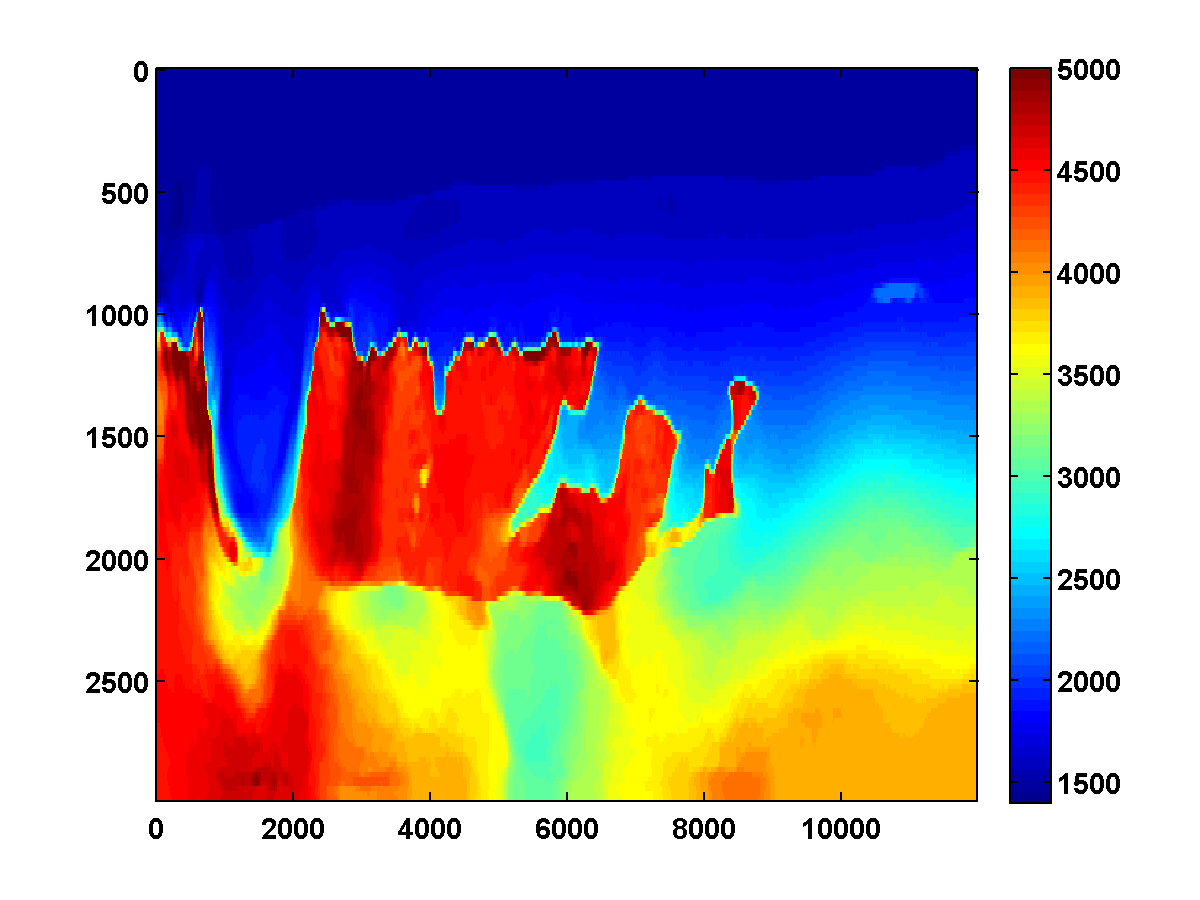}}
\caption{Recovered velocity with a TV constraint from a good smooth initial model after one pass with
$\tau = .9 \tau_{\text{true}}$(a) and a second pass with
$\tau = .99 \tau_{\text{true}}$ (b) using small frequency batches from
$3$ to $20$ Hz.}\label{fig:TVBPTL}
\end{figure}

Shortcomings of WRI with multiple passes include noisy artifacts lack of clear delineation at the top and bottom salt. As shown in Figure~\ref{fig:TVBPTL}, the inversion results improve significantly when we use the TV-norm constraint $\|m\|_{TV} \leq \tau$, together with a relaxation strategy over multiple passes. 
For the first pass, we choose $\tau$ to be $0.9$ times the TV-norm of the ground truth model. To reduce computation time, we again use two simultaneous shots, but now with Gaussian weights resampled every time the model is updated. The estimated model after $25$ outer iterations per frequency batch is shown in Figure~\ref{fig:TVBPTL_a}. Using this result as a warm start for a second pass through the frequency batches yields the improved result in Figure~\ref{fig:TVBPTL_b}. The top salt is much better resolved compared to the examples without the TV-norm constraint. The second pass also leads to a significant improvement of the bottom salt while preserving the velocity low below the salt related to an over-pressured reservoir. Because the model is off centre, issues with illumination remain but overall the interior of the salt body itself is well recovered and includes the small low-velocity inclusions that can be observed in the original model (Figure~\ref{fig:BPTL_vtrue}).

\subsection*{Example--poor starting model}

\begin{figure}
\centering
\subfloat[\label{fig:BPTM_vtrue}]{\includegraphics[width=0.480\hsize, natwidth = 650 ,natheight=642]{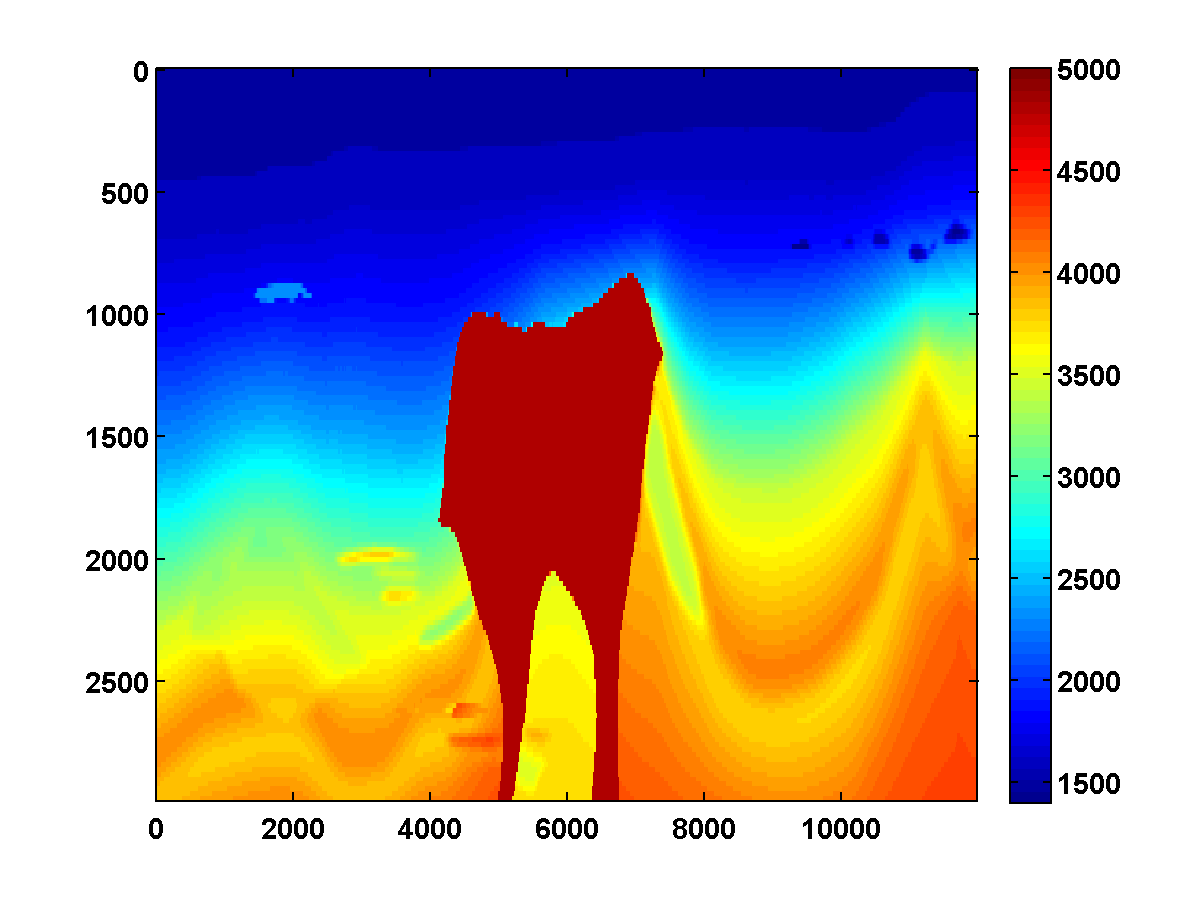}}
\subfloat[\label{fig:BPTM_vinit}]{\includegraphics[width=0.480\hsize, natwidth = 650 ,natheight=642]{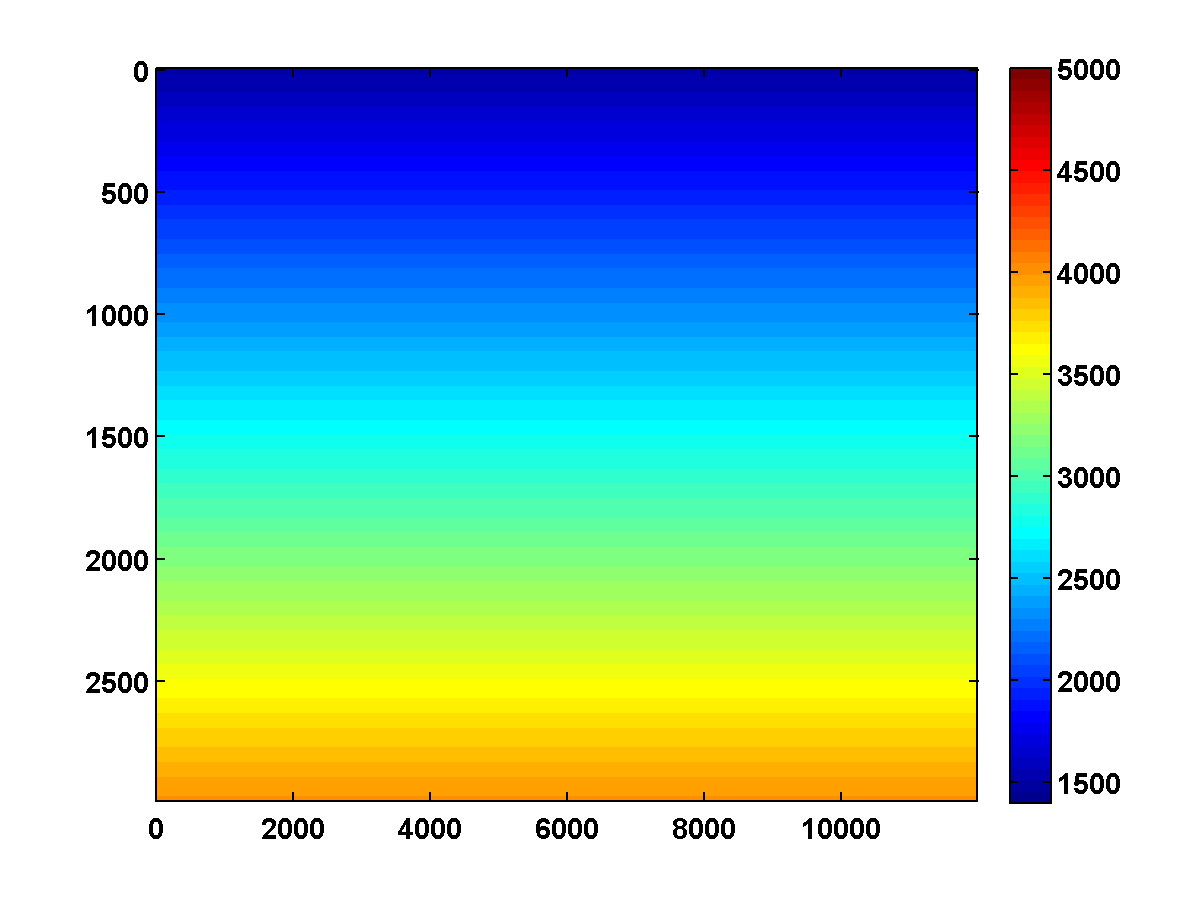}}
\caption{Top middle portion of BP 2004 velocity model (a), source and
receiver locations (b) and initial velocity (c).}\label{fig:BPTMsetup}
\end{figure}

The inversion results shown in Figure~\ref{fig:TVBPTL_b} relied heavily on having good initial models, and 
were obtained by smoothing the ground truth models. Consider the $3$ by $12\,\mathrm{km}$ velocity model shown in Figure~\ref{fig:BPTM_vtrue}, which is the top middle portion of the 2004 BP velocity benchmark data set, also sampled with a gridspacing of $20\,\mathrm{m}$; our acquisition geometry as well as the source signature is the identical to the one defined in the previous example, but with better illumination. The approach we used to recover the top left portion of the model also works well here using a smoothed version of the true model as the starting point. However, starting with a poor initial model as in Figure~\ref{fig:BPTM_vinit}, the method 
obtains a poor inversion result, possibly because of finding a parasitic stationary point. 
This happens despite the fact that the salt body is well illuminated.

\begin{figure}
\centering
\subfloat[\label{fig:BPTM_a}]{\includegraphics[width=0.330\hsize, natwidth = 650 ,natheight=642]{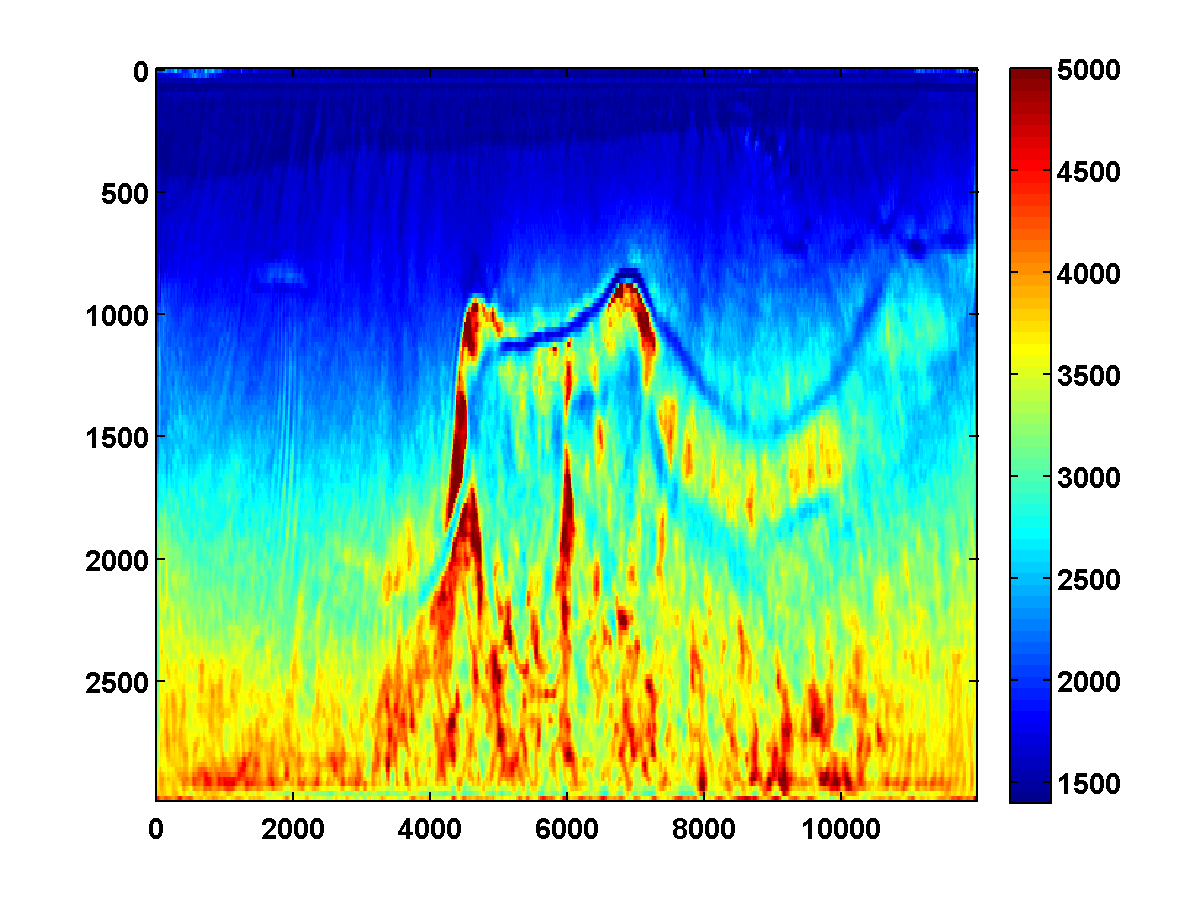}}
\subfloat[\label{fig:BPTM_b}]{\includegraphics[width=0.330\hsize, natwidth = 650 ,natheight=642]{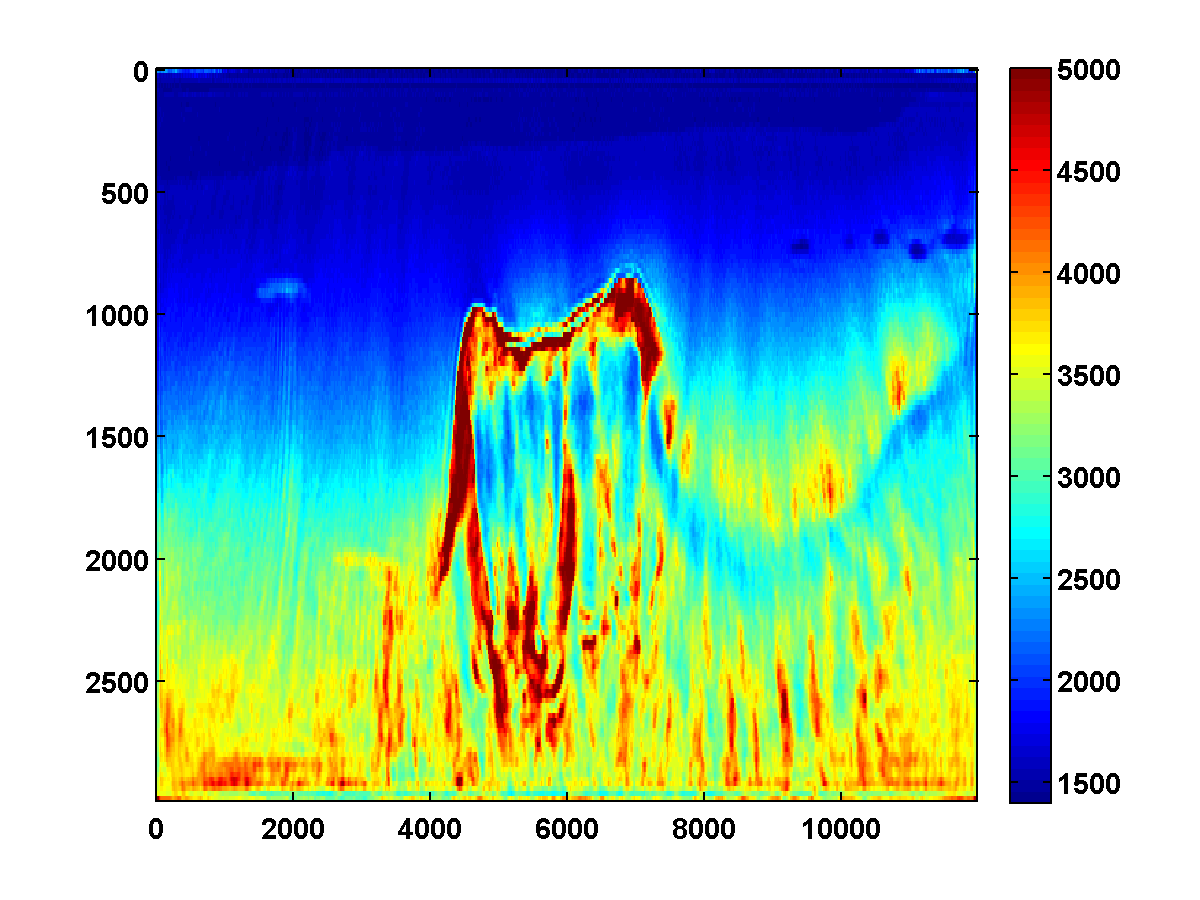}}
\subfloat[\label{fig:BPTM_c}]{\includegraphics[width=0.330\hsize, natwidth = 650 ,natheight=642]{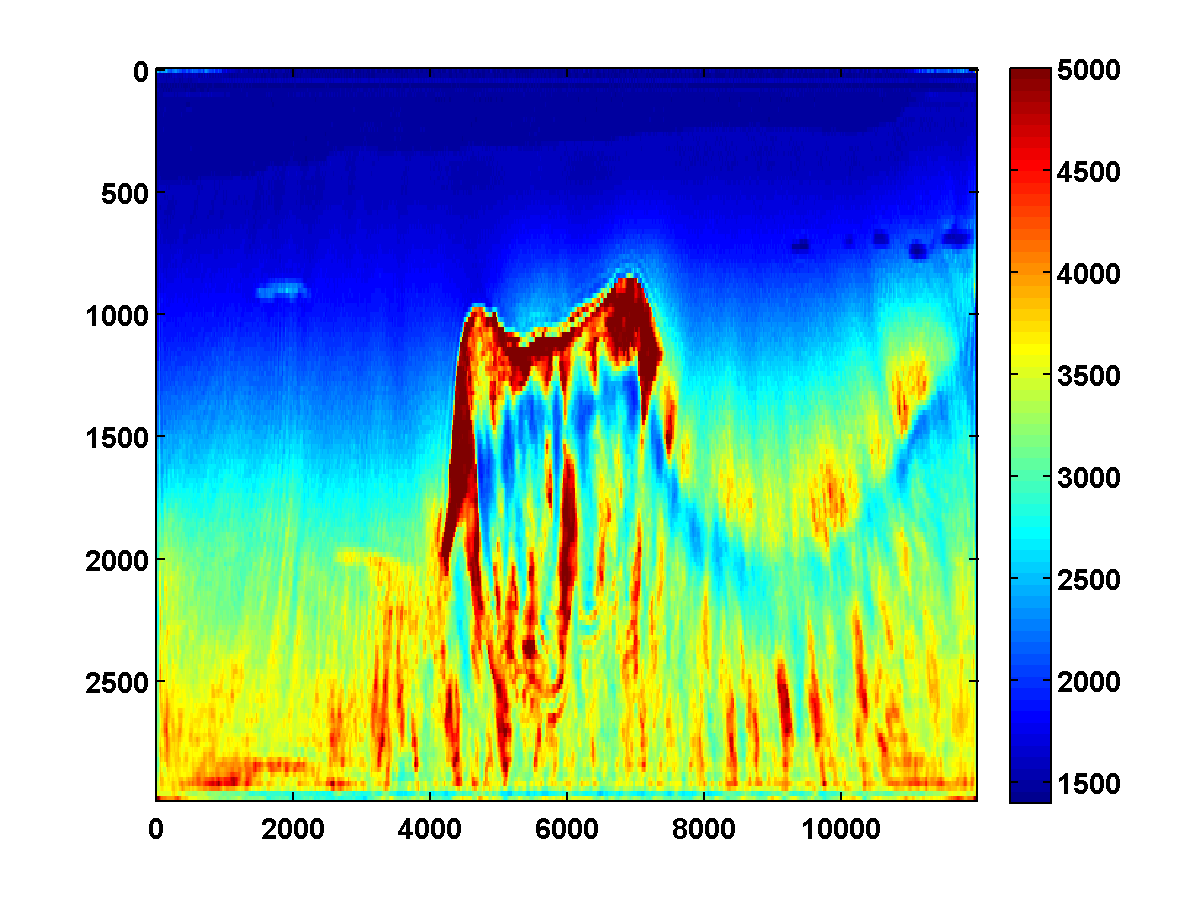}}
\caption{Recovered velocity with a no TV constraint from a poor initial model
after one pass (a) after a second pass (b) and
after a third pass (c) through small frequency batches from $3$ to $20$
Hz.}\label{fig:BPTM}
\end{figure}

As shown in Figure~\ref{fig:BPTM}, the WRI method with bound constraints alone yields a noisy result that fails to improve significantly even after multiple passes through the frequency batches. 
The effect of the initial updates, which decrease the velocity the after stepping into the salt, is persistent
and can not be overcome. This behavior is typical for data that misses the low frequencies, as we also observed in Figure~\ref{fig:example2} of the introduction.

\begin{figure}
\centering
\subfloat[\label{fig:TVBPTM_a}]{\includegraphics[width=0.330\hsize, natwidth = 650 ,natheight=642]{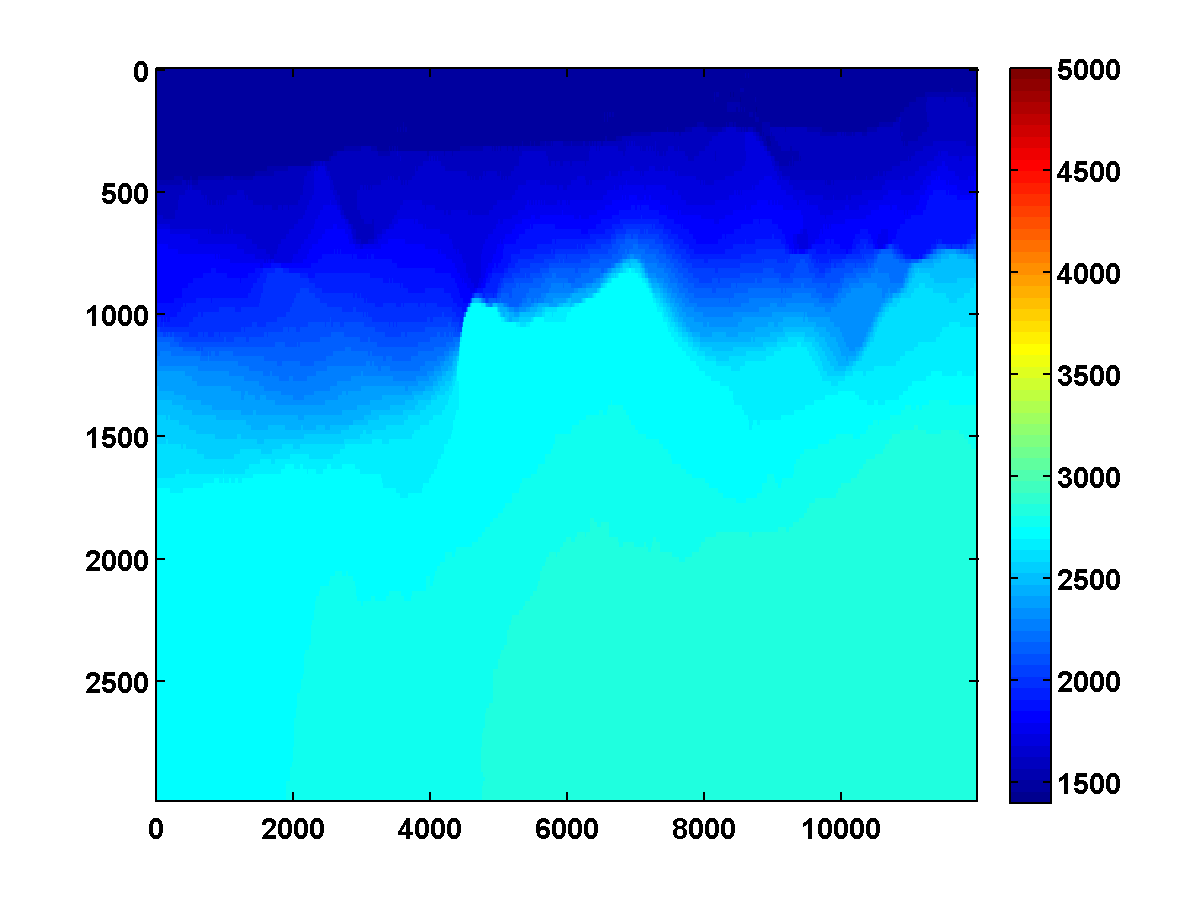}}
\subfloat[\label{fig:TVBPTM_b}]{\includegraphics[width=0.330\hsize, natwidth = 650 ,natheight=642]{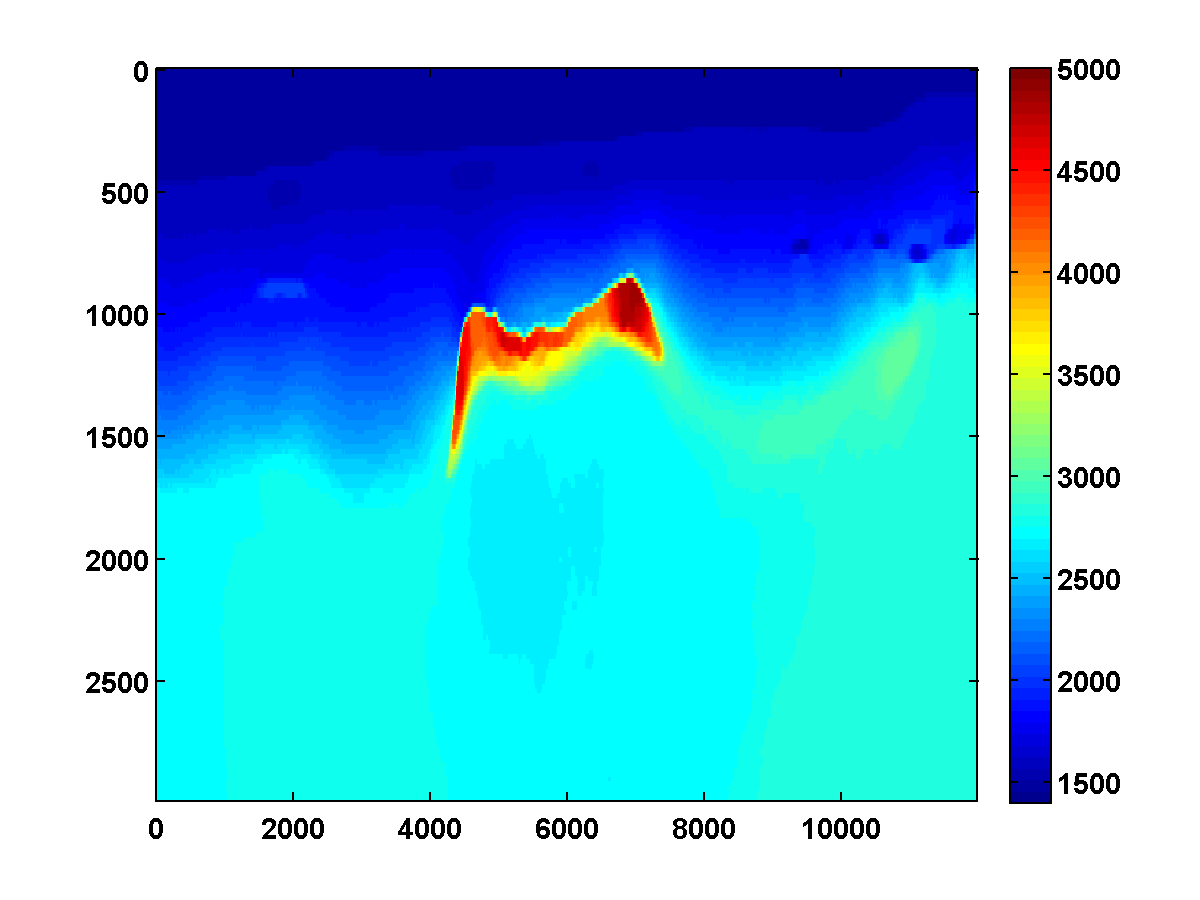}}
\subfloat[\label{fig:TVBPTM_c}]{\includegraphics[width=0.330\hsize, natwidth = 650 ,natheight=642]{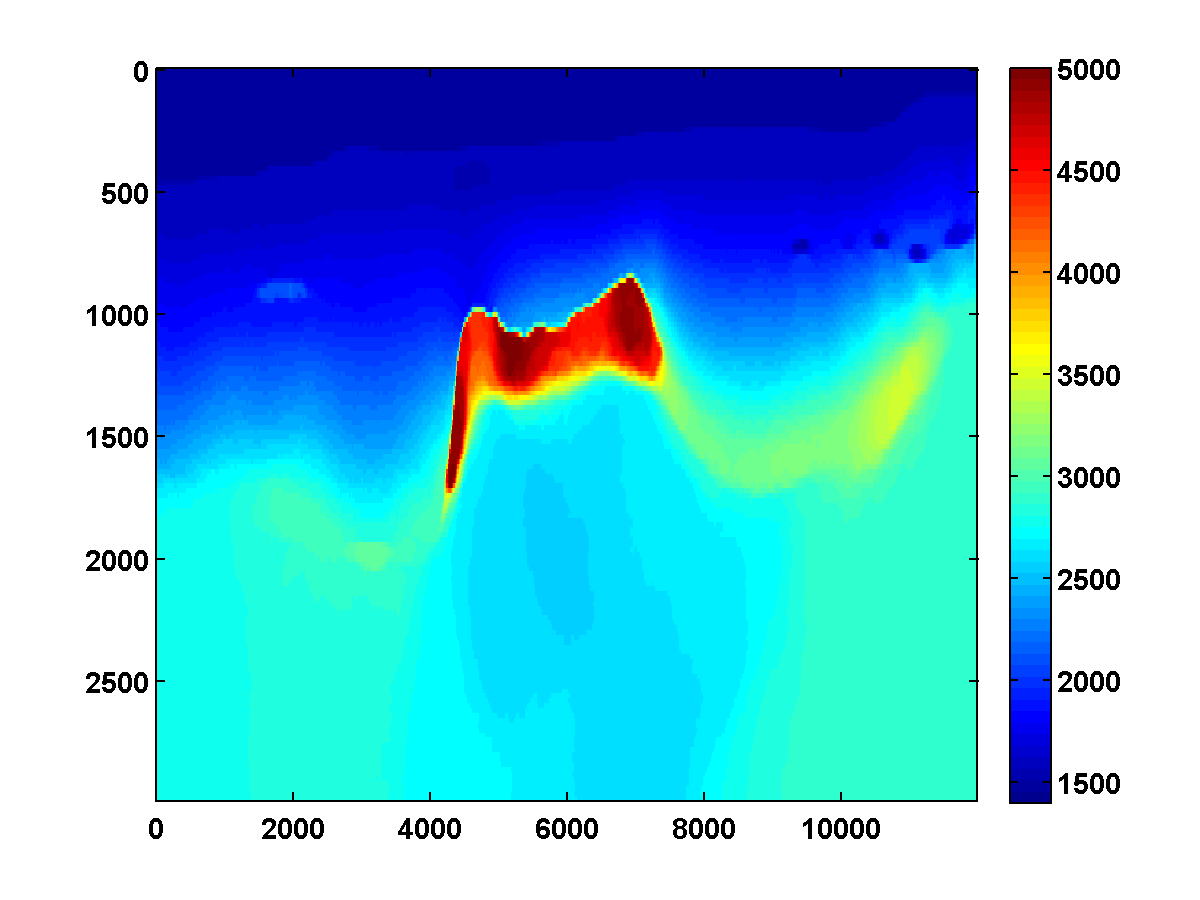}}
\caption{Recovered velocity without TV constraints from a poor initial model after one pass with
$\tau = .75 \tau_{\text{true}}$(a), a second pass with
$\tau = .825 \tau_{\text{true}}$ (b) and a third pass with
$\tau = .9 \tau_{\text{true}}$.}\label{fig:TVBPTM}
\end{figure}

With a TV constraint added, the method still tends to get stuck at a poor solution, even with multiple passes and different choices of $\tau$. Figure~\ref{fig:TVBPTM} shows the estimated velocity models after three passes, where increasing values of $\tau$ were used so that the TV constraint was weakened slightly after each pass. Inclusion of the TV constraint leads to accurate recovery of top salt, which a significant improvement. 
However,  the imprint of the velocity low is still too strong, 
and the results do not improve after three passes.

\section{Asymmetric TV-norm Constraints}\label{one-sided-tv-constraint}

Since velocity generally increases with depth in sedimentary basins while remaining more or less constant in salt, it is natural to penalize downward jumps in velocity. This can be done with a one-sided, asymmetric total-variation constraint that penalizes increases in the slowness squared in the depth direction. Such a constraint naturally fits in the scaled projected gradient framework, and can be imposed along with TV-norm and box constraints. Define a forward difference operator $D_z$ that acts in the depth direction so that $D_z m$ is a vector of differences of the form $\frac{1}{h}(m_{k+1,l} - m_{k,l})$ for $k=1,\dots, M_z$ and $l=1, \dots,  M_x$. To penalize the sum of the positive differences in $m$, we include the constraint
\begin{equation}
\|\max(0,D_z m)\|_1 \leq \xi \ ,
\label{HL}
\end{equation}
where the $\max$ operation is understood in a componentwise sense so that
\[
\|\max(0,D_z m)\|_1 = \sum_{k,l} \max(0,\frac{1}{h}(m_{k+1,l} - m_{k,l})).
\]
This constraint, which limits the size of the asymmetric TV-norm ball to a value $\xi$, is the hinge loss penalty 
applied to $D_z m$. The hinge loss is commonly used in machine learning for support vector machines
and support vector regression.

The constraint in (\ref{HL}) does not penalize model discontinuities in the horizontal direction, only in the depth direction. It is therefore likely to lead to vertical artifacts unless combined with additional constraints that penalize variations in the horizontal direction. We therefore combine this hinge-loss constraint with a TV constraint.

We can apply box constraints, TV-norm constraints, and one-sided TV constraints using 
the same scaled projected gradient framework. 
The set $C$ in~\eqref{Fmincon} is now the intersection of three constraints:
\begin{equation}
\min_m f(m) \ \  \text{subject to } \ \ m_i \in [b_i , B_i] \text{ , } \|m\|_{TV} \leq \tau \text{ and } \|\max\left(0,D_z m\right)\|_1 \leq \xi,
\label{QP3}
\end{equation}
and the convex subproblem~\eqref{dmcon} now has an additional one-sided TV constraint: 
\begin{equation}
\begin{aligned}
\dm & = \argmin_{\dm} \dm^\top  \nabla f(m^n) + \frac{1}{2}\dm^\top  (H^n + c_n I) \dm  \\
& \text{subject to } m^n_i + \dm_i \in [b_i , B_i] \text{ , } \|m^n + \dm\|_{TV} \leq \tau \\
& \text{ and } \|\max(0,D_z (m^n + \dm))\|_1 \leq \xi \ .
\end{aligned}
\label{pTVb3}
\end{equation}
As before, we can use the PDHG algorithm. Analogous to (\ref{Lagrangian}), we want to find a saddle point of the Lagrangian
\begin{equation}
\begin{aligned}
\L(\dm,p_1,p_2) & = \dm^\top  \nabla f(m^n) + \frac{1}{2}\dm^\top  (H^n + c_n I) \dm + g_B(m^n + \dm) \\
& + p_1^\top  D(m^n + \dm) - \tau \|p_1\|_{\infty,2} \\
& + p_2^\top  D_z(m^n + \dm) - \xi \max(p_2) - g_{\geq 0}(p_2) \ ,
\end{aligned}
\label{Lagrangian3}
\end{equation}
where $g_{\geq 0}$ denotes an indicator function defined by
\begin{equation*}
g_{\geq 0}(p_2) = \begin{cases} 0 & \quad \text{if} \quad p_2 \geq 0 \\
\infty & \quad \text{otherwise} \end{cases}.
\end{equation*}
The additional terms in this Lagrangian come from
the conjugate representation of the hinge loss constraint: 
\begin{equation}
\sup_{p_2} p_2^\top  D_z (m^n + \dm) - \xi \max(p_2) - g_{\geq 0}(p_2). 
\label{TV1sp}
\end{equation}
The reader can check that maximizing over $p_2$ we recover the indicator of the convex set
defined by the 1-sided TV constraint:  
\begin{equation*}
\begin{cases}
0 & \quad\text{if} \quad\|\max(0,D_z(m^n + \dm)\|_1 \leq \xi \\
\infty & \quad \text{otherwise.}
\end{cases}
\end{equation*}

The modified PDHG iterations are similar to~\eqref{PDHGexplicit} and given below. 

\begin{equation}
\begin{aligned}
p_1^{k+1} & = p_1^k + \delta D(m^n + \dm^k) - \Pi_{\|\cdot\|_{1,2} \leq \tau \delta}(p_1^k + \delta D(m^n + \dm^k)) \\
p_2^{k+1} & = p_2^k + \delta D_z(m^n + \dm^k) - \Pi_{\|\max(0,\cdot)\|_1 \leq \xi \delta}(p_2^k + \delta D_z(m^n + \dm^k)) \\
\dm^{k+1}_i & =  \max\left( b_i - m^n_i, \min\left( B_i - m^n_i, \widetilde{\dm}_i \right) \right)
\end{aligned}
\label{PDHGexplicit3}
\end{equation}
with 
\[
\widetilde{\dm}=\left(H^n + \left(c_n+\frac{1}{\alpha}\right)\I\right)^{-1}\left(-\nabla
f(m^n) + \frac{\dm^k}{\alpha} - D^\top (2p_1^{k+1}-p_1^k) -
D_z^\top (2p_2^{k+1}-p_2^k)\right).
\] 
The projection $\Pi_{\|\max(0,\cdot)\|_1
\leq \xi \delta}(z)$ is computed by projecting the positive part of $z$,
$\max(0,z)$, onto the simplex defined by $\{z:\, z_k \geq 0 \ , \ \sum_k z_k =
\xi \delta\}$.

\subsection*{Including the asymmetric TV-norm constraint}

\begin{figure}[t!]
\centering
\subfloat[\label{fig:BPTM_vinit}]{\includegraphics[width=0.330\hsize, natwidth = 650 ,natheight=642]{./images/initial_velocity_BPTM.png}}
\subfloat[\label{fig:TVHLBPTMdpw_a}]{\includegraphics[width=0.330\hsize, natwidth = 650 ,natheight=642]{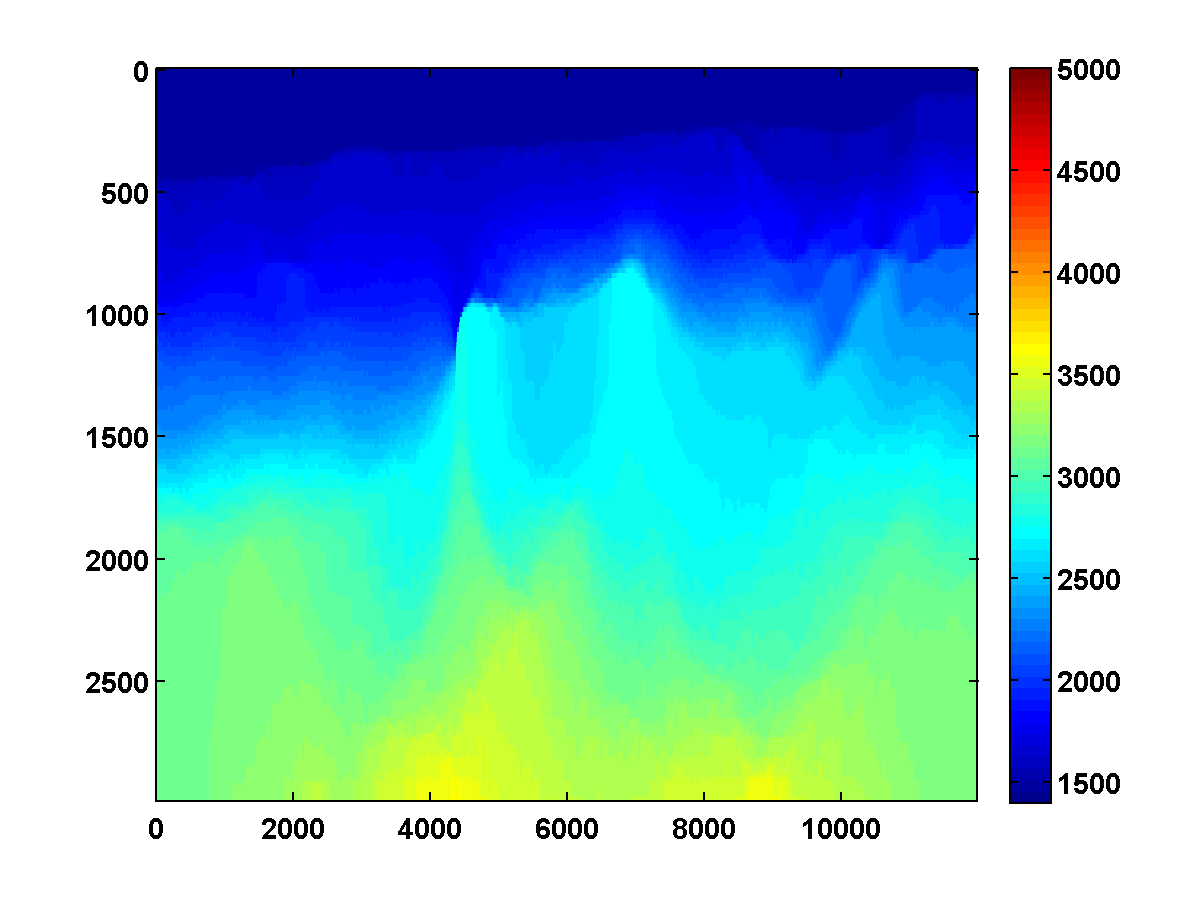}}
\subfloat[\label{fig:TVHLBPTMdpw_b}]{\includegraphics[width=0.330\hsize, natwidth = 650 ,natheight=642]{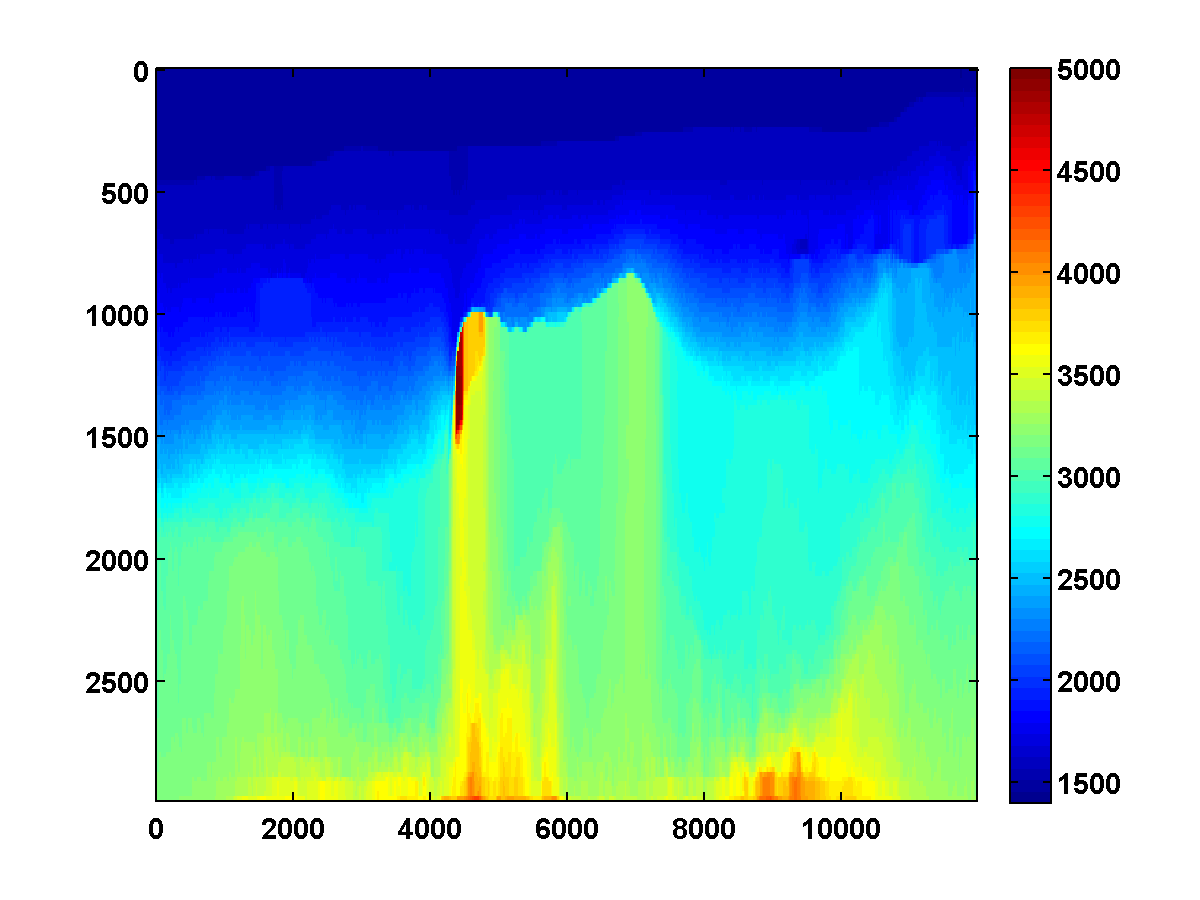}}
\\
\subfloat[\label{fig:TVHLBPTMdpw_c}]{\includegraphics[width=0.330\hsize, natwidth = 650 ,natheight=642]{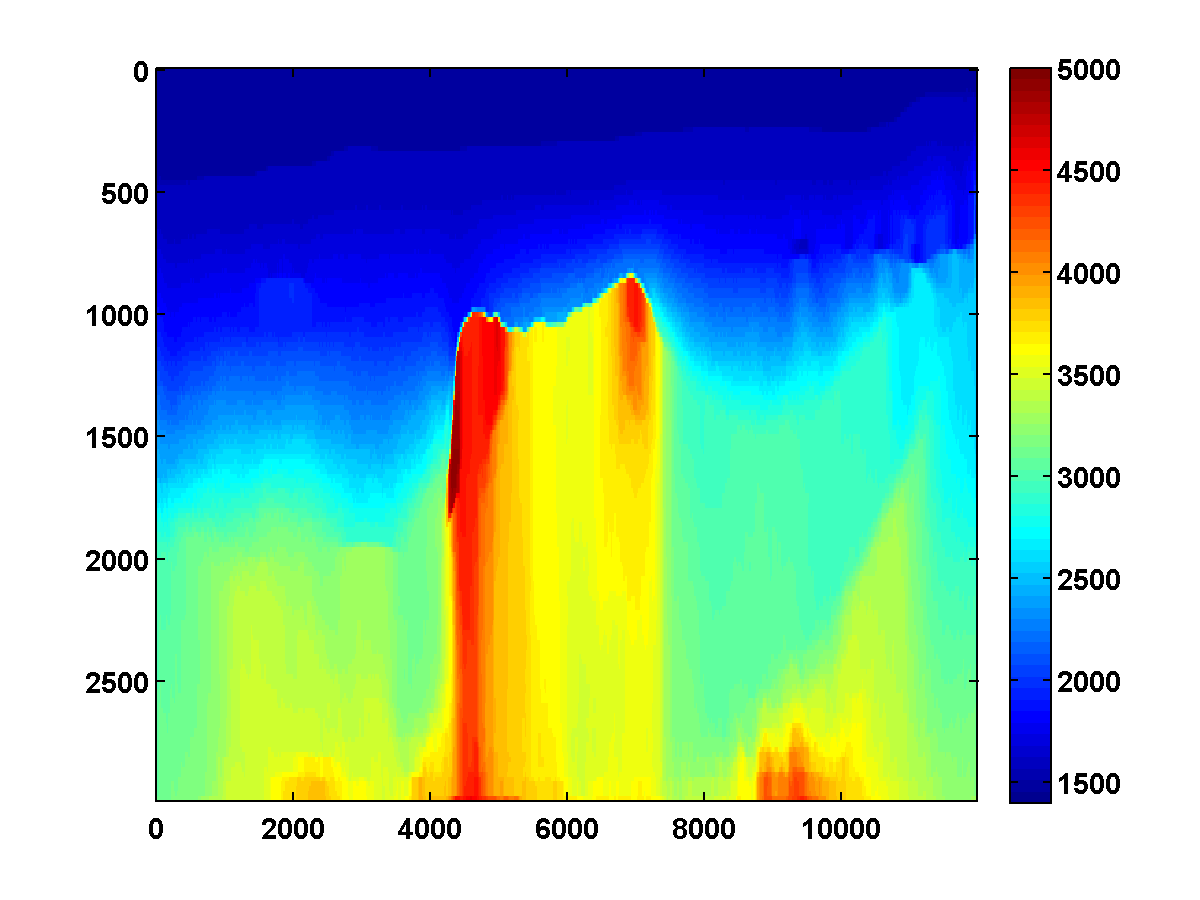}}
\subfloat[\label{fig:TVHLBPTMdpw_d}]{\includegraphics[width=0.330\hsize, natwidth = 650 ,natheight=642]{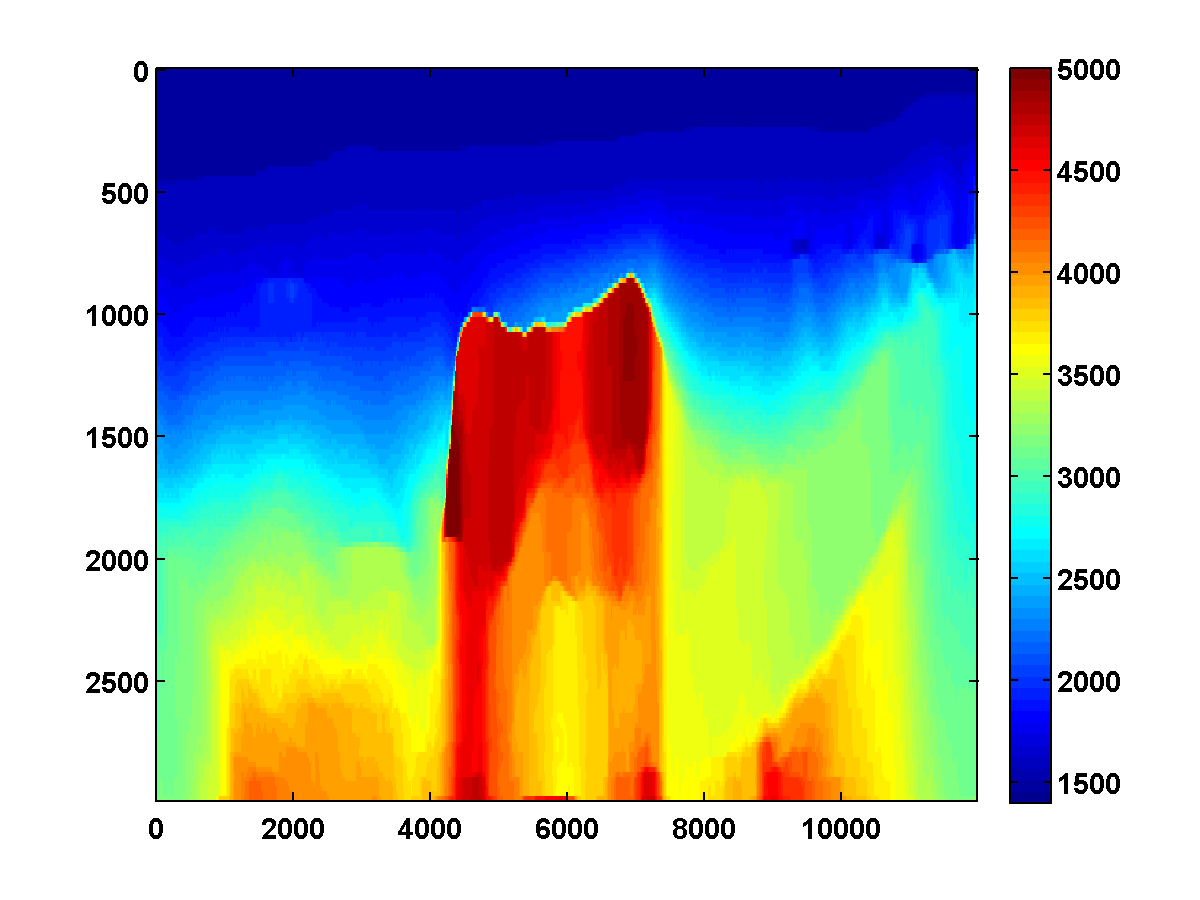}}
\subfloat[\label{fig:TVHLBPTMdpw_e}]{\includegraphics[width=0.330\hsize, natwidth = 650 ,natheight=642]{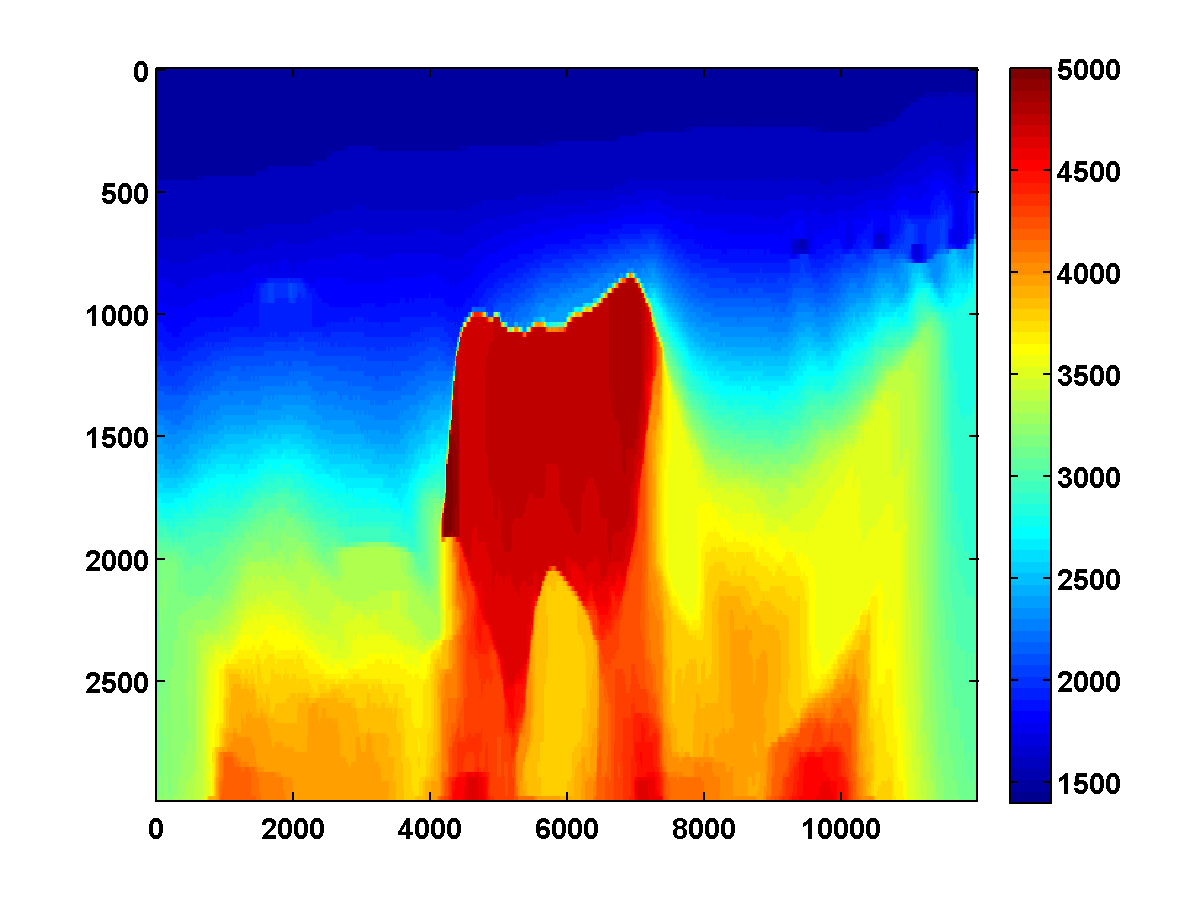}}
\\
\subfloat[\label{fig:TVHLBPTMdpw_f}]{\includegraphics[width=0.330\hsize, natwidth = 650 ,natheight=642]{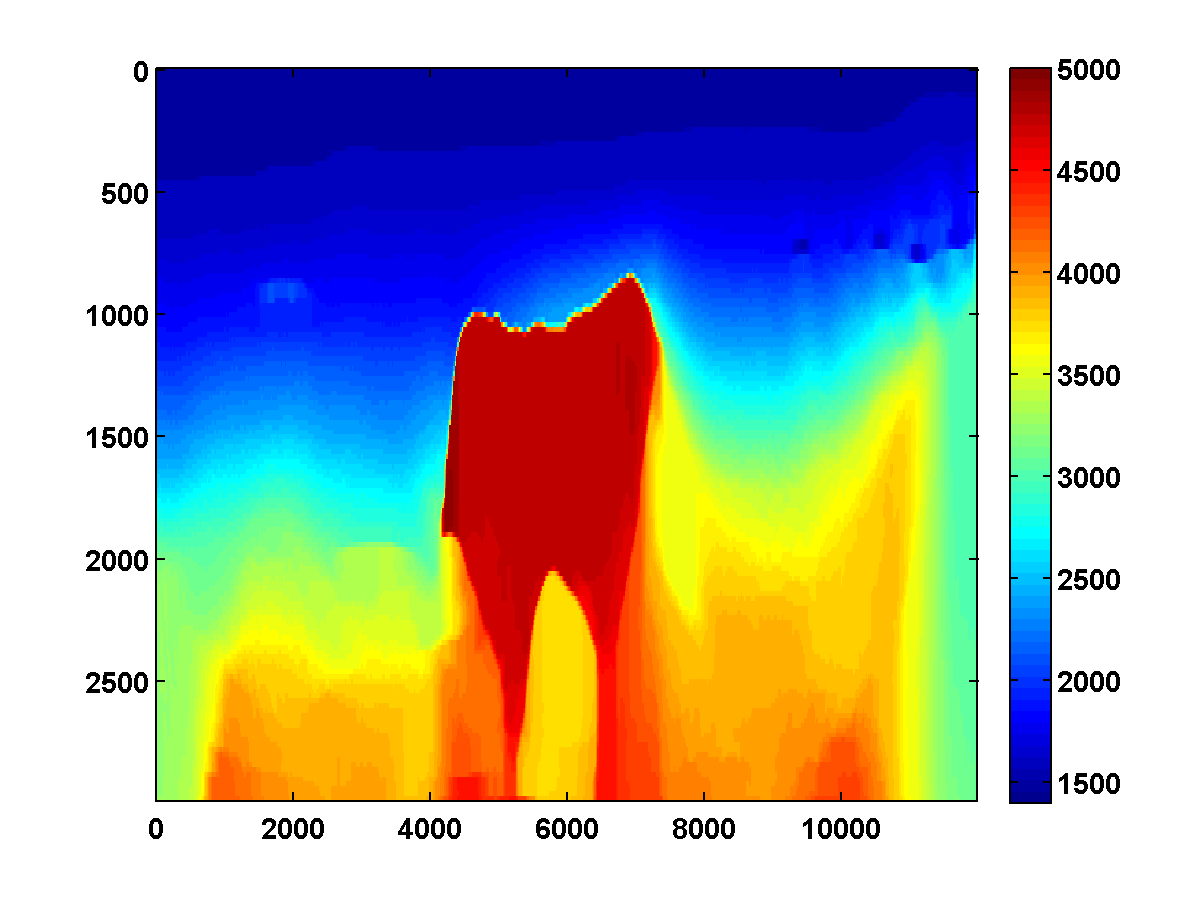}}
\subfloat[\label{fig:TVHLBPTMdpw_g}]{\includegraphics[width=0.330\hsize, natwidth = 650 ,natheight=642]{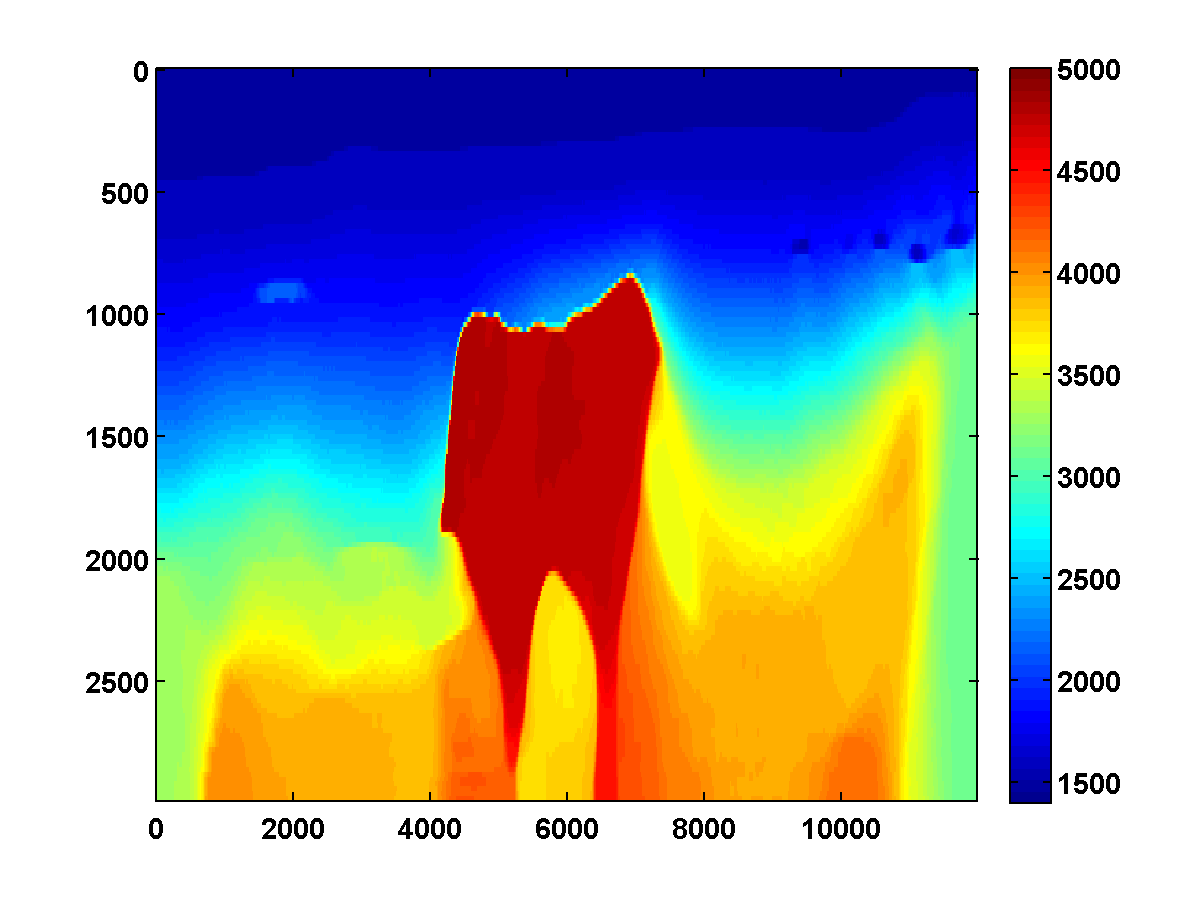}}
\subfloat[\label{fig:TVHLBPTMdpw_h}]{\includegraphics[width=0.330\hsize, natwidth = 650 ,natheight=642]{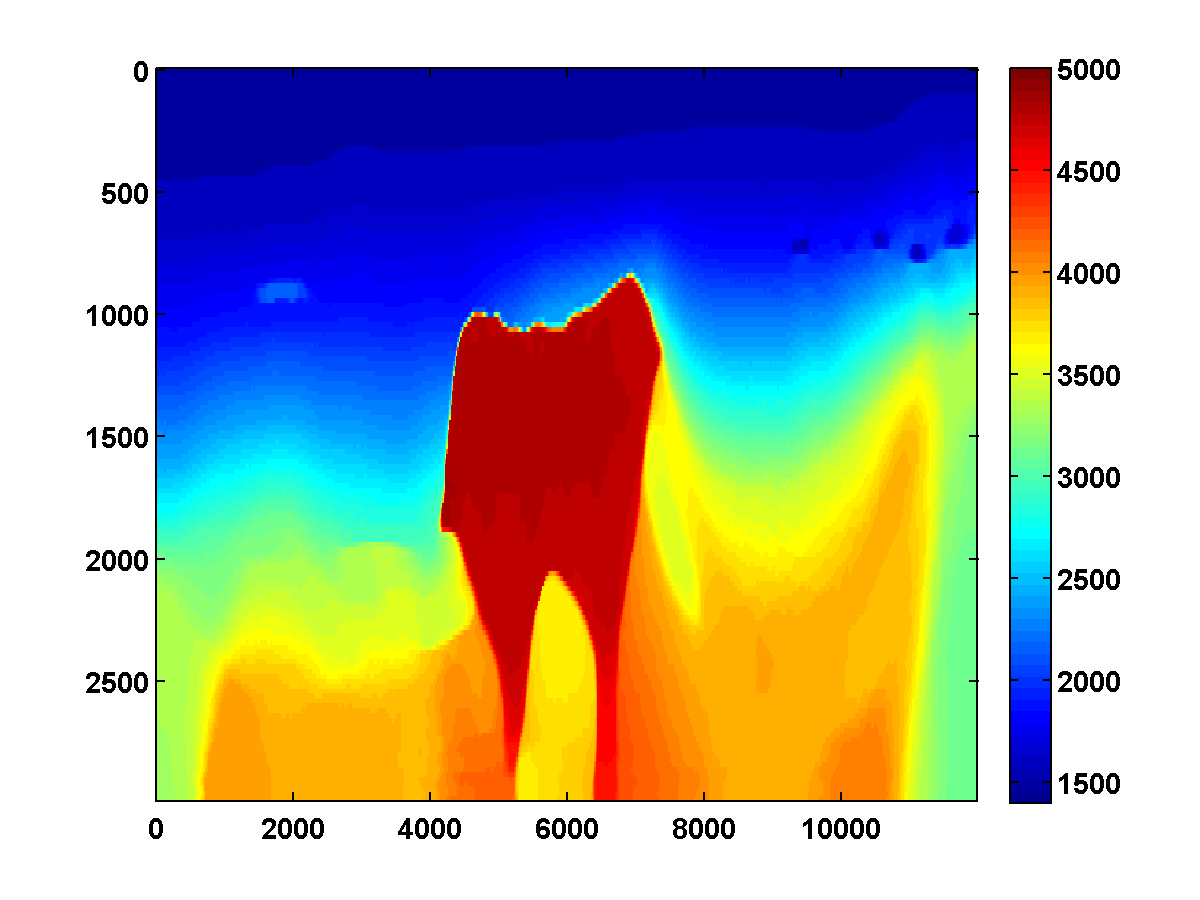}}
\caption{Initial velocity (a) and recovered velocity with asymmetric TV
continuation corresponding to $\frac{\xi}{\xi_{\text{true}}} = $ $.01$
(a), $.05$ (b), $.10$ (c), $.15$ (d), $.20$ (e), $.25$ (f), $.40$ (g),
$.90$ (h). Movies of the solutions with and without the asymmetric TV norm are available as ancillary material.}\label{fig:TVHLBPTMdpw}
\end{figure}

Using WRI, bounds, and TV-norm constraints, our inversion results  still suffered from missing low frequencies when starting from a poor initial model (see Figures~\ref{fig:BPTM} and \ref{fig:TVBPTM}). 
To discourage spurious downward jumps in velocity after entering the salt, we add the asymmetric TV constraint and use a continuation strategy in the $\xi$ parameter. 
We start with a small value for the 1-sided TV-norm level set, and gradually increase 
it over each successive pass through the frequency batches. 
This continuation approach encourages the initial velocity estimates to be nearly monotonically increasing in depth, a notion that corresponds to assuming we are in a sedimentary basin where we can step into the salt but not out of it. 
At later passes, the asymmetric TV constraint is relaxed to allow more downward jumps to fit the observed data. Starting with the poor initial model in Figure~\ref{fig:BPTM_vinit}, Figure~\ref{fig:TVHLBPTMdpw} shows the progression of velocity estimates over eight passes. The sequence of $\xi$ parameters as a fraction of $\|\max(0,D_z m_{\text{true}})\|_1$ is chosen to be $\{.01, .05, .10, .15, .20, .25, .40, .90\}$. We keep the $\tau$ parameter fixed at $.9 \tau_{\text{true}}$ throughout. Although small values of $\xi$ cause some vertical artifacts, the continuation strategy is surprisingly effective at preventing the method from getting stuck at a poor solution. As $\xi$ increases, the bottom of the salt is recovered.

Both the model error and data misfit continue to decrease during each pass, as shown in Figure~\ref{fig:WRI_datamod_error}.
The rugged shape of the objective function curve has two different scales that correspond to the relaxation passes and
the increasing frequency content. The former causes the sudden drops in the functional every 9 iterations,
which occurs because the data fit improves as we relax the TV constraints: due to the extra structure allowed
in the model as we relax the constraints, the data-misfit reduces as more events can be generated and matched.
The second, the objective function increase within each low to high pass through frequency batches is explained by the
increase in frequency: there is more energy present at higher frequencies in the data and therefore the data-misfit increases accordingly.

Even with a poor initial model, the asymmetric TV constraint with
continuation is able to recover the main features of the ground truth
model. Poor recovery near the left and right boundaries is expected
because the sources along the surface start about $1000\,\mathrm{m}$ away from these
boundaries, which affects the illumination.

\subsection*{Conventional FWI approach}\label{constrained-adjoint-state-method-comparisons}

So far, all examples we discusses were based on WRI~\eqref{eq:WRI}.
Recall that unconstrained FWI was unable to produce meaningful results when given a poor starting model such as the one plotted in Figure~\ref{fig:BPTM_vtrue}. However, if we impose the combination of the box, TV and asymmetric constraints in Algorithm~\ref{SGPalg}, the conventional reduced adjoint-state method (\textrm{cf.}~\eqref{eq:FWI}) produces excellent results as long as we replace $H^n + c_n\I$ $c_n(H^n + \nu\I)$ for some small positive $\nu$ and where $H^n$ is defined by \eqref{eq:HFWI}.

\begin{figure}
\centering
\subfloat[\label{fig:BPTM_vinit}]{\includegraphics[width=0.330\hsize, natwidth = 650 ,natheight=642]{./images/initial_velocity_BPTM.png}}
\subfloat[\label{fig:ASTVHLBPTMdpw_a}]{\includegraphics[width=0.330\hsize, natwidth = 650 ,natheight=642]{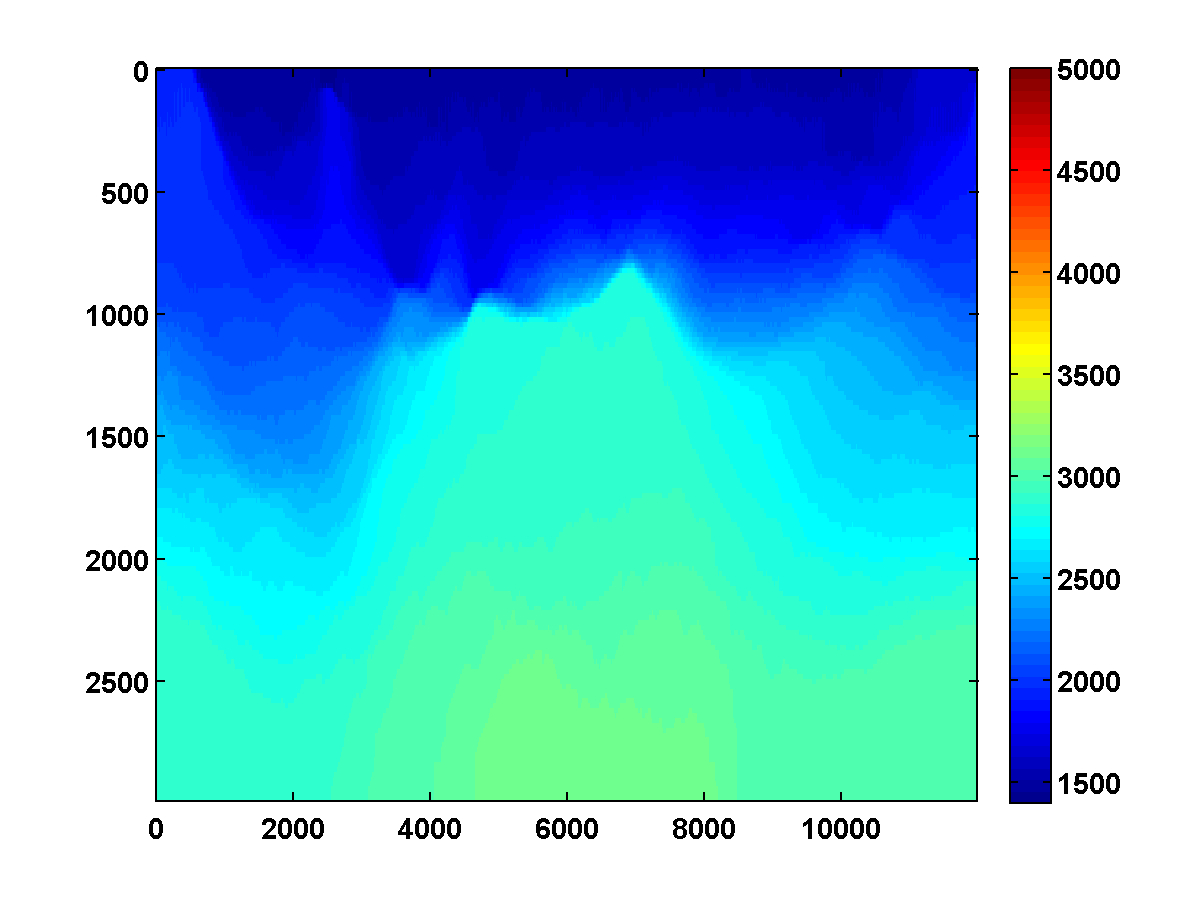}}
\subfloat[\label{fig:ASTVHLBPTMdpw_b}]{\includegraphics[width=0.330\hsize, natwidth = 650 ,natheight=642]{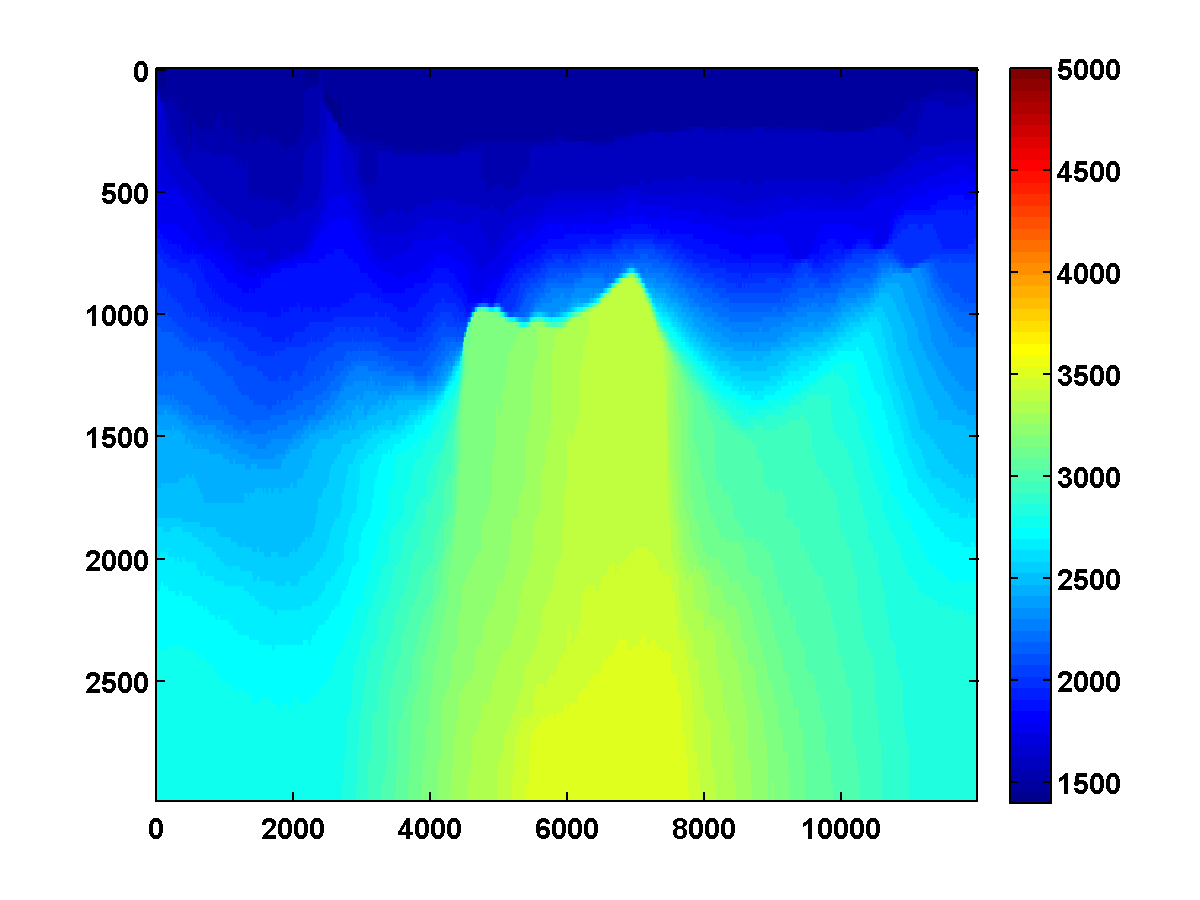}}
\\
\subfloat[\label{fig:ASTVHLBPTMdpw_c}]{\includegraphics[width=0.330\hsize, natwidth = 650 ,natheight=642]{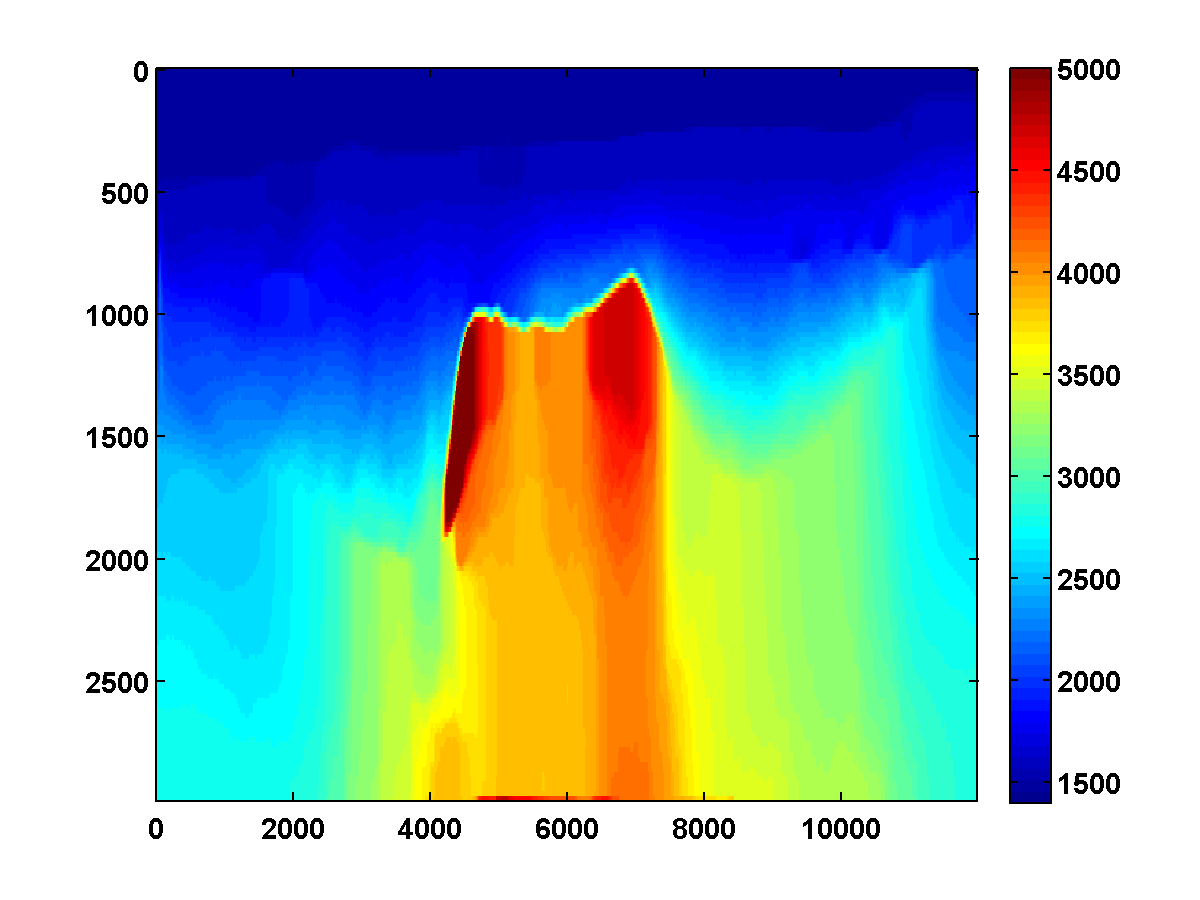}}
\subfloat[\label{fig:ASTVHLBPTMdpw_d}]{\includegraphics[width=0.330\hsize, natwidth = 650 ,natheight=642]{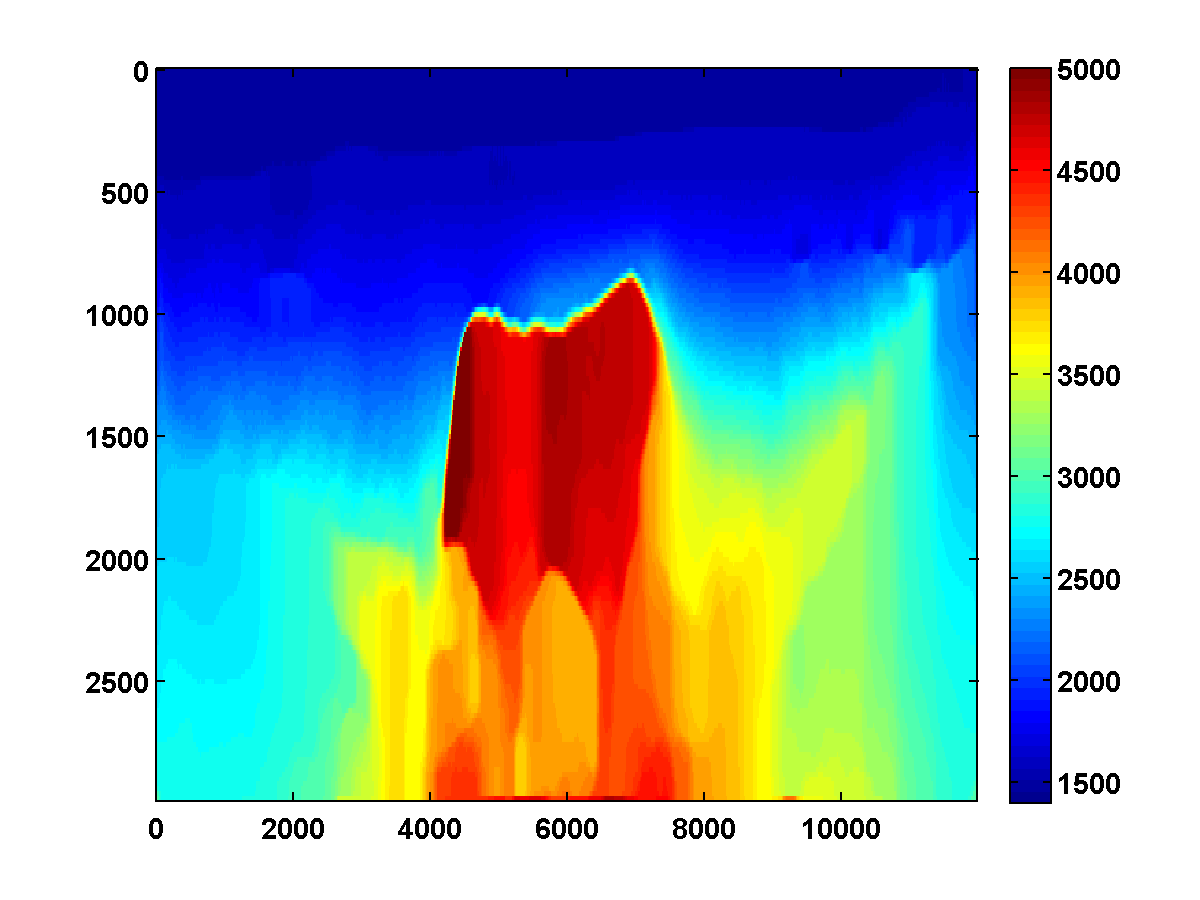}}
\subfloat[\label{fig:ASTVHLBPTMdpw_e}]{\includegraphics[width=0.330\hsize, natwidth = 650 ,natheight=642]{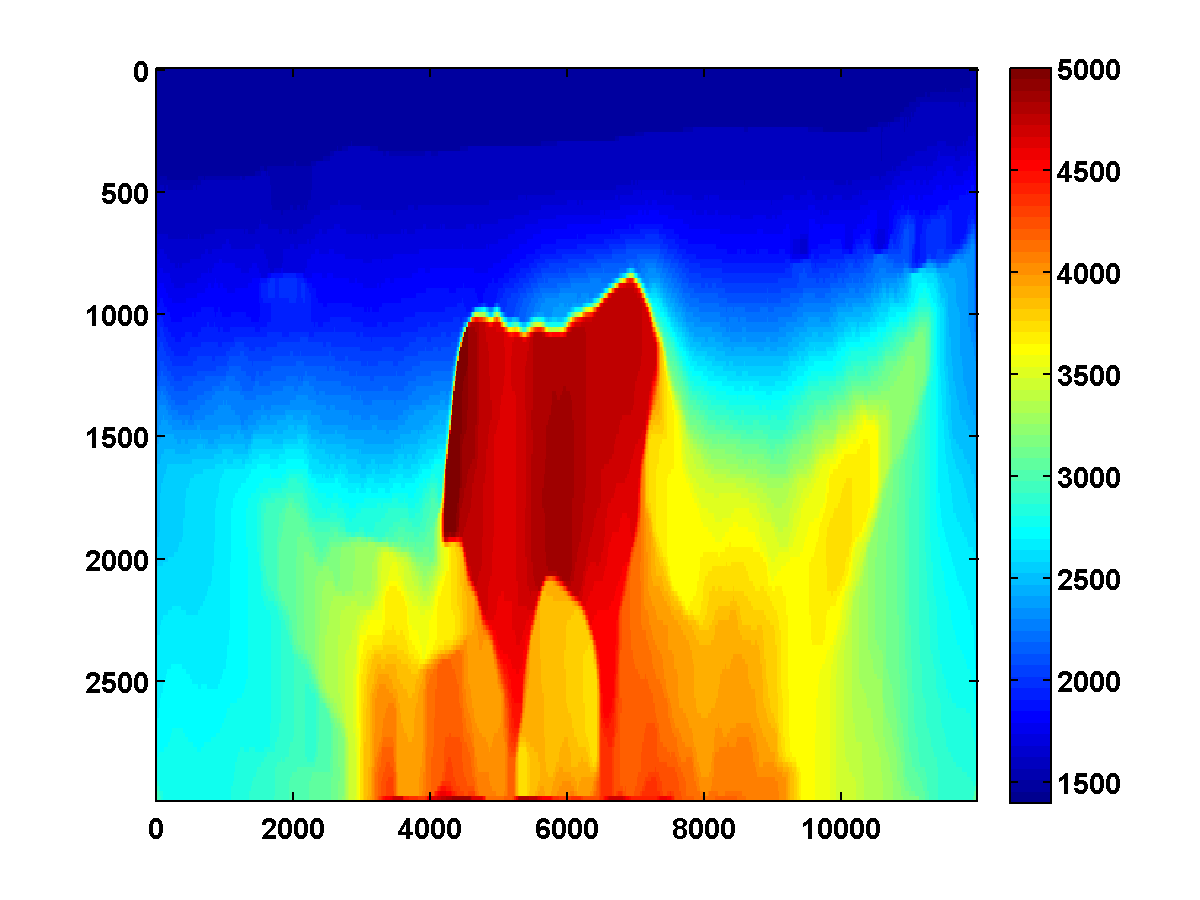}}
\\
\subfloat[\label{fig:ASTVHLBPTMdpw_f}]{\includegraphics[width=0.330\hsize, natwidth = 650 ,natheight=642]{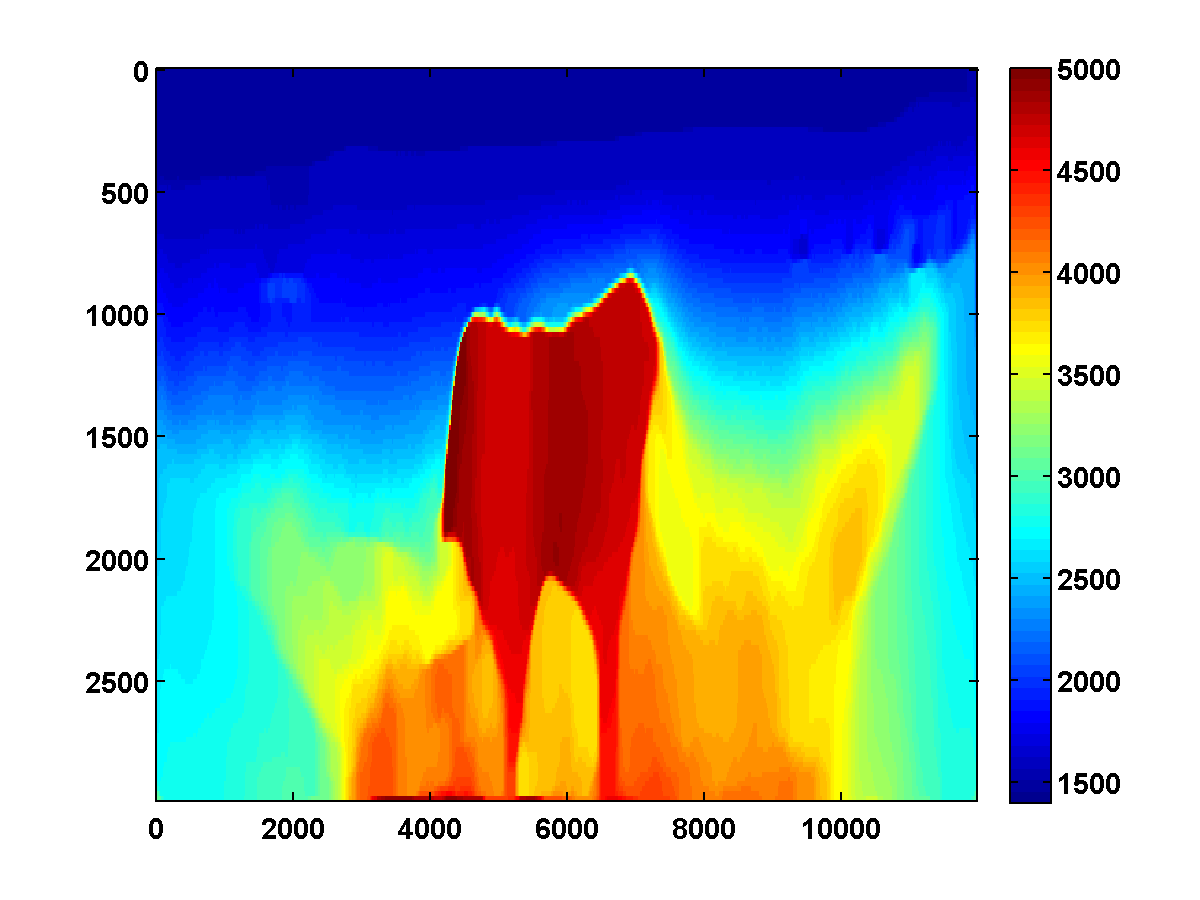}}
\subfloat[\label{fig:ASTVHLBPTMdpw_g}]{\includegraphics[width=0.330\hsize, natwidth = 650 ,natheight=642]{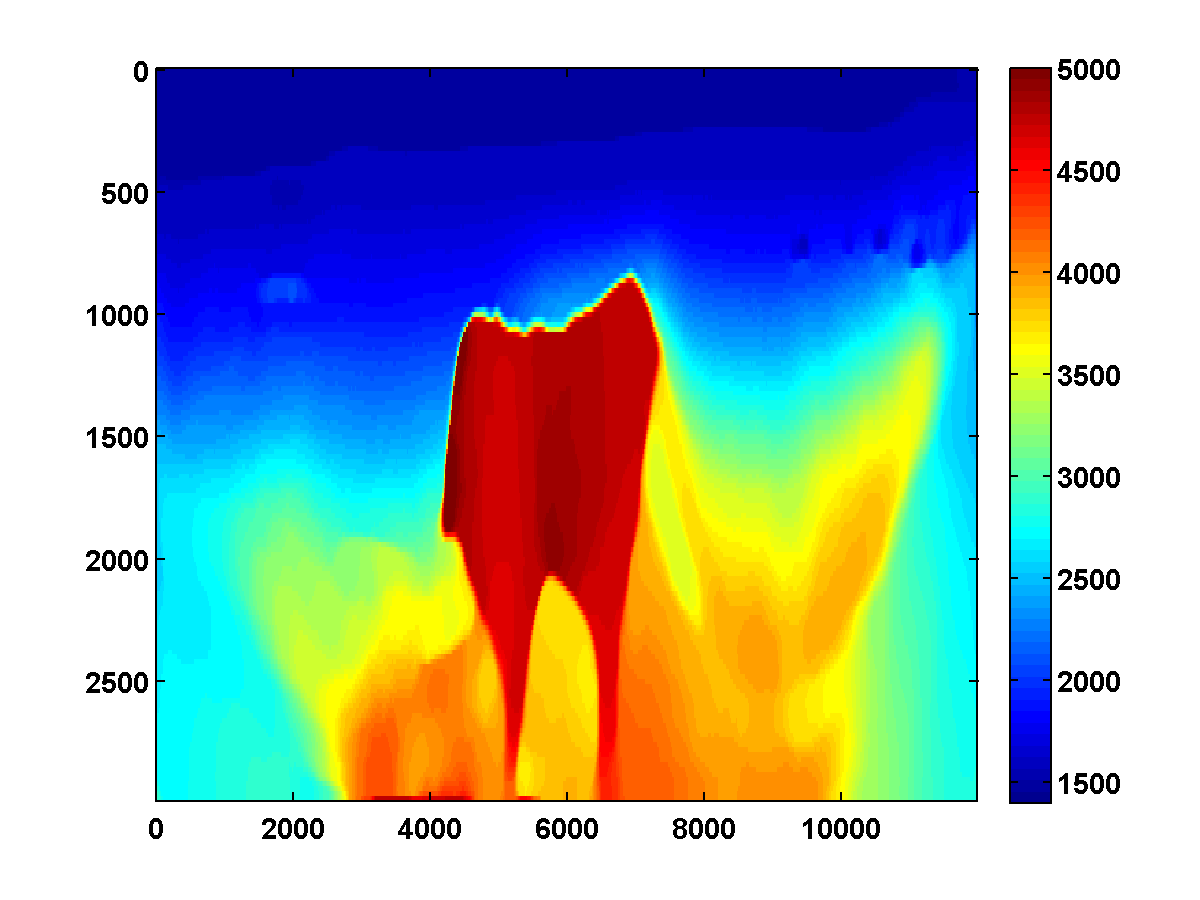}}
\subfloat[\label{fig:ASTVHLBPTMdpw_h}]{\includegraphics[width=0.330\hsize, natwidth = 650 ,natheight=642]{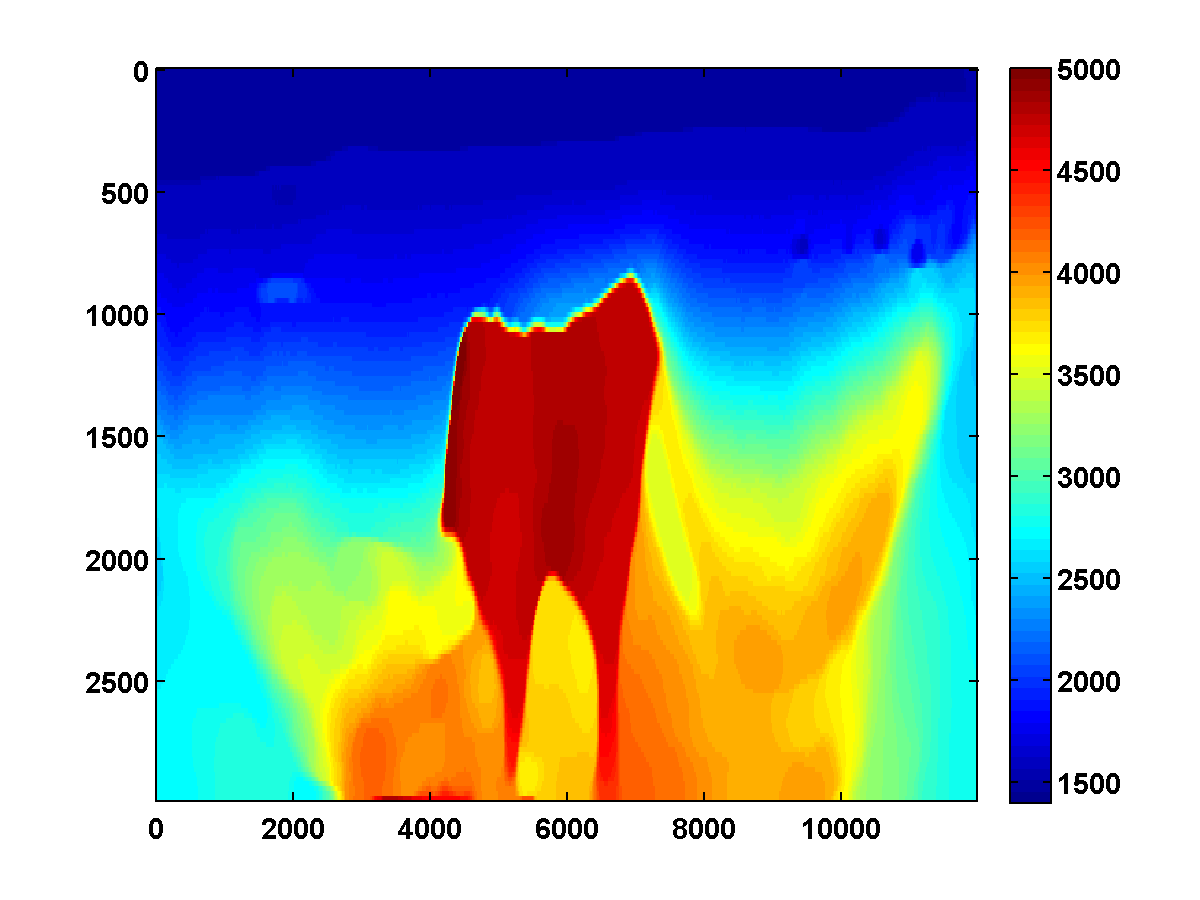}}
\caption{Initial velocity (a) and recovered velocity for the constrained
adjoint state method with one-sided TV continuation corresponding to
$\frac{\xi}{\xi_{\text{true}}} = $ $.01$ (a), $.05$ (b), $.10$ (c),
$.15$ (d), $.20$ (e), $.25$ (f), $.40$ (g), $.90$
(h).}\label{fig:ASTVHLBPTMdpw}
\end{figure}

\begin{figure}
\centering
\subfloat[\label{fig:WRI_data_error}]{\includegraphics[width=0.480\hsize]{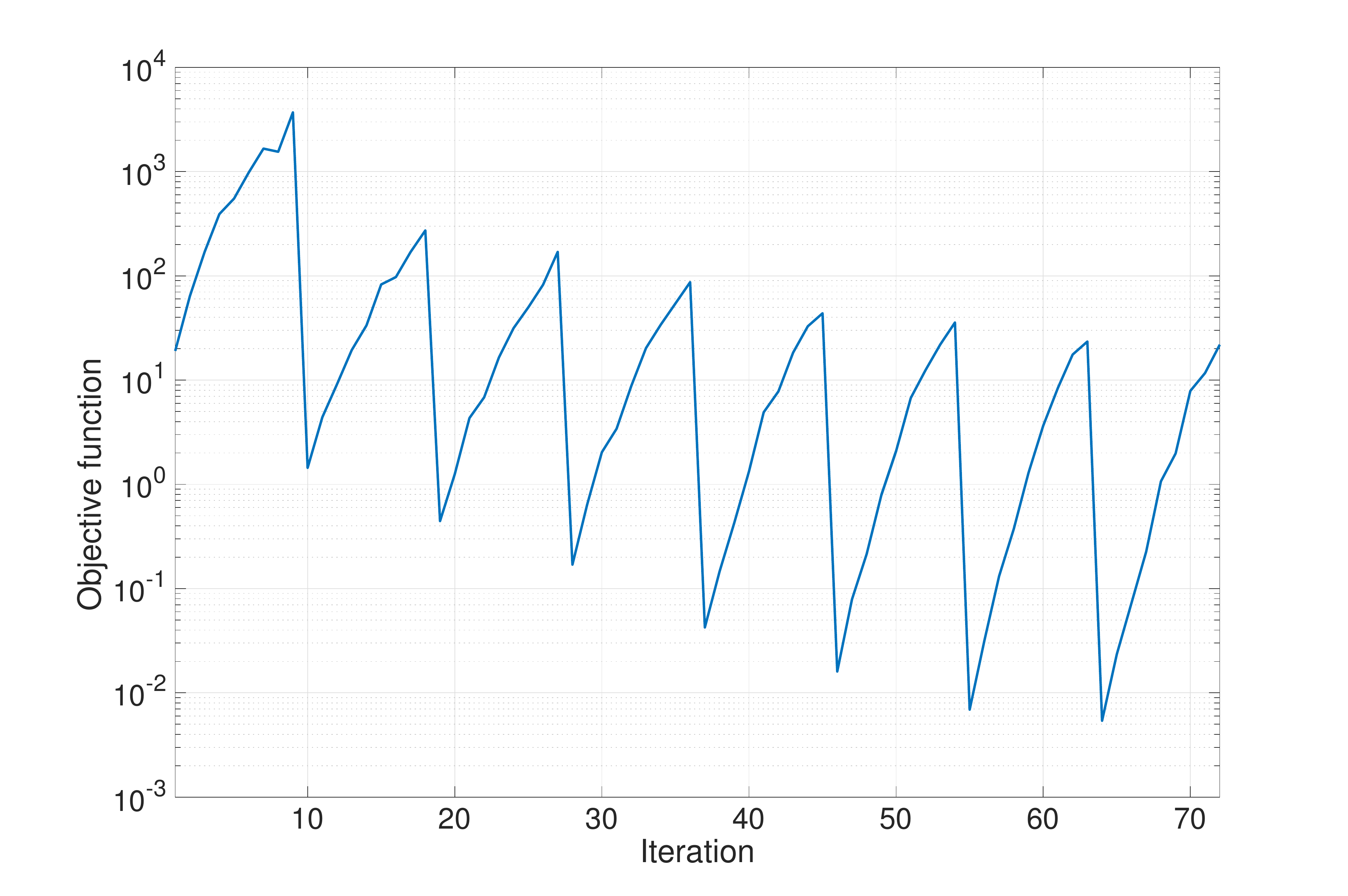}}
\subfloat[\label{fig:WRI_mod_error}]{\includegraphics[width=0.480\hsize]{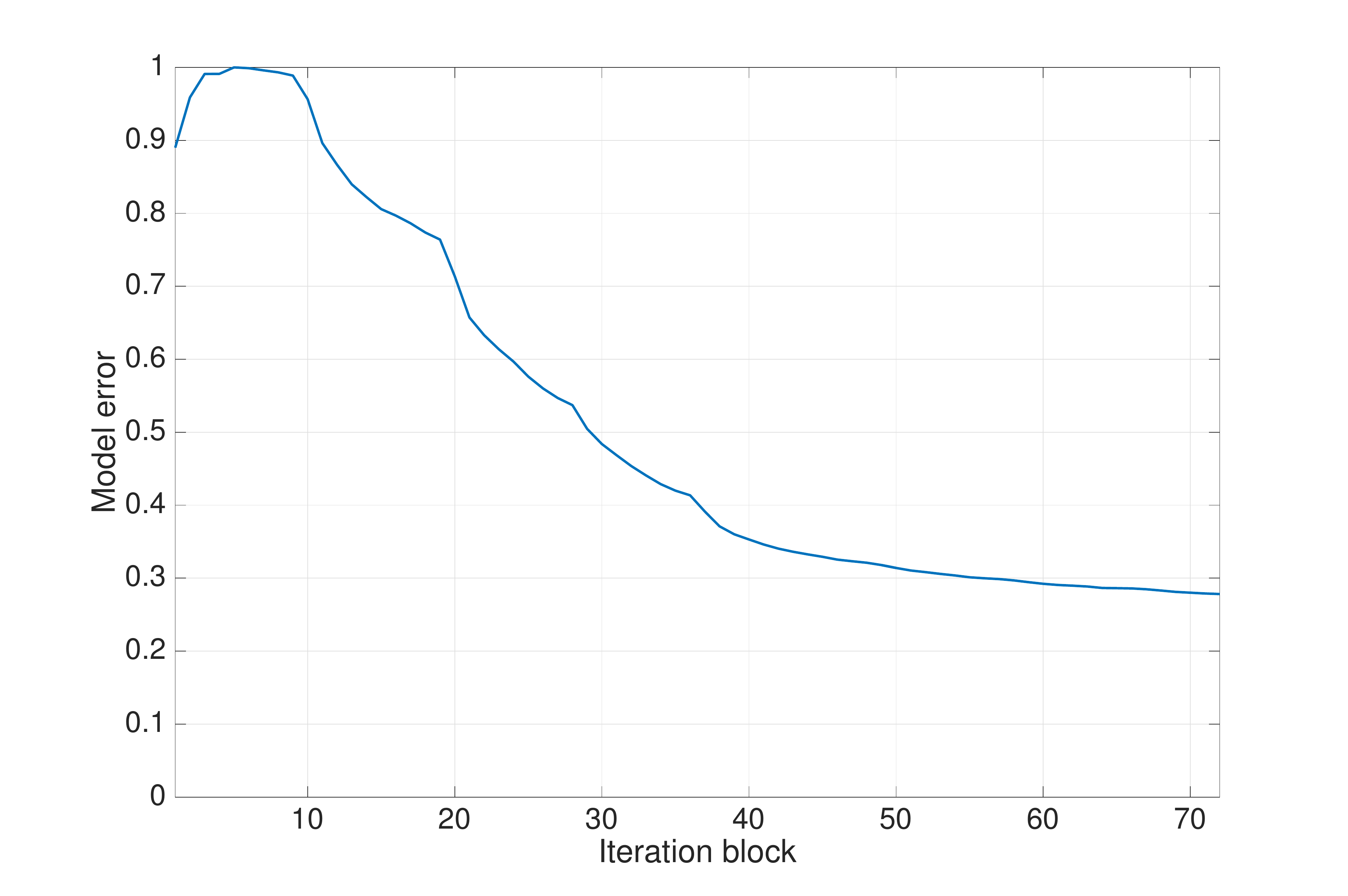}}
\caption{WRI objective function (data misfit) evolution (a). The model error ---normalised RMS of the difference between the true
and the current iterate--- (b) increases initially, due to the effect the strict
constraints coupled with the WRI updates have in the
early iterations, but after the values are relaxed, we observe a descending trend that reduces the model
misfit below 30\%. The 72 iterations correspond to the cumulate effect of all the iterations for 
each frequency.}
\label{fig:WRI_datamod_error}
\end{figure}

\begin{figure}
\centering
\subfloat[\label{fig:FWI_data_error}]{\includegraphics[width=0.480\hsize]{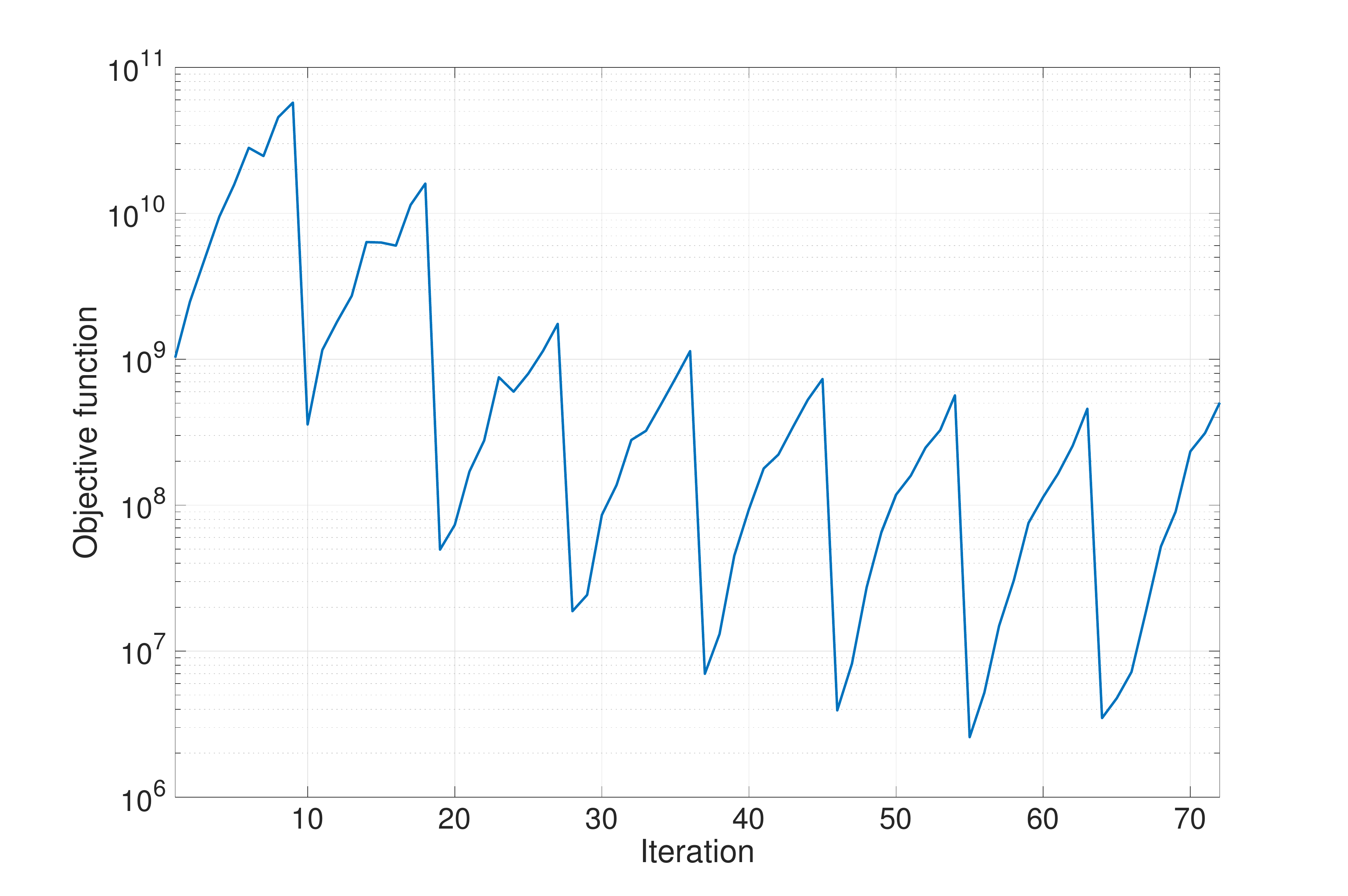}}
\subfloat[\label{fig:FWI_mod_error}]{\includegraphics[width=0.480\hsize]{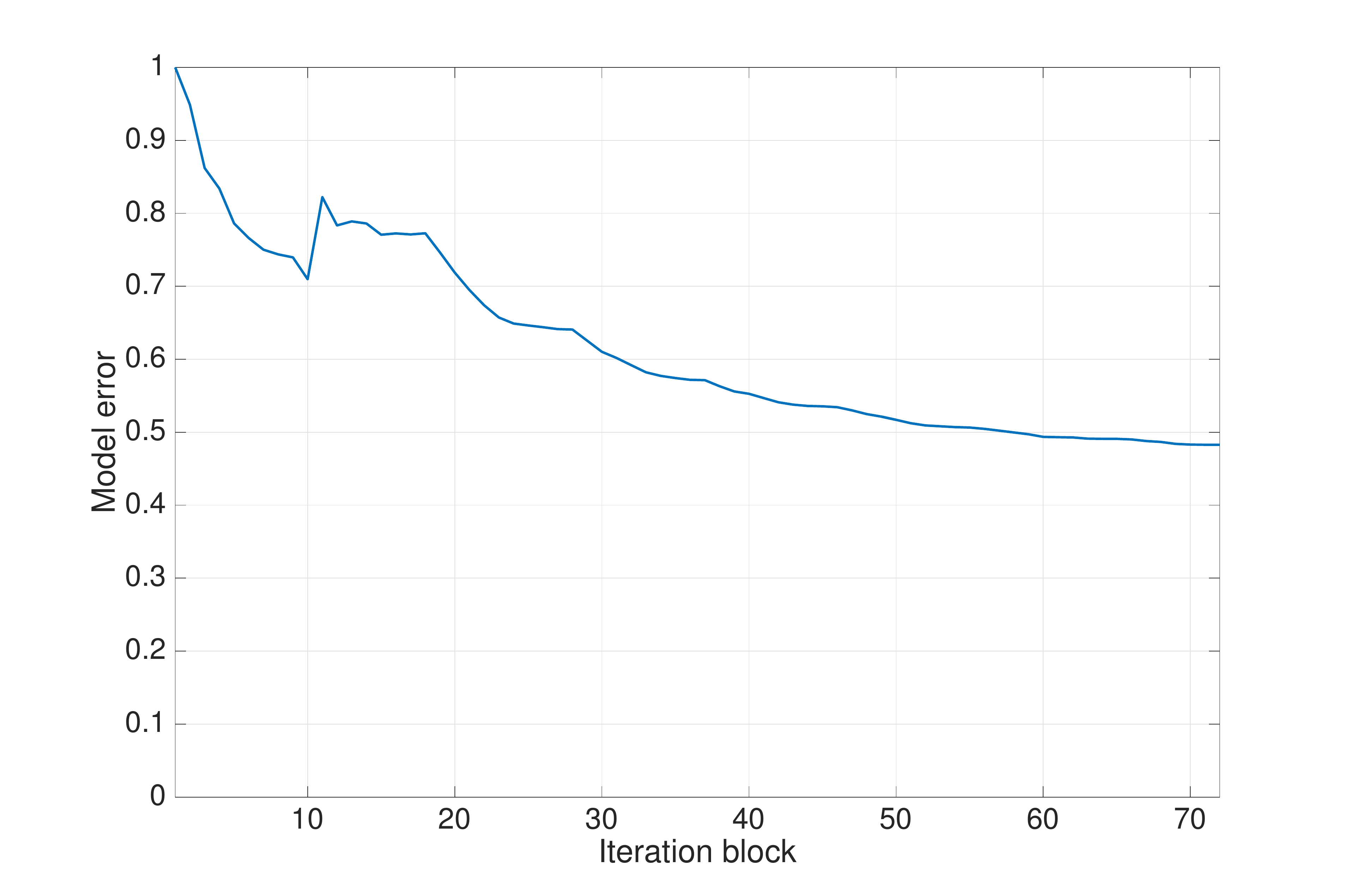}}
\caption{FWI objective function (data misfit) evolution (a). In this case the model error  ---normalised RMS of the difference between the true
and the current iterate--- (b) behaves as expected initially but shows an increase after the first pass;
due to the limitations of FWI (compared to WRI), the model error does not decrease as much as in the
WRI experiment and stalls below 50\%.
Again, the 72 iterations correspond to the cumulate effect of all the iterations for 
each frequency.}
\label{fig:FWI_datamod_error}
\end{figure}

We use the same continuation strategy as was used to generate the results in Figure~\ref{fig:TVHLBPTMdpw}, and the results are nearly as good (see Figure~\ref{fig:ASTVHLBPTMdpw}). Compared to the WRI method, the results are visually slightly worse near the top and sides of the model. Additionally, WRI finds 
a significantly better model error relative to ground truth (compare Figure~\ref{fig:WRI_mod_error} with Figure~\ref{fig:FWI_mod_error}), but it is encouraging to see that once again the error continues to decrease during each pass instead of stagnating at a poor solution. 

In summary, the examples clearly demonstrate that continuation in the $\xi$ parameter for the asymmetric TV constraint appears to be a promising strategy for preventing both the constrained WRI and adjoint-state FWI  from stagnating in a poor solution when starting from a bad starting model.

\section{Discussion}
\label{final-discussion}

 \textit{The work presented so far mainly concerned contributions from the first author Ernie Esser, who passed away preparing his work for publication. Our edits and additions were mostly stylistic and very much done in the spirit of Ernie's original draft\footnote[1]{Available at \href{https://www.slim.eos.ubc.ca/content/total-variation-regularization-strategies-full-waveform-inversion-improving-robustness-noise}{https://www.slim.eos.ubc.ca/content/total-variation-regularization-strategies-full-waveform-inversion-improving-robustness-noise}.}, which contained a section on "ongoing work". We took the liberty to rewrite and annotate this section to reflect 
 the impact of Ernie's contributions to seismic inversion for geologically challenging settings.
 Below, we give a summary of Ernie's main points, along with brief descriptions of 
 research progress Ernie's work inspired.  
 }

\subsection*{Ernie's points}

While the numerical experiments clearly demonstrate the potential benefits of including bound, TV- and asymmetric TV-norm constraints in wave-equation based inversions, the following issues remain:

\begin{itemize}
   
   \item\textbf{Inversion crime.} The inversion results were all obtained from
   data generated with the same forward modelling kernel as the kernel used to
   invert the data. This practice, known as the {\it inversion crime}, is common
   when developing new methodologies. Ernie was aware of this and
   wrote \textit{``The numerical examples should be recomputed without inversion
   crime.''} With the help of others, Ernie's co-authors have redone the
   experiments and we are happy to reports that Ernie's proposed method
   continues to perform well. These results have been presented at the 2016
   EAGE meeting and will be included in a future publications.

   \item{\textbf{Boundary conditions.}} Ernie's original work was in his words
   based on \textit{``simplistic boundary conditions that may conspire to make the inversion problem easier.''} Ernie refers here to the fact that the
   inversions may inadvertently been helped by unnatural reflections emanating
   from the boundaries. These reflected waves could have illuminated bottom
   salt and the salt flanks. Again with assistance from others, we have been
   able to redo the experiments with a more sophisticated modeling including
   more accurate boundary conditions. The conclusion from these experiments is
  that simplistic boundary conditions had little to no
   effect on the final inversion results. Following Ernie's suggestion to
   \textit{``generate the data using more
   sophisticated modeling and boundary conditions while continuing to use the
   simple variants for inversion,''} we found that the
   improvements in the inversion results still stand.
   
   \item{\textbf{Effect of randomness.}} To accelerate the computations, Ernie
   used a source encoding technique where the required number of
   wave-equation solves is reduced via randomized projections. While this
   approach reduces the number of sources to only two, it creates noisy
   crosstalk. Quoting Ernie, \textit{``To remove any effect of randomness, the
   examples should be recomputed using sequential shots instead of simultaneous
   shots with the random weights redrawn every model update. This will be more
   computationally expensive, but it's doable and is not expected to
   significantly alter the results.''} We followed up on this suggestion and we
   found that working with all sources significantly improved the
   results (See \cite{Esser2016CWI}).

   \item{\textbf{Continuation strategy.}} Ernie wrote \textit{``More practical
   methods for selecting the parameters and more principled ways of choosing
   continuation strategies are needed. Currently the $\tau$ and $\xi$
   parameters are chosen to be proportional to values corresponding to the true
   solution. Although the true solution clearly isn't known in advance, it may
   still be reasonable to base the parameter choices on estimates of its
   total-variation (or one-sided TV). It would be even better to develop
   continuation strategies that don't rely on any assumptions about the
   solution but that are still effective at regularizing early passes through
   the frequency batches to prevent the method from stagnating at poor
   solutions.''} This is still a topic of active research. However, our
   experience applying Ernie's heuristic of relaxing the constraint in
   combination with warm restarts suggests to us that it is critical to relax
   the constraints slowly enough so that the algorithm steers free from the
   effects of parasitic local minima.
   
   \item{\textbf{Convex subproblems.}} With regard to imposing constraints on
   the model iterations Ernie pointed out \textit{``There is a lot of room for
   improvement in the algorithm used to solve the convex subproblems. There are
   for example methods such as the one in \cite{Chambolle2011} that have better
   theoretical rates of convergence and are straightforward to implement.''}. We
   leave this suggestion. as well as extensions to 3D seismic models, for future work.
   
\end{itemize}

\subsection*{Geophysical impact}

Inverting for velocity in sedimentary basins interspersed with high-velocity sharp-contrast salt bodies is perhaps one of the most challenging inversion problems in geophysics. Despite numerous efforts from academia as well as from industry, little progress has been made using fully automatic (i.e. void of extensive human interaction) methods to solve this inversion problem in a systematic and reproducible way. Instead the community relies on intricate workflows, which combine reflection tomography, manual {\it salt flooding} (where the salt is continued downwards into the basin), and full-waveform inversion. Quality of results depend rely on geophysical experience and manual intervention, and inversions are costly not reproducible.

To our knowledge, Ernie's work is a first successful attempt to replace these labour intensive manual salt flooding workflows with an automatic heuristic: successive relaxations of the asymmetric TV-norm constraints with warm-started passes through frequency batches. In essence, Ernie encoded salt flooding \cite{Esser2016CWI} into a powerful automatic workflow, a breakthrough result.

Ernie's continuation approach is able to steer free from parasitic stationary points. Numerical experience with the approach highlights the following important considerations: 

\begin{itemize} 
	
   \item The lowest frequency should be sufficiently low to allow 
   progress during each pass through the frequencies.
   
   \item The successively relaxed constraint sets should be large enough to prevent stalling (i.e. to allow model updates), but small enough to help avoid parasitic stationary points.

\end{itemize}

\section{Conclusions and Future Work}\label{conclusions-and-future-work}

We presented a scaled gradient projection algorithm for minimizing WRI and FWI formulations subject to additional convex constraints. We showed in particular how to solve the convex subproblems that arise when adding bound, total-variation (TV), and asymmetric TV constraints to the model space. The proposed framework is general, and the convex constraints can be replaced or augmented by others, as long as the feasibility regions remain nonempty.

Synthetic experiments suggest that for sufficiently accurate starting models,
TV constraints enhance the recovery by improving delineation of the salt and by eliminating spurious artifacts. 
The experiments also show that the asymmetric TV constraint, designed to encourage velocity to increase with depth, leads to major improvements in the recovery of salt structures from poor starting models. 
By pushing the bottom salt down into the basin in the beginning of the inversion procedure, the asymmetric constraint prevents the inversion from creating a velocity-low after the algorithm steps into the salt. In combination with a continuation strategy that gradually weakens the asymmetric TV constraint, we arrive at an approach that avoids getting stuck when starting with a poor initialization. Future work aims to study more realistic numerical experiments and investigate how to better take advantage of the proposed WRI framework.

\section{Acknowledgements}\label{acknowledgements}

Thanks to Bas Peters for insights about the formulation and implementation of the penalty method for full-waveform inversion, including the strategies for frequency continuation and computing multiple passes with warm starts. Thanks also to Polina Zheglova for helpful discussions about the boundary conditions.

\section{In memoriam: John “Ernie” Esser$^\dagger$ (May 19, 1980 -- March 8, 2015)}

This paper is dedicated to its main author, Ernie Esser, who passed away in tragic circumstances during preparation of the manuscript. He was a very promising young scientist, positive, energetic, generous and talented, always a pleasure to work with. He is dearly missed by his family, friends and colleagues.

\begin{figure*}
\centering
\includegraphics[width=0.4\hsize, natwidth = 650 ,natheight=642]{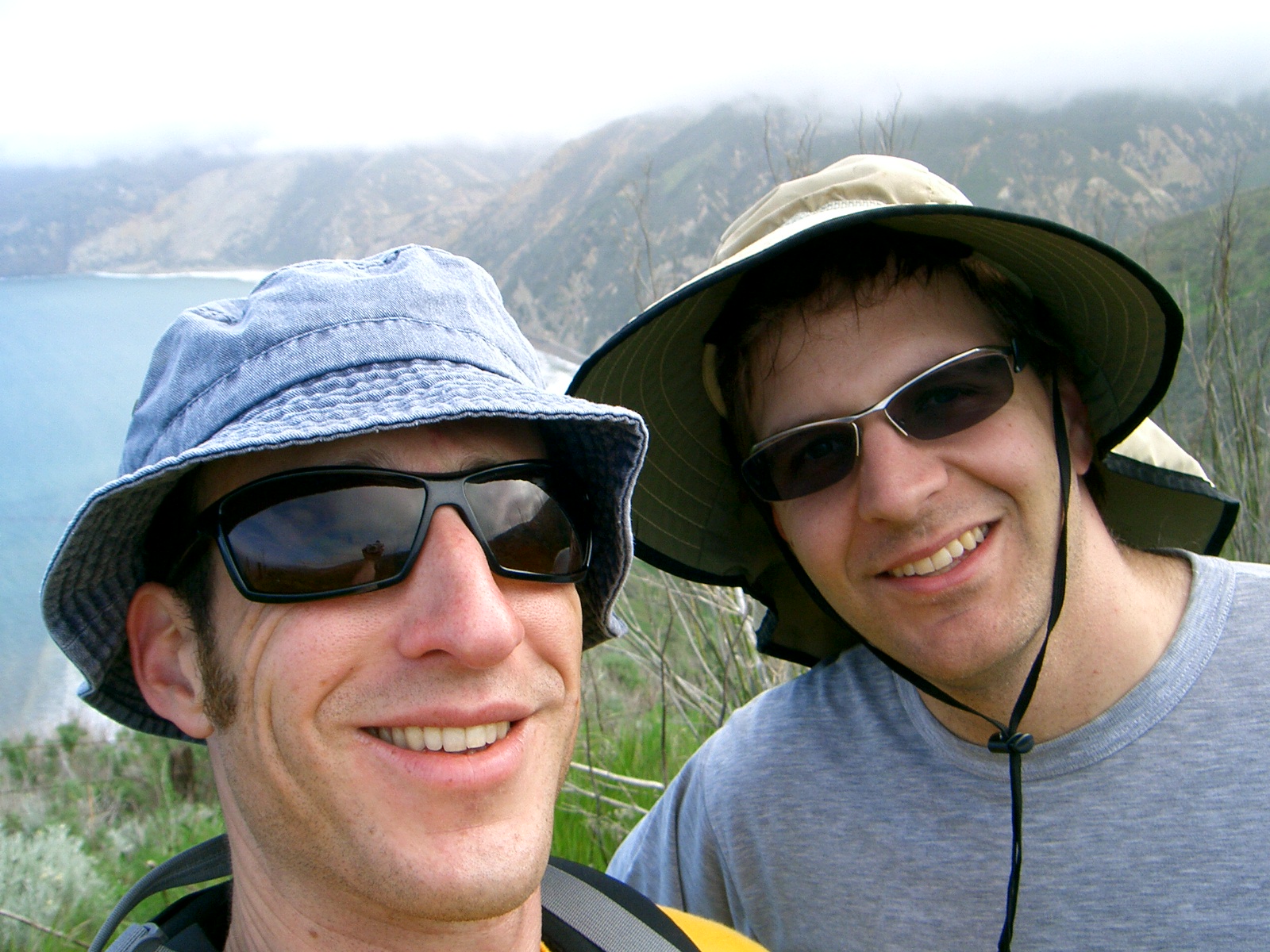}\qquad
\includegraphics[width=0.4\hsize, natwidth = 650 ,natheight=642]{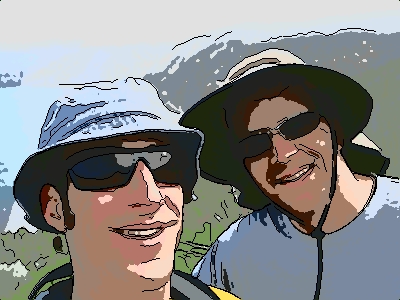}
\caption{Ernie Esser (May 19, 1980 -- March 8, 2015). What Ernie (pictured on the left) loved most in life was hiking with friends, good espresso, and total-variation minimization. The above picture of Ernie with his long-time friend Justin Jacobs was taken on the Channel Islands in 2011. The total-variation constrained version of this picture was made by Ernie.}
\end{figure*}

\bibliography{penaltyTV}


\section{Appendix: PDE-constrained optimization}\label{adjoint-state-formulation}
We develop gradients and Hessian approximations for objectives defined in \eqref{eq:NLLS} and \eqref{eq:ENLLS}. 
These computations also appear in \cite{Haber2000,Virieux2009}.

\subsection*{Gradient and Hessian for the reduced adjoint-state method}

The most commonly used formulation for FWI is obtained by eliminating the
PDE-constraint, in which case the forward modelling operator can be explicitly written as
\[
F(m)q \equiv Pu(m),
\]
with $P$ the sampling operator that models the measurement process. In this formulation, $u(m)$ is the discrete wavefield, computed by solving the discretized PDE 
\[
A(m)u = q.
\]
After solving this system, the objective defined in \eqref{eq:NLLS} becomes
\begin{equation}
\label{eq:FWI}
f(m) = \sum_{j=1}^{N_s} \textstyle{\frac{1}{2}}\|PA(m)^{-1}q_j - d_j\|_2^2. 
\end{equation}
Derivatives with respect to the $i^{\mathrm{th}}$  entry of the model vector $m$ are given by 
\begin{eqnarray}
	\label{eq:gFWI}
\frac{\partial f}{\partial m_i} 
&=& \sum_{j=1}^{N_s}\left(\frac{\partial PA(m)^{-1}q_j}{\partial m_i}\right)^T\left(PA(m)^{-1}q_j - d_j\right) \\
&=& \sum_{j=1}^{N_s} q_j^TA(m)^{-T}\!\left(\frac{\partial A(m)}{\partial m_i}\right)^T\!A(m)^{-T}P^T\!\left(PA(m)^{-1}q_j - d_j\right)
\end{eqnarray}
with the symbol $^T$ denoting the adjoint. 
Introducing intermediate states
\[
u_j = A(m)^{-1}q_j, \quad v_j = A(m)^{-T}P^T\!\left(Pu_j - d_j\right),\quad  j=1\cdots N_s, 
\]
and the matrix
\[
G_i(m) = \frac{\partial A(m)}{\partial m_i},
\]
the derivative may be more succinctly written as
\[
\frac{\partial f}{\partial m_i}  = \sum_{j=1}^{N_s} u_j^TG_i(m)^Tv_j.
\]
The derivative formula requires solving both the forward PDE for each source term to obtain the state variables (wavefields) $u_j$, and the {\it adjoint} PDE to obtain adjoint state variables ("reverse-time" wavefields) $v_j$. 

The full Hessian is dense and its elements are given by
\[
\frac{\partial^2 f}{\partial m_k\partial m_l} = \sum_{j=1}^{N_s} u_j^TR_{kl}^T(m)v_j + u_j^TG_k(m)^TA^{-T}P^T\!PA^{-1}G_l(m)u_j, 
\]
with
\[
R_{kl}(m) = \frac{\partial^2 A(m)}{\partial m_k\partial m_l} + 2G_k^T\!(m)A^{-1}G_l(m).
\]
In practice, the individual elements of the Hessian are never computed for large-scale problems. Instead the action of the Hessian on a given vector can be evaluated at the cost of a few additional PDE-solves. While this approach avoids infeasible explicit storage of the Hessian, the additional PDE solves quickly become too expensive. To avoid these costs, the Hessian is often approximated by the \emph{pseudo-Hessian} approximation, whose elements are given by 
\begin{equation}
	\label{eq:HFWI}
	H_{kl} = \sum_{j=1}^{N_s} u_{j}^TG_k^TG_lu_{j}.
\end{equation}
In our discretization, this approximation turns out to be diagonal, as $G_k$ is a diagonal matrix with a single element at location $(k,k)$. 

\subsection*{Gradient and Hessian for the extended formulation}

While the adjoint-state method undergirds the majority of practical approaches to large-scale FWI, it requires relative accurate starting models to avoid getting stuck in local minima. Extended formulations, where the PDE-constraints are not eliminated but replaced by $\ell_2$-norm penalties, are less prone to these minima because they have more room to fit the data by optimizing over both the model $m$ and state variables $u_j$, $j=1\cdots N_s$. In wave-equation based inversion, the extended formulation \eqref{eq:ENLLS} has explicit form
\begin{equation}
\label{eq:WRI}
f_{\lambda}(m) = \min_{u} \sum_{j=1}^{N_s} \textstyle{\frac{1}{2}}\|Pu_j - d_j\|_2^2 + \textstyle{\frac{\lambda^2}{2}}\|A(m)u_j - q_j\|_2^2,
\end{equation}
where we have made the variable substitution $\Delta q_j \rightarrow A(m)u_j - q_j$. Note 
the objective is written as an optimal value function of $m$. Evaluating $f_{\lambda}$ for a given $m$ requires solving an optimization problem in $u$. However, since both terms are quadratic in the $u_j$'s, a closed-form solution is available:
\[
u_{\lambda,j} = \left(A^T\!A + \lambda^{-2} P^T\!P\right)^{-1}\left(A^Tq_j + \lambda^{-2} P^T\!d_j\right).
\]
This expression can be seen as the equivalent of the PDE-solve needed to evaluate the regular objective but with the important distinction that solutions also aim to fit observed data. Approaches where PDE solves are combined with data-fit objectives are widely known in the field of data assimilation where PDEs appear as $\ell_2$-norm penalties, also known as "weak constraints" \cite{Fisher2005}. Indeed, as $\lambda \uparrow \infty$, we see that the expression reduces to $u_{\lambda,j} = A^{-1}q_j$. The derivatives of this objective are given by 
\begin{eqnarray}
	\label{eq:gWRI}
\frac{\partial f_{\lambda}}{\partial m_k}  = \sum_{j=1}^{N_s} u_{\lambda,j}^TG_k(m)^Tv_{\lambda,j},
\end{eqnarray}
and
\[
\frac{\partial^2 f}{\partial m_k\partial m_l} = \lambda^2 \sum_{j=1}^{N_s} u_{\lambda,j}^TG_k^TG_lu_{\lambda,j} + u_{\lambda,j}^TR_{kl}^Tv_{\lambda,j} - \left(A^TG_ku_{\lambda,j} + G_k^T v_{\lambda,j}\right)^T\!\left(\lambda^2 A^T\!A + P^T\!P\right)^{-1}\!\left(A^TG_lu_{\lambda,j} + G_l^Tv_{\lambda,j}\right)
\]
where 
\[
v_{\lambda,j} = \lambda^2\left(A(m)u_{\lambda,j} - q_j\right).
\]
Moreover, it can be shown that $\nabla f_\lambda$ is Lipschitz continuous, with the bound on the Lipschitz constant independent 
of $\lambda$~\cite{aravkin2016qp}, exactly as required by the implicit trust region framework. By ignoring higher order derivatives and the dependency of $u_{\lambda,i}$ on $m$, we obtain a positive definite approximation of this Hessian with elements
\begin{equation}
\label{eq:HWRI}
\left(H_{\lambda}\right)_{kl} =  \lambda^2 \sum_{i=1}^{N_s} u_{\lambda,i}^TG_k^TG_lu_{\lambda,i}. 
\end{equation}
This approximation is typically sparse and does not involve additional PDE-solves. 

\subsection*{Discretization}
In the numerical experiments we use a finite-difference discretization of the Helmholtz operator with Robin boundary conditions, in which case $A$ is block-diagonal matrix 
\[
A = \left(
\begin{matrix}
A_1&   &      &       \\
   &A_2&      &       \\
   &   &\ddots&       \\ 
   &   &      &A_{N_f}\\
\end{matrix}
\right)
\]
with $N_f$ blocks
\[
A_i = \omega_i^2\mathsf{diag}(b)\mathsf{diag}(m) - \imath\omega_i \mathsf{diag}(1 - b)\mathsf{diag}(m^{1/2}) + L,
\]
where $\omega$ is the angular frequency, $b$ is a vector with $b_i = 1$ in the interior of the domain and $b_i = 0$ on the boundary and $L$ is a 5-point discretization of the Laplace operator with Neumann boundary conditions. Both the wavefields $u$, $v$ and the source vectors $q$ are block-vectors with one block for each frequency as well. A single solve of the system $Au = q$ thus involves solving $N_f$ systems of equations independently. 
The Jacobian matrix $G_k = \frac{\partial A}{\partial m_k}$ is a block-diagonal matrix with blocks
\[
\frac{\partial A_i}{\partial m_k} = \omega_i^2\mathsf{diag}(b)\mathsf{diag}(e_k) - \frac{1}{2}\imath\omega_i \mathsf{diag}(1 - b)\mathsf{diag}(e_k)\mathsf{diag}(m^{-1/2}),
\]
where $e_k$ is the $k^\mathrm{th}$ unit vector.

\end{document}